\documentclass[10pt]{amsart}
      \usepackage[mathscr]{eucal}
      \usepackage{amsmath,amsfonts}

      \parskip\smallskipamount
          \newtheorem{theorem}{Theorem}[section]
      
      \newtheorem{proposition}[theorem]{Proposition}
      \newtheorem{corollary}[theorem]{Corollary}
      \newtheorem{lemma}[theorem]{Lemma}

      \newtheorem{remark}[theorem]{Remark}
      \makeatletter
      \@addtoreset{equation}{section}
      \makeatother

      \newcommand{\CC}{{\mathbb C}}
      \newcommand{\NN}{{\mathbb N}}

      \newcommand{\DD}{{\mathbb D}}
      \newcommand{\RR}{{\mathbb R}}
      \newcommand{\FF}{{\mathbb F}}
      \newcommand{\TT}{{\mathbb T}}

      \newcommand{\cA}{{\mathcal A}}
      \newcommand{\cB}{{\mathcal B}}
      \newcommand{\cC}{{\mathcal C}}
      \newcommand{\cD}{{\mathcal D}}
      \newcommand{\cE}{{\mathcal E}}
      \newcommand{\cF}{{\mathcal F}}
      
      \newcommand{\cH}{{\mathcal H}}
      \newcommand{\cK}{{\mathcal K}}
      \newcommand{\cL}{{\mathcal L}}

      \newcommand{\cP}{{\mathcal P}}
      \newcommand{\cR}{{\mathcal R}}
      \newcommand{\cS}{{\mathcal S}}
      \newcommand{\cT}{{\mathcal T}}

      \newcommand{\cY}{{\mathcal Y}}
      \newcommand{\cX}{{\mathcal X}}

      \newdimen\expt
      \def\boxit#1{\setbox0\hbox{$\displaystyle{#1}$}
            \hbox{\lower.4\expt
       \hbox{\lower3\expt\hbox{\lower\dp0
            \hbox{\vbox{\hrule height.4\expt
       \hbox{\vrule width.4\expt\hskip3\expt
            \vbox{\vskip3\expt\box0\vskip2\expt}%
       \hskip3\expt\vrule width.4\expt}\hrule height.4\expt}}}}}}
      
\begin{document}
       \pagestyle{myheadings}
      \markboth{ Gelu Popescu}{ Noncommutative transforms and free pluriharmonic functions }
      %\pagestyle{plain}
      %\begin{flushright}
       % \it Date of this draft: \today
      %\end{flushright}
      %\bigskip

      \title [Noncommutative transforms and free pluriharmonic functions   ]
      { Noncommutative transforms and free pluriharmonic functions
      }
        \author{Gelu Popescu}
      %\date{\today}
\date{January 12, 2007}
      \thanks{Research supported in part by an NSF grant}

      \address{Department of Mathematics, The University of Texas
      at San Antonio \\ San Antonio, TX 78249, USA}
      \email{\tt gelu.popescu@utsa.edu}

\begin{abstract}
In this paper,  we study free pluriharmonic functions on
noncommutative  balls \ $[B(\cH)^n]_\gamma$, $\gamma>0$,
  and their boundary behavior.   These functions have the form
 $$
 f(X_1,\ldots, X_n)=\sum_{k=1}^\infty \sum_{|\alpha|=k}
 b_\alpha X_\alpha^* + a_0I+ \sum_{k=1}^\infty \sum_{|\alpha|=k}
 a_\alpha X_\alpha,\quad a_\alpha, b_\alpha\in \CC,
 $$
 where the convergence of the series  is in the operator norm topology for any
    $(X_1,\ldots,X_n)\in [B(\cH)^n]_\gamma$,
    and $B(\cH)$ denotes  the algebra
  of all bounded linear operators on a Hilbert space $\cH$.
 The main tools used in this study are certain noncommutative transforms which
  are introduced
 in the present paper and which
 generalize the classical transforms of Berezin, Poisson, Fantappi\` e,
 Herglotz, and
 Cayley. Several classical results from complex analysis have
 free analogues in
 our   noncommutative multivariable setting.

\end{abstract}

      \maketitle

\bigskip

\section*{Contents}
{\it
 Introduction

 \begin{enumerate}
 \item[ 1.]  Multi-Toeplitz operators on Fock spaces and their Fourier representations
 \item[ 2.]  Noncommutative Berezin transforms and free pluriharmonic
 functions
 \item[ 3.]  Bounded free pluriharmonic functions
 \item[ 4.]   Dirichlet extension problem  for free  pluriharmonic functions
\item[ 5.]    Noncommutative  transforms: Fantappi\` e, Herglotz,
and Poisson
\item[ 6.]   The Banach space $Har^1(B(\cH)^n_1)$
 \item[ 7.]    Noncommutative Cayley transforms
 \item[ 8.]   Carath\' eodory interpolation for free holomorphic
 functions with positive real parts
\end{enumerate}

 References
}

  \bigskip

\bigskip

\section*{Introduction}

In recent years, significant progress  has been made in
noncommutative multivariable operator theory regarding
 noncommutative dilation theory, its applications to
    interpolation in several variables,
  and unitary
invariants for $n$-tuples of operators. In \cite{Po-holomorphic}, we
developed a theory of holomorphic functions in several noncommuting
(free) variables and  provide a framework for  the study of
arbitrary
 $n$-tuples of operators. This theory enhances our program to develop a {\it free}
  analogue of
  Sz.-Nagy--Foia\c s theory \cite{SzF-book}, for row contractions.

Let $\FF_n^+$ be the unital free semigroup on $n$ generators
$g_1,\ldots, g_n$ and the identity $g_0$.  The length of $\alpha\in
\FF_n^+$ is defined by $|\alpha|:=0$ if $\alpha=g_0$  and
$|\alpha|:=k$ if
 $\alpha=g_{i_1}\cdots g_{i_k}$, where $i_1,\ldots, i_k\in \{1,\ldots, n\}$.
If $(X_1,\ldots, X_n)\in B(\cH)^n$, where $B(\cH)$ is the algebra of
all bounded linear operators on the Hilbert space $\cH$,    we
denote $X_\alpha:= X_{i_1}\cdots X_{i_k}$  and $X_{g_0}:=I_\cH$.
Free pluriharmonic functions arise in the study of free holomorphic functions
on the
noncommutative  open  unit ball
$$
[B(\cH)^n]_1:=\{(X_1,\ldots, X_n)\in B(\cH)^n: \ \|X_1X_1^*+\cdots +
X_nX_n^*\|^{1/2}<1\}.
$$
We recall that a map $f:[B(\cH)^n]_1\to B(\cH)$ is called free
holomorphic function  with scalar coefficients if \
$\limsup\limits_{k\to\infty} \left(\sum\limits_{|\alpha|=k}
|a_\alpha|^2\right)^{1/2k}\leq 1$ and $ f(X_1,\ldots,
X_n)=\sum_{k=0}^\infty \sum_{|\alpha|=k}
 a_\alpha X_\alpha$,\quad $(X_1,\ldots, X_n)\in [B(\cH)^n]_1.
$
We say that $h:[B(\cH)^n]_1\to B(\cH)$ is a self-adjoint free
pluriharmonic function on $[B(\cH)^n]_1$ if $h=\text{\rm Re}\, f$
for some free holomorphic function $f$. An arbitrary free
pluriharmonic function is a linear combination of self-adjoint free
pluriharmonic functions.

In this paper,  we study free pluriharmonic functions on the noncommutative ball
 $[B(\cH)^n]_1$ and their boundary behavior.
 The main tools used in this study are   noncommutative transforms which
 generalize the classical transforms of Berezin, Poisson, Fantappi\` e,
 Herglotz-Riesz, and
 Cayley (see \cite{Be}, \cite{H},  \cite{Her}, \cite{Co}, \cite{Ri},  \cite{R}, \cite{R2}).
 We show that several classical results from complex analysis have
 free analogues in
 our   noncommutative multivariable setting.

Multi-Toeplitz operators on the full Fock space on $n$ generators
$F^2(H_n)$ have played an important role in multivariable  operator
theory (\cite{Po-multi}, \cite{Po-analytic},  \cite{Po-structure},
\cite{Po-nehari}, \cite{Po-entropy}, \cite{Po-unitary}, \cite{DKP}).
In Section 1, we associate with each  multi-Toeplitz operator
 a formal Fourier series
$$
\varphi(S_1,\ldots, S_n)=\sum_{|\alpha|\geq 1} b_\alpha  S_\alpha^* +a_0I
+\sum_{|\alpha|\geq 1} a_\alpha S_\alpha,
$$
where
$S_1,\ldots, S_n$  are the left creation operators on $F^2(H_n)$.
We show that a multi-Toeplitz operator is determined by its Fourier series
and can be recaptured
from it. The main result of Section 1  is a
 characterization of  the multi-Toeplitz operators
  in terms of  their Fourier representations.
As a consequence, we deduce  that the set of all multi-Toeplitz
operators coincides with $ \overline{\cA_n^*
+\cA_n}^{SOT}=\overline{\cA_n^* +\cA_n}^{WOT}, $ where $\cA_n$ is
the noncommutative disc algebra (\cite{Po-von}, \cite{Po-disc}),
i.e., the norm
 closed algebra generated by $S_1,\ldots, S_n$ and the identity.

Let $Har(B(\cH)^n_1)$ be the set of all free pluriharmonic functions
 on $[B(\cH)^n]_1$ with operator-valued coefficients. When the coefficients
  are scalars, we
 use the notation $Har_\CC(B(\cH)^n_1)$.
 An important role in the study of the  free pluriharmonic
functions  and their boundary behavior  is played by the
noncommutative Berezin transforms  $\cB_\mu$, introduced in
 Section 2,  which are  associated with
completely bounded maps  $\mu$  on $B(F^2(H_n))$.

Throughout this paper,  the Berezin transform $B_\tau$, where $\tau$ is the linear
functional on $B(F^2(H_n))$ defined by $\tau(f):=\left<f(1),1\right>$, will be
called Poisson transform because it coincides with the noncommutative
 Poisson transform introduced
in \cite{Po-poisson}.   If $f\in B(F^2(H_n))$ and $X\in
[B(\cH)^n]_1$, then the Poisson transform of $f$ at $X$ satisfies the equations
$$
P_X[f]= B_\tau(f,X)=K_X^*(f\otimes I_\cH) K_X,
$$
where $K_X$ is the noncommutative Poisson kernel.

The classical characterization of the harmonic functions on the open
unit disc $\DD$ as continuous functions with the mean value property
has a noncommutative analogue  in our setting.  We show that a free
pluriharmonic function $u:[B(\cH)^n]_1\to B(\cH)$ is uniquely
determined  by its radial function
$$[0,1)\ni r\mapsto u(rS_1,\ldots, rS_n)\in \cA_n^* +\cA_n
$$
and  the {\it Poisson mean value property}, i.e., $ u(X_1,\ldots,
X_n)=P_{\frac{1}{r} X} [u(rS_1,\ldots, rS_n)] $ for $X:=(X_1,\ldots,
X_n)\in [B(\cH)^n]_1$ and $r\in (0,1)$. This characterization is
used to obtain a Weierstrass  type convergence theorem for free
pluriharmonic
 functions,
 which enables us to introduce a metric on $Har_\CC(B(\cH)^n_1)$ with respect
  to which
 it becomes a complete metric space.

We prove a Harnack type inequality  (see \cite{Co}, \cite{R} for the
classical result)  for positive free pluriharmonic functions and
obtain a Harnack type convergence theorem for increasing
 sequences of  free pluriharmonic functions, as well as  a maximum
 (resp. minimum)
 principle for free pluriharmonic functions.

In Section 3, we characterize the set $Har_\CC^\infty(B(\cH)^n_1)$ of all bounded
 free pluriharmonic functions on $[B(\cH)^n]_1$ in terms of the boundary functions in
 $ \overline{\cA_n^* +\cA_n}^{SOT}$   and   obtain a Fatou type
result \cite{H} for bounded free
 pluriharmonic functions, which extends the $F_n^\infty$-functional calculus
 for pure row contractions \cite{Po-funct}.

The Dirichlet problem ( \cite{H}, \cite{Co}) for the unit disc $\DD$
states: given a continuous function $f$ on the unit circle
$\TT:=\{z\in \CC:\ |z|=1\}$, find a continuous function $h$
 on
$\overline{\DD}$ such that $h|_\TT=f$ and $h|_\DD$ is harmonic. This
problem is completely solved by the Poisson integral formula. In
Section 4, we consider an analogue of this problem for free
pluriharmonic functions.  We prove that a function
$u:[B(\cH)^n]_1\to B(\cH)$ is free pluriharmonic and has continuous
extension (in the operator norm topology) to the closed ball
$[B(\cH)^n]_1^-$ if and only if there exists $f\in \overline{\cA_n^*
+\cA_n}^{\|\cdot\|}$  such that
$$
 u(X_1,\ldots, X_n)=P_X[f],\qquad X:=(X_1,\ldots, X_n)\in [B(\cH)^n]_1.
 $$
   A similar result is provided  for the
class of  {\it $C^*$-harmonic functions} $u:[B(\cH)^n]_1\to B(\cH)$
 which have
  continuous
extensions (in the operator norm topology) to the closed ball
$[B(\cH)^n]_1^-$.   A version of the maximum  principle for
$C^*$-harmonic functions is also obtained.

In Section 5, we introduce noncommutative versions of Fantappi\` e, Herglotz,
and Poisson
transforms associated with completely bounded maps on the operator system
$\cR_n^* + \cR_n$ (or $B(F^2(H_n))$), where $\cR_n$ is the noncommutative
 disc algebra
generated by the right creation operators $R_1,\ldots, R_n$ on
 $F^2(H_n)$ and the identity.  These transforms are used to obtain  characterizations
for the set of all free holomorphic functions
on $[B(\cH)^n]_1$
with positive real parts, and to study the geometric structure and boundary
 behavior of the free pluriharmonic functions on $[B(\cH)^n]_1$.

In particular, we obtain the following noncommutative analogue of the
 Herglotz-Riesz
representation theorem (\cite{Her}, \cite{Ri}): if  $f:[B(\cH)^n]_1
\to B(\cH)$ is a free holomorphic function with $\text{\rm Re}\,
f\geq 0$ on $[B(\cH)^n]_1$, then there is a positive linear  map
$\mu$ on the Cuntz-Toeplitz algebra $C^*(R_1,\ldots, R_n)$ such that
$$
f(X_1,\ldots, X_n)=(H\mu)(X_1,\ldots, X_n) +i(\text{\rm Im}\,  f(0)),
$$
where the noncommutative Herglotz transform $H\mu: [B(\cH)^n]_1\to B(\cH)$ is
 defined
by
$$
(H\mu)(X_1,\ldots, X_n):=(\mu\otimes \text{\rm id})\left[
2(I-R_1^*\otimes X_1-\cdots -R_n^*\otimes X_n)^{-1}-I\right]
$$
and $(\mu\otimes \text{\rm id})(f\otimes Y):=\mu(f)Y$ for $f\in
C^*(R_1,\ldots, R_n)$ and $Y\in B(\cH)$.

In Section 5, we also introduce the noncommutative Poisson transform
of a completely
bounded linear map on $B(F^2(H_n))$ and show that it is a particular case of the Berezin transform of
Section 2.
In the particular case  when $\mu$ is a bounded linear functional on
 $C^*(R_1,\ldots, R_n)$, the Poisson transform
 $\cP\mu:[B(\cH)^n]_1\to B(\cH)$ is defined by
$$
(\cP\mu)(X_1,\ldots, X_n):=(\mu\otimes \text{\rm
id})\left[P(R,X)\right], \qquad X:=(X_1,\ldots,X_n)\in [B(\cH)^n]_1,
$$
where the free pluriharmonic Poisson kernel is given by
$$
P(R,X):=\sum_{k=1}^\infty\sum_{|\alpha|=k}
R_{\widetilde\alpha}\otimes X_\alpha^*
+I+\sum_{k=1}^\infty\sum_{|\alpha|=k} R_{\widetilde\alpha}^*\otimes
X_\alpha
$$
and the series are convergent in the operator norm topology.

We show that the map $\mu\mapsto \cP\mu$ is a  linear and one-to-one
correspondence between the space of all completely positive linear
maps on the operator system $\cR_n^* + \cR_n$ and the space of all
positive free pluriharmonic functions on the open noncommutative
ball $[B(\cH)^n]_1$ with operator-valued coefficients. In
particular, any positive free pluriharmonic function on
$[B(\cH)^n]_1$  is the Poisson transform of   a completely positive
linear  map  on the Cuntz-Toeplitz algebra $C^*(R_1,\ldots, R_n)$.
Moreover, we show that a  free pluriharmonic function $h:
[B(\cH)^n]_1\to B(\cH)$ is positive if and only is there exists an
$n$-tuple of isometries $(V_1,\ldots, V_n)$
on a Hilbert space $\cK$, with orthogonal ranges,  and a vector $\xi\in \cK$ such that
$$
h(X_1,\ldots, X_n)= (\omega_{\xi}\otimes  \text{\rm id})\left[
B_X(V_1,\ldots, V_n)^* B_X(V_1,\ldots, V_n) \right],
$$
where $B_X(V_1,\ldots, V_n)$ is the noncommutative  Berezin kernel
 defined in Section 2 and $\omega_{\xi}$ is the linear functional
 defined by $\omega_{\xi}(Y):=\left<Y\xi,\xi\right>$.

In Section 6,  we introduce  the space $Har_\CC^1(B(\cH)^n_1)$ of all
 free pluriharmonic functions
$h$ such that $\|h\|_1:=\sup\limits_{0\leq
r<1}\|\nu_{h,r}\|<\infty$, where $\{\nu_{h,r}\}$ are bounded linear
functionals associated with the radial function $ [0,1)\ni r\mapsto
h(rR_1,\ldots, rR_n)\in \cR_n^*+\cR_n. $ We show that \
$(Har_\CC^1(B(\cH)^n_1), \|\cdot\|_1)$ \  is a Banach space that can
be identified, through  the noncommutative Poisson transform, with
the dual of
 the operator
system $\cR_n^*+\cR_n$. As a consequence, we characterize the self-adjoint free
 pluriharmonic functions $u$ which admit a Jordan type decomposition $u=u^+-u^-$,
  where
 $u^+$, $u^-$ are positive free pluriharmonic functions on $[B(\cH)^n]_1$.
Another consequence of the above-mentioned result is that the space of free holomorphic functions
$$H_\CC^1(B(\cH)^n_1):=Hol_\CC(B(\cH)^n_1)\cap Har_\CC^1(B(\cH)^n_1)$$
is a Banach space (with respect to $\|\cdot\|_1$) which can be
identified with the annihilator of $\cR_n$ in the dual of the
operator system $\cR_n^*+\cR_n$.

  In Section 7, we introduce a noncommutative Cayley transform which turns out
  to be a bijection between the set of all contractive free holomorphic functions
  $f$  on
  $[B(\cH)^n]_1$ with $f(0)=0$, and  the set of all free holomorphic functions $g$
  with $g(0)=0$ and
  $$g(X_1,\ldots, X_n)^*+I+g(X_1,\ldots,X_n)\geq 0 \quad \text{ for any } \
  (X_1,\ldots, X_n)\in[B(\cH)^n]_1.
  $$
  This result and its consequences concerning truncated Cayley transforms
  are used,
  in Section 8, to solve the Carath\' eodory interpolation problem for
  free holomorphic
  functions with positive real parts on $[B(\cH)^n]_1$.
  We show that given a sequence of complex numbers
   $\{b_\alpha\}_{|\alpha|\leq m}$ with $b_0\geq 0$, there exists a sequence
   $\{b_\alpha\}_{|\alpha|\geq m+1}\subset \CC$ such that
   $$
   g(X_1,\ldots, X_n):=\frac{b_0}{2} +\sum_{k=1}^\infty \sum_{|\alpha|=k}
    b_\alpha X_\alpha, \quad
  (X_1,\ldots, X_n)\in[B(\cH)^n]_1,
   $$
   is a free holomorphic function with Re\,$g(X_1,\ldots, X_n)\geq 0$
   for any $(X_1,\ldots, X_n)\in [B(\cH)^n]_1$
   if and only if
   $$
   \sum_{1\leq |\alpha|\leq m} \bar b_\alpha (S_\alpha^{(m)})^*+ b_0 I+
   \sum_{1\leq |\alpha|\leq m}  b_\alpha S_\alpha^{(m)}\geq 0,
   $$
where $S_1^{(m)},\ldots, S_m^{(m)}$ are the compressions
 of the left creation operators $S_1,\ldots, S_n$ to the subspace $\cP^{(m)}$ of all
 polynomials in $F^2(H_n)$ of degree $\leq m$. We also show that the condition above is
 equivalent to the existence of a positive linear map $\nu$ on
 $C^*(S_1,\ldots, S_n)$ such that
 $$\nu(S_\alpha)=\bar b_\alpha\quad   \text{ for } \ |\alpha|\leq m,
 $$
 i.e., $\nu$ solves the noncommutative trigonometric moment problem
 for the operator system $\cA_n^*+ \cA_n$, with data
  $\{\bar b_\alpha\}_{|\alpha|\leq m}$.

  We also show that  the  Carath\' eodory  interpolation
problem for free holomorphic functions  with positive real parts on
$[B(\cH)^n]_1$ is equivalent to  the Carath\' eodory-Fej\' er
interpolation  problem  for multi-analytic operators
\cite{Po-analytic} and to the Carath\' eodory  interpolation problem
for positive semidefinite multi-Toeplitz kernels on free semigroups
\cite{Po-structure} (see \cite{Ca}, \cite{CaFe}, \cite{Sc}, \cite{S}
for the classical results). This result together with
\cite{Po-structure} provide a
 parametrization of all solutions of the  Carath\' eodory interpolation problem
 for free holomorphic functions  with
positive real parts, in terms  of generalized Schur sequences.

Finally, we should mention that all the results of this paper are
presented in the more general setting of free pluriharmonic
functions
 with operator-valued coefficients.

\section{ Multi-Toeplitz operators on Fock spaces and their Fourier representations}

There are three fundamental questions about multi-Toeplitz operators on Fock spaces
  and the associated Fourier
series.
\begin{enumerate}
\item[(1)] Is a multi-Toeplitz operator $A$ determined by its Fourier series $?$
\item[(2)] If so, how can we recapture $A$, given the Fourier series
$?$
\item[(3)] Given
 $\{A_{(\alpha)}\}_{\alpha\in \FF_n^+}$ and
 $\{B_{(\alpha)}\}_{\alpha\in \FF_n^+\backslash\{g_0\}}$,
 two sequences of operators on a Hilbert space $\cE$, when is the   formal series
 associated with them the formal  Fourier representation of a multi-Toeplitz
 operator  on $\cE\otimes F^2(H_n)$ $?$
\end{enumerate}
We will answer these questions  in this section. The results will
play an important role in our investigation.

Let $H_n$ be an $n$-dimensional complex  Hilbert space with
orthonormal
      basis
      $e_1$, $e_2$, $\dots,e_n$, where $n\in\{1,2,\dots\}$.
       We consider the full Fock space  of $H_n$ defined by
      $$F^2(H_n):=\bigoplus_{k\geq 0} H_n^{\otimes k},$$
      where $H_n^{\otimes 0}:=\CC 1$ and $H_n^{\otimes k}$ is the (Hilbert)
      tensor product of $k$ copies of $H_n$.
      Define the left  (resp. right) creation
      operators  $S_i$ (resp.~$R_i$), $i=1,\ldots,n$, acting on $F^2(H_n)$  by
      setting
      $
       S_i\varphi:=e_i\otimes\varphi, \quad  \varphi\in F^2(H_n),
      $
       (resp.~$
       R_i\varphi:=\varphi\otimes e_i, \quad  \varphi\in F^2(H_n).
      $)
The noncommutative disc algebra $\cA_n$ (resp.~$\cR_n$)
is the norm closed algebra
generated by the left (resp.~right) creation operators and the identity.
The   noncommutative analytic Toeplitz algebra $F_n^\infty$ (resp.~$\cR_n^\infty$)
 is the the weakly
closed version of $\cA_n$ (resp.~$\cR_n$). These algebras were
introduced in \cite{Po-von} in connection with a noncommutative von
Neumann inequality (see \cite{vN} for the classical case). They
 have  been studied
    in several papers
   \cite{Po-multi},  \cite{Po-funct}, \cite{Po-analytic},
\cite{Po-disc}, \cite{Po-poisson},
  \cite{DP1}, \cite{DP2},
    \cite{Po-curvature},  \cite{DKP},  \cite{PPoS},
       and \cite{Po-holomorphic}.

Let $\FF_n^+$ be the unital free semigroup on $n$ generators
      $g_1,\dots,g_n$, and the identity $g_0$.
       We denote $e_\alpha:=
e_{i_1}\otimes\cdots \otimes  e_{i_k}$ and $e_{g_0}:=1$. Note that
$\{e_\alpha\}_{\alpha\in \FF_n^+}$ is an orthonormal basis for
$F^2(H_n)$. An operator $A\in B(\cE\otimes F^2(H_n))$ is called
multi-Toeplitz with  respect to the right creation operators
$R_1,\ldots, R_n$ if and only if  $$(I_\cE\otimes
R_i^*)A(I_\cE\otimes R_j)=\delta_{ij} A\quad \text{  for } \
i,j=1,\ldots, n.
$$
When $n=1$  and $\cE=\CC$ we find again the classical Toeplitz operators on the
Hardy space $H^2(\DD)$. Define the formal  Fourier representation of
$A$  by setting
$$
\varphi(S_1,\ldots, S_n):= \sum_{|\alpha|\geq 1} B_{(\alpha)}\otimes  S_\alpha^*
+ A_{(0)}\otimes  I +\sum_{|\alpha|\geq 1} A_{(\alpha)}\otimes S_\alpha,
$$
where the coefficients are given by
\begin{equation}\label{f-coef}
\begin{split}
\left<A_{(\alpha)}x,y\right> &:=\left< A(x\otimes 1), y\otimes
e_\alpha\right>, \quad \alpha\in \FF_n^+,\\
 \left<B_{(\alpha)}x,
y\right>&:=\left< A(x\otimes e_\alpha), y\otimes 1\right>, \quad
\alpha\in \FF_n^+\backslash\{g_0\},
\end{split}
\end{equation}
for any $x,y\in \cE$. We also set $A_{(0)}:=A_{(g_0)}$.

A few more notations are necessary. If $\omega, \gamma\in \FF_n^+$,
we say that $\omega
>_{r}\gamma$ if there is $\sigma\in
\FF_n^+\backslash\{g_0\}$ such that $\omega=\sigma \gamma$. In this
case  we set $\omega\backslash_r \gamma:=\sigma$. Similarly, we say
that $\omega
>_{l}\gamma$ if there is $\sigma\in
\FF_n^+\backslash\{g_0\}$ such that $\omega= \gamma \sigma$ and set
$\omega\backslash_l \gamma:=\sigma$.
 We denote by
$\tilde\alpha$  the reverse of $\alpha\in \FF_n^+$, i.e.,
  $\tilde \alpha= g_{i_k}\cdots g_{i_k}$ if
   $\alpha=g_{i_1}\cdots g_{i_k}\in\FF_n^+$.
  Notice that $\omega>_{r}\gamma$ if and only if $\tilde \omega
  >_l\tilde \gamma$. In this case we have
 $ \widetilde{\omega\backslash_r\gamma}=\tilde\omega \backslash_l
  \tilde\gamma$.

\begin{theorem}\label{Fourier}
If $A\in B(\cE\otimes F^2(H_n))$ is a multi-Toeplitz operator and
$\varphi(S_1,\ldots, S_n)$ is its formal Fourier  representation,
then $Aq=\varphi(S_1,\ldots, S_n)q$ for any vector-valued polynomial
$q=\sum_{|\alpha|\leq m} h_\alpha\otimes e_\alpha$,
$h_\alpha\in \cE$, and $m\in \NN$. If $A$,  $B$ are multi-Toeplitz operators  having the
same formal Fourier  representation, then $A=B$.
\end{theorem}

\begin{proof}

 Notice that, since
$A(x\otimes 1)=A_{(0)}x\otimes 1+\sum_{|\alpha|\geq 1} (A_{(\alpha)}x\otimes
 e_\alpha) \in \cE\otimes F^2(H_n)$, we
deduce that the series \  $\sum_{|\alpha|\geq 1}
A_{(\alpha)}^*A_{(\alpha)}$ is convergent in the weak operator
 topology (WOT). Similarly,
since we have $A^*(x\otimes 1)=A^*_{(0)}x\otimes
1+\sum_{|\alpha|\geq 1} (B^*_{(\alpha)}x\otimes
 e_\alpha) \in \cE\otimes F^2(H_n)$,
  we deduce that
$\sum_{|\alpha|\geq 1} B_{(\alpha)}B_{(\alpha)}^*$ is WOT convergent. This implies that
$$\varphi(S_1,\ldots, S_n)q:= \sum_{|\alpha|\geq 1}
(B_{(\alpha)}\otimes  S_\alpha^*)q + (A_{(0)}\otimes I) q +\sum_{|\alpha|\geq 1}
(A_{(\alpha)}\otimes
S_\alpha) q$$
 makes sense as a vector in the  Hilbert space tensor product
$\cE\otimes F^2(H_n)$.
Since  $A$ is a multi-Toeplitz  operator, we deduce that
\begin{equation*}
\begin{split}
\left< A(x\otimes e_\gamma), y\otimes e_\omega\right>&= \left<
(I\otimes R_{\tilde\omega}^*) A(I\otimes R_{\tilde \gamma})(x\otimes
1),y\otimes 1\right>
=
\begin{cases}
\left< (I\otimes R^*_{\tilde \omega \backslash_l
\tilde\gamma})A(x\otimes 1,y\otimes1\right>;&\quad
\tilde\omega>_l\tilde\gamma\\
\left<A(x\otimes 1),y\otimes 1\right>;&\quad  \tilde\omega=\tilde\gamma\\
\left<A(I\otimes R_{\tilde\gamma\backslash_l\tilde\omega})
(x\otimes 1,y\otimes 1\right>;&\quad
\tilde\gamma>_l \tilde\omega\\
0;&\quad \text{otherwise}
\end{cases}\\
& =
\begin{cases}
\left<  A(x\otimes 1),y\otimes e_{\omega \backslash_r \gamma}\right>;&\quad
\omega>_r \gamma\\
\left<A(x\otimes 1),y\otimes 1\right>;&\quad  \omega=\gamma\\
\left<A (x\otimes e_{\gamma\backslash_r\omega}),y\otimes 1\right>;&\quad
\gamma>_r \omega\\
0;&\quad \text{otherwise}
\end{cases}
=
\begin{cases}
 \left<A_{(\omega \backslash_r \gamma)}x,y\right>;&\quad
\omega>_r \gamma\\
\left<A_{(0)}x,y\right>;&\quad  \omega=\gamma\\
 \left<B_{(\gamma\backslash_r\omega)}x,y\right>;&\quad
\gamma>_r \omega\\
0;&\quad \text{otherwise}
\end{cases}
\end{split}
\end{equation*}
for any $x,y\in \cE$ and  $\gamma,\omega\in \FF_n^+$.
On the other hand, since   $S_j^* S_i=\delta_{ij} I$ for
$i,j=1,\ldots,n$,  and $\{e_\alpha\}_{\alpha\in \FF_n^+}$ is an
orthonormal basis for $F^2(H_n)$, we  have
\begin{equation*}
\begin{split}
&\left<\varphi(S_1,\ldots, S_n)(x\otimes e_\gamma),y\otimes e_\omega\right>=
\left<\varphi(S_1,\ldots, S_n) (I\otimes S_\gamma)(x\otimes 1), (I\otimes
S_\omega)(y\otimes 1)\right>\\
 &= \left<(I\otimes S_\omega ^*)\left( \sum_{|\alpha|\geq 1} B_{(\alpha)}
(I\otimes S_\alpha^*)
+ A_{(0)}\otimes  I +\sum_{|\alpha|\geq 1} A_{(\alpha)}\otimes
 S_\alpha \right)(I\otimes  S_\gamma)(x\otimes 1),
y\otimes 1\right>\\
&=\left< \sum_{|\alpha|\geq 1}(B_{(\alpha)}\otimes S_\omega ^*
S_\alpha^* S_\gamma) (x\otimes 1),y\otimes 1\right>
+\left<A_{(0)}x,y\right>\left< S_\omega^* S_\gamma 1,1\right>
 + \left< \sum_{|\alpha|\geq 1}(A_{(\alpha)}\otimes
S_\omega ^* S_\alpha S_\gamma)
(x\otimes 1),y\otimes 1\right>\\
&=
\begin{cases} \left<A_{(\omega\backslash_r \gamma)}x,y\right>;&\quad  \omega
>_r\gamma\\
\left<A_{(0)}x,y\right> ;&\quad  \omega=\gamma\\
 \left<B_{(\gamma\backslash_r\omega)}x,y\right>;&\quad
\gamma>_r \omega\\
0;&\quad \text{otherwise}.
\end{cases}
\end{split}
\end{equation*}
Therefore, $$ \left< A(x\otimes e_\gamma), y\otimes e_\omega\right>=
\left<\varphi(S_1,\ldots, S_n)(x\otimes e_\gamma), y\otimes e_\omega\right>
$$
 for any $x,y\in \cE$ and  $\gamma,\omega\in \FF_n^+$.
 Hence, we deduce that  $Aq=\varphi(S_1,\ldots, S_n)q$ for
  any vector-valued polynomial
$q=\sum_{|\alpha|\leq m} h_\alpha\otimes e_\alpha$,
$h_\alpha\in \cE$ and $m\in \NN$.
   The last part of the theorem follows now   easily. The
proof is complete.
\end{proof}

 It is easy to see that if $A$ is a multi-Toeplitz operator, then
 $A=A^*$ if and only if $A_{(0)}=A_{(0)}^*$ and
 $B_{(\alpha)}=A_{(\alpha)}^*$ for any $\alpha\in \FF_n^+\backslash\{g_0\}$.

\bigskip

An $n$-tuple  $T:=(T_1,\dots, T_n)$ of bounded linear  operators
acting on a common Hilbert space $\cH$
   is called
contractive (or row contraction) if $$ T_1T_1^*+\cdots +T_nT_n^*\leq
I_\cH. $$
The defect operators associated with
$T$  are
$$
D_{T^*}:=\left(I_\cH-\sum_{i=1}^n T_iT_i^*\right)^{1/2}\in
B(\cH)\quad \text{ and }\quad
D_{T}:=\left(\left[\delta_{ij}I_\cH-T_i^*T_j\right]_{n\times
n}\right)^{1/2}\in B(\cH^{(n)}),
$$
while the defect spaces  of $T$ are
$\cD_*=\cD_{T^*}:=\overline{D_{T^*}\cH}$ and
$\cD=\cD_{T}:=\overline{D_{T}\cH^{(n)}}$, where
$\cH^{(n)}:=\oplus_{i=1}^n \cH$ denotes the direct sum of $n$ copies
of $\cH$.
We say that an $n$-tuple  $V:=(V_1,\dots, V_n)$ of isometries on a
Hilbert space $\cK\supset \cH$ is a  minimal isometric dilation of
$T$ if the following properties are satisfied:
\begin{enumerate}
\item[(i)]  $V_1V_1^*+\cdots +V_nV_n^*\le I_\cK;$
\item[(ii)]  $V_i^*|_\cH=T_i^*, \ i=1,\dots,n;$
\item[(iii)] $\cK=\bigvee_{\alpha\in \FF^+_n} V_\alpha \cH.$
\end{enumerate}
The isometric dilation theorem for row contractions (see \cite{Bu},
\cite{F}, \cite{Po-isometric})
 asserts that
  every  row contraction $T$ has a minimal isometric
  dilation $V$, which is uniquely
determined up to an isomorphism.
 Let $\Delta_i:\cH\to F^2(H_n)\otimes \cD$ be defined by
$$
\Delta_i h:= 1\otimes D_{T}(
 \underbrace{0,\ldots, 0,}_{\text{$i-1$}\ times}
 h, 0,\dots,0)\oplus 0\oplus 0\cdots.
$$
Consider the Hilbert space $\cK:=\cH\oplus (F^2(H_n)\otimes \cD)$
and embed $\cH$ and  $\cD$ in $\cK$ in the natural way. For each
$i=1,\ldots, n$,  define the operator  $V_i:\cK\to\cK$ by
\begin{equation}
\label{iso-dil}
 V_i(h\oplus (\xi\otimes d)):= T_ih \oplus [\Delta_i
h +(S_i\otimes I_{\cD})(\xi\otimes d)]
\end{equation}
for any $h\in \cH, \xi\in F^2(H_n), d\in \cD$, where $S_1,\ldots,
S_n$ are the left creation operators on the full Fock space
$F^2(H_n)$. The $n$-tuple  $V:=(V_1,\ldots, V_n)$,  is a realization
of the minimal isometric dilation of $T$. According to
\cite{Po-isometric},
\begin{equation}
\label{L*}
\cL_*:=\overline{\left(I_\cK-\sum_{i=1}^n
V_iT_i^*\right)\cH}
\end{equation}
 is wandering subspace  for $V$, i.e., $V_\alpha\cL_*\perp V_\beta \cL_*$ for
 any $\alpha, \beta\in \FF_n^+$ with $\alpha\neq \beta$. Moreover, there
 is a unitary operator $\Phi_*:\cL_*\to \cD_*$ defined by
 \begin{equation}\label{fistea}\Phi_*\left(I-\sum_{j=1}^n
V_jT^*_j\right)h=D_{T^*}h,\quad h\in\cH.
\end{equation}
We recall that $\cK=M_V(\cL_*):=\oplus_{\alpha\in \FF_n^+} V_\alpha
\cL_*$ if and only if $T$ is a pure row contraction, i.e.,
$\sum\limits_{|\alpha|=k}\|T_\alpha^*h\|^2\to 0$ as $k\to\infty$, for any
$h\in \cH$.

We denote by $\cA_n(\cE)$ the spatial tensor product
$B(\cE)\otimes_{min}\cA_n$, where $\cA_n$ is the noncommutative disc
algebra. The main result of this section is the following
characterization of the  multi-Toeplitz operators in terms of their
Fourier representations.

\begin{theorem}\label{Toeplitz}
Let $\{A_{(\alpha)}\}_{\alpha\in \FF_n^+}$ and $\{B_{(\alpha)}\}_{\alpha\in
\FF_n^+\backslash \{g_0\}}$ be two sequences of  operators on a Hilbert space $\cE$.
Then
$$
\varphi(S_1,\ldots, S_n):=\sum_{|\alpha|\geq 1} B_{(\alpha)}\otimes S_\alpha^*
+ A_{(0)}\otimes  I +\sum_{|\alpha|\geq 1} A_{(\alpha)}\otimes  S_\alpha
$$
is the Fourier representation of a multi-Toeplitz operator $A\in
B(\cE\otimes F^2(H_n))$ if and only if
\begin{enumerate}
\item[(i)]
$\sum_{|\alpha|\geq 1} A_{(\alpha)}^* A_{(\alpha)}$ and
$\sum_{|\alpha|\geq 1} B_{(\alpha)} B_{(\alpha)}^*$ are WOT
convergent series, and
\item[(ii)]
$\sup\limits_{0\leq r<1} \|\varphi(rS_1,\ldots, rS_n)\|<\infty$.
\end{enumerate}
Moreover, in this case,
\begin{enumerate}
\item[(a)]  for each $r\in [0,1)$, the operator
$$\varphi(rS_1,\ldots, rS_n):=\sum_{k=1}^\infty \sum_{|\alpha|=k}
B_{(\alpha)}\otimes  r^{|\alpha|} S_\alpha^* + A_{(0)}\otimes I
 +\sum_{k=1}^\infty
\sum_{|\alpha|=k} A_{(\alpha)}\otimes r^{|\alpha|} S_\alpha $$ is in the operator space
$\cA_n(\cE)^* +\cA_n(\cE) $, where the series are convergent in the operator
norm topology;
\item[(b)]
$A=\text{\rm SOT-}\lim\limits_{r\to 1} \varphi(rS_1,\ldots, rS_n)$,
and
\item[(c)]
$ \|A\|=\sup\limits_{0\leq r<1} \|\varphi(rS_1,\ldots,
rS_n)\|=\lim\limits_{r\to 1} \|\varphi(rS_1,\ldots,
rS_n)\|=\sup\limits_{q\in \cE\otimes\cP, \|q\|\leq 1}
\|\varphi(S_1,\ldots, S_n)q\|.$
\end{enumerate}
\end{theorem}

\begin{proof}
 Assume that $A\in B(\cE\otimes F^2(H_n))$ is a multi-Toeplitz  and let
 $$
\varphi(S_1,\ldots, S_n):=\sum_{k=1}^\infty\sum_{|\alpha|=k}
B_{(\alpha)}\otimes S_\alpha^* + A_{(0)}\otimes  I +
\sum_{k=1}^\infty\sum_{|\alpha|=k} A_{(\alpha)}\otimes  S_\alpha $$
 be its Fourier
representation, where the coefficients are given by \eqref{f-coef}.
Part (i) of this theorem  follows from the proof of Theorem
\ref{Fourier}. To prove part (ii), notice first that the operator
$\varphi(rS_1,\ldots, rS_n)$ is in   $\cA_n(\cE)^*+\cA_n(\cE)$.
Indeed,  since $S_i^* S_j=\delta_{ij} I$, $i,j=1,\ldots, n$,  one
can easily see that
 $$\left\|\sum_{|\alpha|=k}
A_{(\alpha)}\otimes  r^{|\alpha|}
S_\alpha\right\|=r^k\left\|\sum_{|\alpha|=k}
A_{(\alpha)}^*A_{(\alpha)}\right\|^{1/2}$$ and a similar equality
holds for the coefficients $B_{(\alpha)}$. Due to part (i), we
deduce  that the series above  are convergent in the operator norm.
This proves part (a).

Now, we prove that
 $\|\varphi(rS_1,\ldots, rS_n)\|\leq \|A\|$ for
$0\leq r<1$. Define the row contraction $T:=(T_1,\ldots, T_n)$,
where $T_i=rS_i$, $i=1,\ldots, n$. Let $V:=(V_1,\ldots, V_n)$ be the
minimal isometric dilation of $T$ on the Hilbert space
$\cK:=\cH\oplus [F^2(H_n)\otimes \cD]$, where $\cH:=F^2(H_n)$.
According to equation \eqref{iso-dil}, we have
\begin{equation}
\label{Mat} V_i=\left[
\begin{matrix}rS_i&0\\
\Delta_i& S_i\otimes I_\cD
\end{matrix}
\right],\qquad i=1,\ldots, n,
\end{equation}
 with respect to the decomposition $\cK=\cH\oplus
[F^2(H_n)\otimes \cD]$. Since $T$ is a pure row contraction,  we
must have  $\cK=M_V(\cL_*)$, where $\cL_*$ is the wandering subspace
defined by  relation \eqref{L*}.  Due to  Proposition 2.10 from
\cite{Po-isometric}, we have $\lim_{k\to\infty} \sum_{|\alpha|=k}
V_\alpha T_\alpha^*h=0$ for any $h\in \cH$.  This implies
\begin{equation}
\label{h=} h=\sum_{k=0}^\infty \sum_{|\alpha|=k} V_\alpha
\left(I_\cK-\sum_{i=1}^n V_iT_i^*\right)T_\alpha^*h
\end{equation}
for any $h\in \cH$. Define the unitary operator $U:\cK\to
F^2(H_n)\otimes \cD_*$ by setting
\begin{equation}
\label{U} U\left( \sum_{\alpha \in \FF_n^+} V_\alpha
\ell_\alpha\right):= \sum_{\alpha\in \FF_n^+} e_\alpha\otimes
\Phi_*(\ell_\alpha),
\end{equation}
where $\sum_{\alpha\in \FF_n^+} |\ell_\alpha|^2<\infty$,
$\ell_\alpha\in \cL_*$, and $\Phi_*$ is defined by relation
\eqref{fistea}. Notice that
\begin{equation}\label{UVi}
UV_i=(S_i\otimes I_{\cD_*})U,\quad i=1,\ldots, n.
\end{equation}
Now, we prove that
\begin{equation}
\label{PHU} P_{\cE\otimes\cH}\left[ (I_\cE\otimes U^*)(A\otimes I_{\cD_*})(I_\cE\otimes U)
\right]
|_{\cE\otimes\cH}=
 \varphi(rS_1,\ldots, rS_n), \quad  0\leq r<1.
\end{equation}
Since both sides are bounded operators, it is enough to prove the
equality on a dense subset of $\cE\otimes \cH=\cE\otimes F^2(H_n)$. Taking
 $h=e_\beta$,
$\beta\in \FF_n^+$, in relation \eqref{h=}, we obtain
$$
e_\beta=\sum_{k=0}^\infty \sum_{|\alpha|=k} V_\alpha
\left(I_\cK-r\sum_{i=1}^n V_iS_i^*\right)r^{|\alpha|}S_\alpha^*
e_\beta.
$$
According to \eqref{U}, the definition of $\Phi_*$, and the fact
that $S_\alpha^*e_\beta=0$ if $|\alpha|>|\beta|$, we deduce that
$$
U(e_\beta)=\sum_{\stackrel{\alpha\in \FF_n^+}{|\alpha|\leq
|\beta|}}e_\alpha\otimes r^{|\alpha|}D_{T^*}S_\alpha^* e_\beta.
$$
Notice that, for any $\gamma, \beta\in \FF_n^+$, we have
\begin{equation*}\begin{split}
& \left< (I_\cE\otimes U^*)(A\otimes I_{\cD_*})(I_\cE\otimes U)(x\otimes  e_\gamma),
 y\otimes e_\beta\right> \\
 &\qquad=
\left<(A\otimes I_{\cD_*})(I_\cE\otimes U)(x\otimes e_\gamma), y\otimes
\sum_{\stackrel{\alpha\in \FF_n^+}{|\alpha|\leq |\beta|}}
e_\alpha\otimes r^{|\alpha|} D_{T^*}S_\alpha^* e_\beta\right>\\
&\qquad= \left<(M_{|\gamma|,|\beta|}(S_1,\ldots, S_n)\otimes I_{\cD_*})
(I_\cE\otimes U)(x\otimes
e_\gamma), (I_\cE\otimes U)(y\otimes e_\beta)\right>,
\end{split}
\end{equation*}
where $M_{|\gamma|,|\beta|}(S_1,\ldots,
S_n):=\sum\limits_{\stackrel{\alpha\in \FF_n^+}{1\leq|\alpha|\leq
|\beta|}} B_{(\alpha)}\otimes S_\alpha^*+ A_{(0)}\otimes I+
\sum\limits_{\stackrel{\alpha\in \FF_n^+}{1\leq|\alpha|\leq
|\beta|}} A_{(\alpha)}\otimes  S_\alpha$.
Now, using relations \eqref{UVi} and \eqref{Mat}, we deduce that
\begin{equation*}
\begin{split}
&\left<(M_{|\gamma|,|\beta|}(S_1,\ldots, S_n)\otimes I_{\cD_*})(I_\cE\otimes U)
(x\otimes
e_\gamma), (I_\cE\otimes U) (y\otimes e_\beta)\right>\\
&\qquad= \left<(I_\cE\otimes U^*)
M_{|\gamma|,|\beta|} (S_1\otimes I_{\cD_*},\ldots, S_n\otimes I_{\cD_*})
(I_\cE\otimes U)(x\otimes e_\gamma), y\otimes e_\beta\right>\\
&\qquad= \left< M_{|\gamma|,|\beta|} (V_1,\ldots, V_n)(x\otimes
e_\gamma), y\otimes
e_\beta\right>\\
&\qquad= \left< M_{|\gamma|,|\beta|} (rS_1,\ldots, rS_n)(x\otimes
e_\gamma),
 y\otimes e_\beta\right>\\
&\qquad= \left< \left(\sum_{k=0}^\infty \sum_{|\alpha|=k}
B_{(\alpha)}\otimes r^{|\alpha|} S_\alpha^* + A_{(0)}\otimes I+
\sum_{k=0}^\infty \sum_{|\alpha|=k} A_{(\alpha)}\otimes
r^{|\alpha|} S_\alpha\right)(x\otimes e_\gamma), y\otimes
e_\beta\right>
\end{split}
\end{equation*}
for any $\gamma,\beta\in \FF_n^+$. Therefore,   relation \eqref{PHU}
holds. Hence, we deduce that
\begin{equation}
\label{sup-fi}\sup_{0\leq r<1}\|\varphi(rS_1,\ldots, rS_n)\| \leq
\|A\|,
\end{equation}
which proves part (ii).
Now,
 we prove that
\begin{equation}\label{a-SOT}
  A=\text{\rm SOT-}\lim\limits_{r\to 1}
\varphi(rS_1,\ldots, rS_n). \end{equation}
First notice that, since $\sum_{|\alpha|\geq 1}A_{(\alpha)}^*
A_{(\alpha)} $ and $\sum_{|\alpha|\geq 1}B_{(\alpha)}
B_{(\alpha)}^*$ are  WOT convergent,  we have
\begin{equation}\label{var-var}
\|\varphi(rS_1,\ldots, rS_n) p-\varphi(S_1,\ldots, S_n) p\|\to
0,\quad\text{as}\quad r\to 1,
\end{equation}
 for any  vector-valued polynomial $~p\in\cE\otimes \cP\subset \cE \otimes F^2(H_n)$,
 where $\cP\subset F^2(H_n)$ is the set of all polynomials in $e_1,\ldots, e_n$.
Given  $~\epsilon>0~$ and $~h\in \cE\otimes F^2(H_n)$, there exists
a polynomial $~p\in \cE\otimes \cP$ such that
$
\|h-p\|\le \frac{\epsilon}{2\|A\|}.
$
Hence, and  using  the fact that   $\|\varphi(rS_1,\ldots,
rS_n)\|\leq \|A\|$ for $0\leq r<1$,  we deduce that
\begin{equation*}
\begin{split}
\|\varphi(rS_1,\ldots, rS_n)h-A h\| &\le \|\varphi(rS_1,\ldots,
rS_n)(h-p)\|+
\|(\varphi(rS_1,\ldots, rS_n)-\varphi(S_1,\ldots, S_n)) p\|
+ \|Ap- Ah\|\\
&\le\|\varphi(rS_1,\ldots, rS_n)\| \|h-p\| +\|(\varphi(rS_1,\ldots,
rS_n)-\varphi(S_1,\ldots, S_n)) p\|
+ \|A\| \|h-p\|\\
&\leq 2\|A\| \|h-p\|+ \|(\varphi(rS_1,\ldots,
rS_n)-\varphi(S_1,\ldots, S_n))
 p\|\\
&\le \epsilon+\|(\varphi(rS_1,\ldots, rS_n)-\varphi(S_1,\ldots,
S_n)) p\|.
\end{split}
\end{equation*}
Therefore, due to \eqref{var-var}, we obtain $
\limsup_{r\to1}\|\varphi(rS_1,\ldots, rS_n) h- A h\|\le\epsilon $
 for any $~\epsilon >0$.
Hence, $ \lim_{r\to 1}\|\varphi(rS_1,\ldots, rS_n)h -Ah\|=0, $
which  implies  relation \eqref{a-SOT}  and, therefore, part (b)
holds.

Conversely, assume that the coefficients
$\{A_{(\alpha)}\}_{\alpha\in \FF_n^+}$,
 $\{B_{(\alpha)}\}_{\FF_n^+\backslash\{g_0\}}$
 satisfy the conditions
(i) and (ii). Let us show that
 $\sup\limits_{q\in \cE\otimes \cP, \|q\|=1}\|\varphi(S_1,\ldots, S_n)q\|<\infty$. If
 this was not the case,
  then, for any $M>0$, there would be  a polynomial $q\in
\cE\otimes \cP$ with $\|q\|=1$ such that $\|\varphi(S_1,\ldots,
S_n)q\|>M$. Since $\|\varphi(rS_1,\ldots, rS_n)q-\varphi(S_1,\ldots,
S_n)q\|\to 0$ as $r\to 1$, there is $r_0\in (0,1)$ such that
$\|\varphi(r_0S_1,\ldots, r_0S_n)q\|> M$. Hence
$\|\varphi(r_0S_1,\ldots, r_0S_n)\|\geq\|\varphi(r_0S_1,\ldots,
r_0S_n)q\|>M,$ which  contradicts (ii).
  Consequently, $\sup\limits_{q\in \cE\otimes \cP, \|q\|=1}\|\varphi(S_1,\ldots,
S_n)q\|<\infty$, and, therefore, there is a unique operator  $A\in
B(\cE\otimes F^2(H_n))$ such that $Aq=\varphi(S_1,\ldots, S_n)q$ for
any polynomial $q\in \cE\otimes \cP$. As in the  the proof of part
(b), one can show that $A=\text{\rm{SOT-}}\lim\limits_{r\to
1}\varphi(rS_1,\dots,rS_n) $. Hence and using that
$\varphi(rS_1,\dots,rS_n)$ is a multi-Toeplitz  operator, i.e.,
$$
(I_\cE\otimes R_i^*)\varphi(rS_1,\dots,rS_n) (I_\cE\otimes R_j)=\delta_{ij}
\varphi(rS_1,\dots,rS_n), \quad i,j=1,\ldots, n,
$$
we deduce that $A$ is also a multi-Toeplitz  operator, which
completes the proof of the converse.

Now, we prove   part (c) of the theorem. If  $\epsilon>0$, then
there exists a polynomial $q\in \cE\otimes \cP$ with $\|q\|=1$ such
that $ \|A q\|=\|\varphi(S_1,\ldots, S_n)q\|>\|A\|-\epsilon.$
Since
$A=\text{\rm{SOT-}}\lim\limits_{r\to 1}\varphi(rS_1,\dots,rS_n) $,
there exists $r_0\in (0,1)$ such that $\|\varphi(r_0S_1,\ldots,
r_0S_n)q\|>\|A\|-\epsilon$. Using now relation \eqref{sup-fi}, we
deduce that
\begin{equation}
\label{sup-norm}
 \sup_{0\leq r<1}\|\varphi(rS_1,\ldots,
rS_n)\|=\|A\|.
\end{equation}

Now,  let $r_1,r_2\in [0,1)$ with $r_1<r_2$. Since  the operator
$g(S_1,\ldots, S_n):=\varphi(r_2S_1,\ldots, r_2S_n)$ is in the
operator system $\cA_n(\cE)^*+ \cA_n(\cE)$, the noncommutative von
Neumann inequality  \cite{Po-von} (see \cite{vN} for the classical
case) implies
  $\|g(rS_1,\ldots, rS_n)\|\leq \|g(S_1,\ldots, S_n)\|$ for
any $0\leq r<1$. In particular, when $r:=\frac {r_1}{r_2}$, we
deduce that $ \|\varphi(r_1S_1,\ldots, r_1S_n)\|\leq
\|\varphi(r_2S_1,\ldots, r_2S_n)\|.$ Consequently,  the   function
$[0,1]\ni r\to \|\varphi(rS_1,\ldots, rS_n)\|\in \RR^+$ is
increasing. Hence, and using  relation \eqref{sup-norm}, we complete
the proof.
\end{proof}

\begin{corollary}\label{multi}
The set of all multi-Toeplitz operators  on $\cE\otimes F^2(H_n)$
coincides with
$$
 \overline{\cA_n(\cE)^*+\cA_n(\cE)}^{WOT}= \overline{\cA_n(\cE)^*+\cA_n(\cE)}^{SOT},
 $$
 where $\cA_n(\cE):= B(\cE)\otimes_{min} \cA_n$ and $\cA_n$ is the noncommutative  disc algebra.
\end{corollary}

\begin{proof}
If $A$ is a multi-Toeplitz operator and $\varphi(S_1,\ldots, S_n)$
is its Fourier representation, then, according to  Theorem
\ref{Toeplitz}, $\varphi(rS_1,\ldots, rS_n)$ is in $\cA_n(\cE)^*+
\cA_n(\cE)$ for any $r\in [0,1)$ and $A=\text{\rm SOT-}\lim_{r\to
1}\varphi(rS_1,\ldots, rS_n)$. Therefore, $A$ is in
$\overline{\cA_n(\cE)^*+\cA_n(\cE)}^{SOT}$. Conversely,
  since any
operator $X\in \cA_n(\cE)^*+\cA_n(\cE)$  satisfies the equation
$(I_\cE\otimes R_i^*)X(I_\cE\otimes R_j)=\delta_{ij} X$ for
$i,j=1,\ldots n$, so does any operator  $T\in
\overline{\cA_n(\cE)^*+\cA_n(\cE)}^{SOT}$. Therefore, $T$ is a
multi-Toeplitz operator. If  $T\in
\overline{\cA_n(\cE)^*+\cA_n(\cE)}^{WOT}$, an argument as above
shows that $T$ is a multi-Toeplitz operator and, due to the first
part of the proof, we deduce that  $T\in
\overline{\cA_n(\cE)^*+\cA_n(\cE)}^{SOT}$. Since the other inclusion
is clear, the  proof is complete.
\end{proof}

We  remark that all the results of this section have appropriate versions for the
multi-Toeplitz  operators with respect to the left
 creation operators on the full Fock space.

\bigskip

\section{ Noncommutative Berezin transforms and
  free pluriharmonic functions}

We introduce  noncommutative Berezin transforms associated with
(completely)  bounded linear maps on $B(F^2(H_n))$, which will play
an important role in the study of free pluriharmonic functions and
their boundary behavior. First, we present some of their  properties
and connections to the classical case \cite{Be} and the
noncommutative Poisson transform \cite{Po-poisson}. Then we work out
some basic properties of the free pluriharmonic functions on the
noncommutative ball $[B(\cH)^n]_1$, including a Poisson mean value
property, Weierstrass type convergence  theorem, Harnack type
inequality (resp.~convergence theorem), and a maximum (resp.~minimum) principle. The free holomorphic functional calculus for
$n$-tuples of operators \cite{Po-holomorphic} is extended to free
pluriharmonic function.

  Let $\cH$ be a Hilbert space and  identify $M_m(B(\cH))$, the set of
$m\times m$ matrices with entries from $B(\cH)$, with
$B( \cH^{(m)})$, where $\cH^{(m)}$ is the direct sum of $m$ copies
of $\cH$.
Thus we have a natural $C^*$-norm on
$M_m(B(\cH))$. If $\cX$ is an operator space, i.e., a closed
subspace of $B(\cH)$, we consider $M_m(\cX)$ as a subspace of
$M_m(B(\cH))$ with the induced norm.
Let $\cX, \cY$ be operator spaces and $u:\cX\to \cY$ be a linear
map. Define the map
$u_m:M_m(\cX)\to M_m(\cY)$ by
$ u_m ([x_{ij}]):=[u(x_{ij})]. $
We say that $u$ is completely bounded  if
$ \|u\|_{cb}:=\sup_{m\ge1}\|u_m\|<\infty. $
If $\|u\|_{cb}\leq1$
(resp. $u_m$ is an isometry for any $m\geq1$) then $u$ is completely
contractive (resp. isometric),
and if $u_m$ is positive for all $m$, then $u$ is called
 completely positive. For basic results concerning  completely bounded maps
 and operator spaces we refer to \cite{P}, \cite{Pi}, and \cite{ER}.

Let $\cK$ be a Hilbert space and let  $\mu: B(F^2(H_n))\to B(\cK)$
be a completely bounded map. It is well-known (see e.g. \cite{P})
that  there exists a completely bounded linear map
$$\widetilde
\mu:=\mu\otimes \text{\rm id} : B(F^2(H_n)) \otimes_{min} B(\cH)\to
B(\cK)\otimes_{min} B(\cH)
$$
such that $ \widetilde \mu(f\otimes Y):= \mu(f)\otimes Y$ for $f\in
B(F^2(H_n)) $ and  $Y\in B(\cH)$. Moreover, $\|\widetilde
\mu\|_{cb}=\|\mu\|_{cb}$ and, if $\mu$ is completely positive, then
so is $\widetilde \mu$.
 We introduce  a  {\it
noncommutative Berezin transform} associated with $\mu$ as the map
$$\cB_\mu:
B(F^2(H_n))\times [B(\cH)^n]_1\to B(\cK)\otimes_{min} B(\cH)$$
defined by
\begin{equation}
\label{Berezin}
 \cB_\mu(f,X):=\widetilde{\mu}\left[ B_X^*(f\otimes
I_\cH)B_X\right], \qquad f\in B(F^2(H_n)), \ X:=(X_1,\ldots,X_n)\in
[B(\cH)^n]_1,
\end{equation}
  where the operator  $B_X \in B(F^2(H_n)\otimes
\cH)$  is defined
  by
\begin{equation}
\label{KTR}
 B_X := (I_{F^2(H_n)} \otimes \Delta_X)
 \left(I-R_1\otimes X_1^*-\cdots -R_n\otimes X_n^*\right)^{-1}
\end{equation}
and  $\Delta_X:=(I_\cH-\sum_{i=1}^n X_iX_i^*)^{1/2}$. We remark that
the {\it reconstruction operator} $$R_X:=R_1\otimes X_1^*+\cdots
+R_n\otimes X_n^*$$ has played an important role in noncommutative
multivariable operator theory (see \cite{Po-unitary},
\cite{Po-holomorphic}).
 Note that, due
to the fact that $R_1,\ldots, R_n$ are isometries with orthogonal
ranges, we have  $\|R_X\|= \|X\| $
 and, therefore, the
operator $B_X$ is well-defined. We also remark that the
noncommutative Berezin transform is well-defined even if the $n$-tuple
$X:=(X_1,\ldots, X_n)\in [B(\cH)^n]_1^-$  has  joint spectral
radius $r(X_1,\ldots, X_n)<1$. We recall   that the joint spectral
radius is defined by
$$
r(X_1,\ldots,X_n):=\lim_{k\to\infty} \left\|\sum_{|\alpha|=k}
X_\alpha X_\alpha^*\right\|^{\frac{1} {2k}}
$$
  and it is also equal to the spectral radius of the reconstruction operator
 $R_X$ (see \cite{Po-unitary}).
Consequently, $r(T_1,\ldots,T_n)<1$ if and only if  the spectrum of
$  R _X$  is included in $ \DD$.

\begin{theorem}\label{berezin}
Let \ $\cB_\mu$ be the  noncommutative Berezin transform  associated
with a completely bounded linear map $\mu: B(F^2(H_n))\to B(\cK)$.
\begin{enumerate}
\item[(i)] If $X\in [B(\cH)^n]_1$ is fixed, then
$$\cB_\mu(\cdot, X):B(F^2(H_n))\to B(\cK)\otimes_{min} B(\cH)$$ is a
completely bounded linear map with $ \|\cB_\mu(\cdot, X)\|_{cb}\leq
\|\mu\|_{cb}  \|B_X\|^2. $
\item[(ii)] If $\mu$ is selfadjoint, then
$\cB_\mu(f^*,X)=\cB_\mu(f,X)^*. $
 Moreover, if $\mu$ is completely positive, then so is the map \ $\cB_\mu(\cdot,
 X)$.
\item[(iii)] If $f\in B(F^2(H_n))$ is fixed, then the map
$$\cB_\mu(f, \cdot): [B(\cH)^n]_1 \to B(\cK)\otimes_{min} B(\cH)$$ is
continuous  and $\|\cB_\mu(f,X)\|\leq \|\mu\|_{cb}\|f\| \|B_X\|^2 $
for any    $ X\in [B(\cH)^n]_1$.
\end{enumerate}
\end{theorem}
\begin{proof}
The items  (i) and (ii) follow easily from the definition of the
noncommutative Berezin transform. To prove part (iii), let $X,Y\in
[B(\cH)^n]_1$ and notice that
\begin{equation*}
\begin{split}
\|\cB_\mu(f,X)-\cB_\mu(f,Y)\|&\leq \|\mu\|\|B_X^*(f\otimes
I_\cH)(B_X-B_Y)\|+ \|\mu\|\|(B_X^*-B_Y^*)(f\otimes I_\cH)B_Y\|\\
&\leq \|\mu\|\|f\|\|B_X-B_Y\|\left(\|B_X\|+\|B_Y\|\right).
\end{split}
\end{equation*}
The continuity of the map $X\mapsto \cB_\mu(f, X)$ will follow once
we prove that $X\mapsto B_X$ is  a continuous  map on
$[B(\cH)^n]_1$. To this end, notice that
\begin{equation}\label{B-B}
\|B_X-B_Y\|\leq \|\Delta_X\|
\|(I-R_X)^{-1}-(I-R_Y)^{-1}\|+\|\Delta_X-\Delta_Y\|\|(I-R_X)^{-1}\|.
\end{equation}
Since $\|R_X-R_Y\|=\|X-Y\|$, the map $X\mapsto R_X$ is continuous on
$[B(\cH)^n]_1$. Taking into account that $\|R_X\|<1$ for any $X\in
[B(\cH)^n]_1$, we deduce that $X\mapsto (I-R_X)^{-1}$  is also a
continuous map  on $[B(\cH)^n]_1$. Due to \eqref{B-B}, it remains to
show that the function $X\mapsto \Delta_X$ is continuous on
$[B(\cH)^n]_1$. By Weierstrass approximation theorem, for any
$\epsilon>0$ there exists a polynomial $p$ of $\lambda$ such that
$\sup_{\lambda\in [0,1]}
|p(\lambda)-\sqrt{\lambda}|<\frac{\epsilon}{3}$. Due to the
representation theorem for normal operators, we have
\begin{equation}
\label{pdpd} \|p(\Delta_X^2)-\Delta_X\|<\frac{\epsilon}{3} \ \text{
and } \ \|p(\Delta_Y^2)-\Delta_Y\|<\frac{\epsilon}{3}.
\end{equation}
Note also that
$$
\|\Delta_X^2-\Delta_Y^2\|\leq \|X-Y\|(\|X\|+\|Y\|)\leq 2\|X-Y\|.
$$
Consequently, since $p$ is a polynomial, there exists $\delta>0$
such that $\|p(\Delta_X^2)-p(\Delta_Y^2)\|<\frac{\epsilon}{3}$ if
$\|X-Y\|<\delta$ and $X,Y\in [B(\cH)^n]_1$. Now, using relation
\eqref{pdpd}, we deduce that
\begin{equation*}
\|\Delta_X-\Delta_Y\|\leq
\|\Delta_X-p(\Delta_X^2)\|+\|p(\Delta_X^2)-p(\Delta_Y^2)\|+\|p(\Delta_Y^2)-\Delta_Y\|\leq
\epsilon
\end{equation*}
if $\|X-Y\|<\delta$, which proves the continuity of  the map
$X\mapsto \Delta_X$. Therefore, the map $X\mapsto \cB_\mu(f, X)$ is
continuous on $[B(\cH)^n]_1$. The inequality in  (iii) is obvious.
The proof is complete.
\end{proof}

In what follows we present two  particular  cases of the
noncommutative Berezin transform which will play an important role
in this paper.

\smallskip

\leftline{\bf The Berezin transform $B_\mu(I, \cdot \,)$.}

If  $f=I$, the identity on $F^2(H_n)$, then the Berezin transform
$B_\mu(I, \cdot \,)$ coincides with the noncommutative Poisson
transform $\cP\mu$ associated with $\mu$, which  will be discussed
in Section 5.   We will show that, for any $X:=(X_1,\ldots,X_n)\in [B(\cH)^n]_1$,
$$\cB_\mu(I,X)=
\sum_{k=1}^\infty\sum_{|\alpha|=k} \mu(R_{\widetilde \alpha})\otimes
X_\alpha^* + \mu(I)\otimes  I+\sum_{k=1}^\infty\sum_{|\alpha|=k}
\mu(R_{\widetilde \alpha}^*) \otimes X_\alpha,
$$
where the convergence is in the operator norm topology of
$B(\cK\otimes \cH)$. Consider  the particular case when $n=1$,
$\cH=\cK=\CC$, $X=re^{i\theta}\in \DD$,  and $\mu$ is a complex
Borel measure on $\TT$. Since $\mu$ can be seen as a bounded linear
functional on $C(\TT)$, there is a unique bounded linear functional
$\hat\mu$ on the operator system $A(\DD)^*+A(\DD)$ (here $A(\DD)$ is
the disc algebra generated by the unilateral shift $S$ acting on the
Hardy space $H^2(\DD)$) such that $\hat\mu(S^k)=\mu(e^{ikt})$ if
$k\geq 0$, and $\hat\mu({S^*}^k)=\mu(e^{-ikt})$ if $k\geq 1$.
Indeed, if $p$ is any polynomial of the form $p(\overline{\lambda},
\lambda)= \sum_{k=1}^q b_k \overline{\lambda}^k+\sum_{k=0}^r a_k
\lambda^k$, then, using the noncommutative  von Neumann inequality
(when $n=1$), we obtain
\begin{equation*}
\begin{split}
|\hat\mu(p(S^*,S))|&=|\mu(p(e^{-it}, e^{it}))|\leq
\|\mu\|\sup_{e^{it}\in\TT}|p(e^{-it}, e^{it})| \leq
\|\mu\|\|p(S^*,S)\|,
\end{split}
\end{equation*}
which proves our assertion. Now, it is easy to see that the
noncommutative  Berezin transform $B_{\hat\mu}(I, \cdot \,)$
coincide with the classical Poisson transform of $\mu$, i.e., $
 \frac{1}{2\pi}\int_{-\pi}^\pi
P_r(\theta-t)d\mu(t), $ where
$P_r(\theta-t)=\frac{1-r^2}{1-2r\cos(\theta- t)+ r^2}$ is the
Poisson kernel.

Throughout  this  paper, the Berezin transform  $B_\mu(I, \cdot \,)$
will be denoted by $\cP\mu$ and called
  the (noncommutative) Poisson transform
 of $\mu$.

Next, we show that the noncommutative Poisson transform introduced
in \cite{Po-poisson} is in fact a particular case of the
noncommutative Berezin transform.

\smallskip
\leftline{\bf The Berezin transform $B_\tau$.} Let $\tau$ be the
linear functional on $B(F^2(H_n))$ defined by $\tau(f):=\left<
f(1),1\right>$. If $X\in [B(\cH)^n]_1$ is fixed, then
$\cB_\tau(\cdot, X): B(F^2(H_n))\to   B(\cH)$ is a completely
contractive linear map and
\begin{equation*}
\begin{split}
\left<\cB_\tau(f, X)x,y\right> =\left< B_X^*( f \otimes I_\cH)
B_X(1\otimes x),1\otimes y\right>,\quad x,y\in \cH.
\end{split}
\end{equation*}
We remark that $\cB_\tau(\cdot, X)$
   coincides with the noncommutative Poisson transform $P_X$
 introduced in
 \cite{Po-poisson}.  More precisely,  we have
 $$\cB_\tau(f, X)=P_X(f):=K_X^*(f\otimes I)K_X,
 $$
 where  $K_X=B_X |_{1\otimes \cH}:\cH\to F^2(H_n)\otimes
\cH$.
 We recall  from \cite{Po-poisson} that the restriction of $P_X$ to the Cuntz-Toeplitz $C^*$-algebra
$C^*(S_1,\ldots, S_n)$  (see \cite{Cu}) can be extended to the
closed
 ball $[B(\cH)^n]_1^-$ by setting
 \begin{equation}\label{Po-tran}
 P_X(f):=\lim_{r\to 1} K_{rX}^* (f\otimes I)K_{rX},\qquad X\in
 [B(\cH)^n]_1^-,\ f\in C^*(S_1,\ldots, S_n),
 \end{equation}
 where $rX:=(rX_1,\ldots, rX_n)$ and
    the limit exists in the operator
 norm topology of $B(\cH)$.
In this case we have
 \begin{equation}
 \label{SSXX}
 P_X(S_\alpha S_\beta^*)=X_\alpha X_\beta^*\ \text{ for any
  } \  \alpha,\beta\in \FF_n^+.
  \end{equation}
 When $X:=(X_1,\ldots, X_n)$  is a pure $n$-tuple, i.e., $\sum_{|\alpha|=k} X_\alpha
 X_\alpha^*\to 0$,  as $k\to\infty$,  in the strong operator topology,
 then  we have $P_X(f)=K_X^*(f\otimes I)K_X$. In particular, if $X=0$,
 then
 $P_0(f)= \left<f(1),1\right>I_\cH$.
We refer to \cite{Po-poisson}, \cite{Po-curvature},
  and \cite{Po-unitary} for more on
noncommutative Poisson transforms on $C^*$-algebras generated by
isometries.

If $f\in   B(F^2(H_n))$ is fixed, then $\cB_\tau(f, \cdot\, ):
[B(\cH)^n]_1 \to   B(\cH)$ is a bounded continuous map and
$\|\cB_\tau (f,X)\|\leq  \|f\| $ for any    $ X\in [B(\cH)^n]_1$.
If $n=1$, $\cH=\CC$, $X=\lambda\in \DD$, we recover the Berezin
transform of a bounded linear operator on the Hardy space
$H^2(\DD)$, i.e.,
$$
B_\tau(f,\lambda)=(1-|\lambda|^2)\left<f k_\lambda,
k_\lambda\right>,\quad f\in B(H^2(\DD)),
$$
where $k_\lambda(z):=(1-\overline{\lambda} z)^{-1}$ and  $z,
\lambda\in \DD$.

Throughout this paper,  the Berezin transform $B_\tau$ will be
called Poisson transform, to be in accord with the terminology used
in our  previous papers. If $f\in B(F^2(H_n))$ and $X\in
[B(\cH)^n]_1$, then the Poisson transform of $f$ at $X$ is given by
$$
P_X[f]:= B_\tau(f,X)=K_X^*(f\otimes I_\cH) K_X.
$$
This induces a completely contractive linear map
$$\text{\rm id}\otimes_{min} P_X:B(\cE)\otimes_{min} B(F^2(H_n))
\to B(\cE)\otimes_{min} B(\cH)
$$
 such that $(\text{\rm id}\otimes_{min} P_X)(Y\otimes f)=Y\otimes P_X[f]$ for any $Y\in B(\cE)$
 and $f\in B(F^2(H_n))$, where $\cE$ is a Hilbert space.
It is easy to see that
$$
(\text{\rm id}\otimes_{min} P_X)(u)=(I_\cE\otimes K_X^*)(u\otimes
I_\cH)(I_\cE\otimes K_X)
$$
for any $u\in B(\cE)\otimes_{min} B(F^2(H_n))$.
Given $X\in [B(\cH)]_1$, we define the operator-valued Poisson transform at $X$ to be
 the map
${\bf P}_X:B(\cE\otimes F^2(H_n))\to B(\cE\otimes \cH)$  defined by
\begin{equation}
\label{boldP}
 {\bf P}_X[u]:=(I_\cE\otimes K_X^*)(u\otimes
I_\cH)(I_\cE\otimes K_X)
\end{equation}
for any $u\in B(\cE\otimes F^2(H_n))$. It is clear that $ {\bf P}_X$
is an extension of the map  $\text{\rm id}\otimes_{min} P_X. $ In
the particular case when $\cE$ is finite dimensional, they coincide.

Now, we need to recall  from \cite{Po-holomorphic} a few facts
concerning free holomorphic functions on noncommutative balls.
Let $\{A_{(\alpha)}\}_{\alpha\in \FF_n^+}$ be a sequence
 of bounded linear operators  on a Hilbert space   $\cE$ and
 define
 $R\in [0,\infty]$ by setting
$$ \frac {1} {R}:= \limsup_{k\to\infty} \left\|\sum_{|\alpha|=k}
A_{(\alpha)}^* A_{(\alpha)}\right\|^{\frac{1} {2k}}. $$
The number  $R$ is called
 {\it radius of
convergence} of the formal power series $\sum_{\alpha\in \FF_n^+}
A_{(\alpha)} \otimes  Z_\alpha $
 in  noncommuting indeterminates
  $Z_1,\ldots, Z_n$,
where $Z_\alpha:=Z_{i_1}\cdots Z_{i_k}$ if $\alpha=g_{i_1}\cdots
g_{i_k}$
 and $Z_{g_0}:=I$.
Define the open noncommutative ball of radius $\gamma>0$,
 $$
 [B(\cH)^n]_\gamma:=\left\{ (X_1,\ldots, X_n)\in B(\cH)^n: \
 \left\|\sum_{i=1}^n X_iX_i^*\right\|^{1/2}<\gamma\right\}.
 $$
A map $F:[B(\cH)^n]_{\gamma}\to B(\cE)\otimes_{min} B( \cH)$ is
called a
  {\it free
holomorphic function} on  $[B(\cH)^n]_{\gamma }$  with coefficients
in $B(\cE)$ if there exist $A_{(\alpha)}\in B(\cE)$, $\alpha\in
\FF_n^+$, such that the formal power series $\sum_{\alpha\in
\FF_n^+} A_{(\alpha)} \otimes  Z_\alpha$ has radius of convergence
$\geq \gamma$ and such that $ F(X_1,\ldots,
X_n)=\sum\limits_{k=0}^\infty
\sum\limits_{|\alpha|=k}A_{(\alpha)}\otimes  X_\alpha, $ where the
series converges in the operator  norm topology  for any
$(X_1,\ldots, X_n)\in [B(\cH)^n]_{\gamma}$. We recall
(\cite{Po-holomorphic})
 that the following statements are equivalent:
 \begin{enumerate}
 \item[(i)] the series $\sum\limits_{k=0}^\infty
\sum\limits_{|\alpha|=k}A_{(\alpha)}\otimes X_\alpha$ is convergent
in the operator norm for any $(X_1,\ldots, X_n)\in
[B(\cK)^n]_\gamma$ and any Hilbert space $\cK$;
\item[(ii)] $\limsup\limits_{k\to\infty}\left\|\sum\limits_{|\alpha|=k}
A_{(\alpha)}^* A_{(\alpha)}\right\|^{1/2k}\leq \frac{1}{\gamma} $;
\item[(iii)]  the series
 $\sum_{k=1}^\infty \sum_{|\alpha|=k} A_{(\alpha )}\otimes r^{|\alpha|}S_\alpha$
 is convergent in the operator norm  for any $r\in [0,\gamma)$.
\end{enumerate}
The set of all free holomorphic functions on $[B(\cH)^n]_\gamma$
with coefficients in $B(\cE)$ is denoted by $Hol(B(\cH)^n_\gamma)$.
If the coefficients are scalars, we use the notation
$Hol_\CC(B(\cH)^n_\gamma)$.
 We say that $G$ is a  self-adjoint  {\it free pluriharmonic
function} on $[B(\cH)^n]_\gamma$ if there exists a free holomorphic
function $F$ on $[B(\cH)^n]_\gamma$ such that $G=\text{\rm Re}\, F$,
i.e.,
$$
G(X_1,\ldots, X_n)=\text{\rm Re}\, F(X_1,\ldots,
X_n):=\frac{1}{2}(F(X_1,\ldots, X_n)+ F(X_1,\ldots, X_n)^*)
$$
for  $(X_1,\ldots, X_n)\in [B(\cH)^n]_\gamma$.   $H$ is a free
pluriharmonic function on $[B(\cH)^n]_\gamma$ if
 $H:=H_1+iH_2$, where $H_1$ and $H_2$ are self-adjoint
free harmonic functions on $[B(\cH)^n]_\gamma$.
 According to \cite{Po-holomorphic}, any
free holomorphic on $[B(\cH)^n]_\gamma$ is continuous and uniformly
continuous on $[B(\cH)^n]_r^-$, $0\leq r <\gamma$. This implies
similar properties for free pluriharmonic functions. We remark that
in the particular case when $n=1$, a function is free pluriharmonic
on $[B(\cH)]_1$ if and only if it is harmonic on the open unit disc
$\DD$. Let $Har(B(\cH)^n_\gamma)$ be the set of all free
pluriharmonic functions on $[B(\cH)^n]_\gamma$ with operator-valued
coefficients. When the coefficients are scalars, we use the notation
$Har_\CC(B(\cH)^n_\gamma)$.

The following  result is an immediate consequence of the
above-mentioned properties.

\begin{proposition}\label{plu-ha}  A map $G:[B(\cH)^n]_\gamma\to
 B(\cE)\otimes_{min} B( \cH)$ is a
free pluriharmonic function  on $[B(\cH)^n]_\gamma$ with
coefficients in $B(\cE)$ if and only if there exist two sequences
 $\{A_{(\alpha)}\}_{\alpha\in \FF_n^+}\subset B(\cE)$ and
 $\{B_{(\alpha)}\}_{\alpha\in \FF_n^+\backslash \{g_0\}}\subset B(\cE)$
such that
$$\limsup\limits_{k\to\infty}\left\|\sum\limits_{|\alpha|=k}
A_{(\alpha)}^* A_{(\alpha)}\right\|^{1/2k}\leq \frac{1}{\gamma}
,\quad  \limsup\limits_{k\to\infty}\left\|\sum\limits_{|\alpha|=k}
B_{(\alpha)} B_{(\alpha)}^*\right\|^{1/2k}\leq \frac{1}{\gamma},
$$
and
\begin{equation} \label{ha} G(X_1,\ldots,
X_n)=\sum_{k=1}^\infty \sum_{|\alpha|=k} B_{(\alpha)}\otimes
X_\alpha^* +A_{(0)}\otimes  I+ \sum_{k=1}^\infty
 \sum_{|\alpha|=k} A_{(\alpha)}\otimes  X_\alpha,
\end{equation}
where the series are convergent in the  operator norm topology
  for any $ (X_1,\ldots, X_n)\in [B(\cH)^n]_\gamma$.
\end{proposition}

We remark that  two sequences of  operators
 $\{A_{(\alpha)}\}_{\alpha\in \FF_n^+}$ and
 $\{B_{(\alpha)}\}_{\alpha\in \FF_n^+\backslash \{g_0\}}$ in $B(\cE)$
 generate, by  relation \eqref{ha},  a  free pluriharmonic function $G:[B(\cH)^n]_\gamma\to
 B(\cE)\otimes_{min} B(\cH)$
if and only if
  the series
 $\sum_{k=1}^\infty
  \sum_{|\alpha|=k} A_{(\alpha)}\otimes  r^{|\alpha|}S_\alpha$ and
  $\sum_{k=1}^\infty
   \sum_{|\alpha|=k} B_{(\alpha)}\otimes  r^{|\alpha|}S^*_\alpha
   $
  are convergent    in the
   operator norm topology for any $r\in [0,\gamma)$. Moreover, if
   $\cH$ is infinite dimensional, then it is enough to assume the
   convergence in the operator norm of   the series in \eqref{ha}.
   Notice also that a
   free pluriharmonic function is uniquely determined by its
   representation on an infinite dimensional Hilbert space, in
   particular, on the full Fock space $F^2(H_n)$.

\begin{lemma}
\label{PPP} If $\gamma_1>0$ and  \ $0\leq\gamma_j\leq 1$ for
$j=2,\ldots, k$, then
$$
P_{\gamma_1\cdots \gamma_k X}=P_{\gamma_1 X}\circ P_{\gamma_2
S}\circ \cdots \circ P_{\gamma_k S}
$$
for any $X\in [B(\cH)^n]_{\frac{1}{\gamma_1}}$, where
$S:=(S_1,\ldots, S_n)$ is the $n$-tuple of left creation operators
on the Fock space $F^2(H_n)$. Moreover,
$$
{\bf P}_{\gamma_1\cdots \gamma_k X}[g]=\left({\bf P}_{\gamma_1
X}\circ {\bf P}_{\gamma_2 S}\circ \cdots \circ {\bf P}_{\gamma_k
S}\right)[g]
$$
for any $g\in B(\cE)\otimes_{min} B(F^2(H_n)$, where ${\bf P}_Y$ is
defined by \eqref{boldP}.
\end{lemma}

\begin{proof}
We recall that the Poisson transform of $f\in B(F^2(H_n))$ at $Y\in
 [B(\cH)^n]_1$ is given by
 $$P_Y[f]:=B_\tau(f,Y)=K_Y^*(f\otimes I_\cH) K_Y.
 $$
Now, we prove the result  for $k=2$. First, we show that
\begin{equation}
\label{PPP1}
 P_{\gamma_1\gamma_2 X}[g]=(P_{\gamma_1 X}\circ
P_{\gamma_2 S})[g]
 \end{equation}
 for any $g\in C^*(S_1,\ldots,
S_n)$. Let $p_m(S_1,\ldots, S_n):=\sum a_{\alpha,\beta}^{(m)}
S_\alpha S_\beta^*$, $m\in \NN$, be a sequence of polynomials in
 $C^*(S_1,\ldots, S_n)$ such that
$p_m(S_1,\ldots, S_n)\to g$ in the operator norm, as $m\to \infty$.
Due to the properties of the Poisson transform, we have
\begin{equation*}
\begin{split}
P_{\gamma_1 X}\left\{ P_{\gamma_2 S}[p_m(S_1,\ldots, S_n)]\right\}&
= K^*_{\gamma_1 X}\left\{\left[K^*_{\gamma_2 S}\left(p_m(S_1,\ldots,
S_n)\otimes I_{F^2(H_n)}\right)K_{\gamma_2 S}\right]\otimes
I_{\cH}\right\}K_{\gamma_1 X}\\
&=K^*_{\gamma_1 X}\left[ p_m\left(\gamma_2 S_1,\ldots, \gamma_2
S_n\right)\otimes I_\cH\right]K_{\gamma_1 X}
= p_m\left(\gamma_1 \gamma_2 X_1,\ldots, \gamma_1 \gamma_2 X_n\right)\\
&=K^*_{\gamma_1 \gamma_2 X}\left( p_m(S_1,\ldots, S_n)\otimes
I_\cH\right)K_{\gamma_1 \gamma_2 X} =P_{\gamma_1 \gamma_2
X}[p_m(S_1,\ldots, S_n)].
\end{split}
\end{equation*}
Since  the Poisson transform  is continuous in the operator norm
topology, we deduce relation \eqref{PPP1}. Recall that
$C^*(S_1,\ldots, S_n)$ contains the compact operators in
$B(F^2(H_n))$  (see \cite{Cu}) and any finite rank operator is
compact. Therefore,  $Q_{\leq m}f \in C^*(S_1,\ldots, S_n)$ for any
$f\in B(F^2(H_n))$, where $Q_{\leq m}:=I-\sum_{|\alpha|=m+1}
S_\alpha S_\alpha^*$ is the orthogonal projection of $F^2(H_n)$ onto
the set of all polynomials of degree $\leq m$. Due to the first part
of the proof, we have $ P_{\gamma_1\gamma_2 X}[Q_{\leq m}
f]=(P_{\gamma_1 X}\circ P_{\gamma_2 S})(Q_{\leq m} f). $ Notice
also that    $\|P_{\gamma_2 S}[Q_{\leq m} f]\|\leq \|Q_{\leq m}
f\|\leq \|f\|$  for any $m\in \NN$, and
 $\text{\rm
SOT-}\lim\limits_{m\to\infty} Q_{\leq m} f=f$. Since the map
$A\mapsto A\otimes I$ is SOT-continuous on bounded subsets of
$B(F^2(H_n))$, the above equality implies $ P_{\gamma_1\gamma_2
X}[f]=(P_{\gamma_1 X}\circ P_{\gamma_2 S})[ f] $ for any $f\in
B(F^2(H_n))$.   The general result follows easily by iteration. Now,
the second equality  can be easily  deduced.  This completes the
proof.
\end{proof}

Let $\cP_n$  be the set of all polynomials in $S_1,\ldots, S_n$ and
the identity, and   denote by $\cP_n(\cE)$ the spatial tensor
product $B(\cE)\otimes \cP_n$. The next result shows that a free
pluriharmonic function is uniquely determined by the Poisson mean
value property and the radial function.

\begin{theorem}\label{MVP}
If   $u:[B(\cH)^n]_\gamma\to B(\cE)\otimes_{min} B( \cH)$ is  a free
pluriharmonic function , then
\begin{enumerate}
\item[(i)] $u(rS_1,\ldots, rS_n)\in \overline{\cP_n(\cE)^*+
\cP_n(\cE)}^{\|\cdot\|}$ for any $r\in [0,\gamma)$, and
\item[(ii)]
$u$ has the Poisson mean value property, i.e., $ u(X_1,\ldots,
X_n)={\bf P}_{\frac{1}{r} X}[u(rS_1,\ldots, rS_n)] $ for any
$X:=(X_1,\ldots, X_n)\in [B(\cH)^n]_r$ and $r\in (0,\gamma)$.
\end{enumerate}

Conversely,  if there exists a map $\varphi:[0,\gamma)\to
\overline{\cP_n(\cE)^*+ \cP_n(\cE)}^{\|\cdot\|}$ such that
\begin{equation}
\label{fipifi}
 \varphi(r)={\bf  P}_{\frac{r}{t} S}[\varphi(t)]
\end{equation}
  for  $0\leq r<t<\gamma$,
   then the map
$v:[B(\cH)^n]_\gamma\to B(\cE)\otimes_{min} B( \cH)$ defined by
\begin{equation}\label{v-def}
v(X_1,\ldots, X_n):={\bf  P}_{\frac{1}{r} X} [\varphi(r)]
 \end{equation}
   for
any  $ X:=(X_1,\ldots, X_n)\in [B(\cH)^n]_r$   and  $ r\in
(0,\gamma)$,
  is a free pluriharmonic function. Moreover,
$v(rS_1,\ldots, rS_n)=\varphi(r)$ for any $r\in [0,\gamma)$.
\end{theorem}

\begin{proof}
Assume that $u$ is a free pluriharmonic function  and has the
representation
$$
u(Y_1,\ldots, Y_n)=\sum_{k=1}^\infty \sum_{|\alpha|=k}
B_{(\alpha)}\otimes  Y_\alpha^* + A_{(0)}\otimes  I
+\sum_{k=1}^\infty \sum_{|\alpha|=k} A_{(\alpha)}\otimes  Y_\alpha
$$
for any $(Y_1,\ldots, Y_n)\in [B(\cH)^n]_\gamma$. Since the series
above are convergent in the operator norm topology, one can easily
see that $u(rS_1,\ldots, rS_n)\in \cA_n(\cE)^*+\cA_n(\cE)\subset
\overline{\cP_n(\cE)^*+ \cP_n(\cE)}^{\|\cdot\|}$ for any $r\in
[0,\gamma)$. Denote $q_m(S_1,\ldots, S_n):=\sum_{0<|\alpha|\leq m}
B_{(\alpha)}\otimes  S_\alpha^*+ \sum_{0\leq|\alpha|\leq m}
A_{(\alpha)}\otimes  S_\alpha$, \ $m\in\NN$, and notice that
relation \eqref{SSXX} implies $ q_m(Y_1,\ldots, Y_n)={\bf
P}_{\frac{1}{r} Y}[q_m(rS_1,\ldots, rS_n)] $ for any
$Y:=(Y_1,\ldots, Y_n)\in [B(\cH)^n]_r$ and $r\in (0,\gamma)$. Taking
into account that  the Poisson transform is completely contractive,
that $q_m(Y_1,\ldots, Y_n)\to u(Y_1,\ldots, Y_n)$ and
$q_m(rS_1,\ldots, rS_n)\to u(rS_1,\ldots, rS_n)$, as $m \to \infty$,
we deduce item (ii).

Conversely, assume that the map $ \varphi$ has the properties stated
in the theorem and fix $r\in (0,\gamma)$.
 Due to Corollary \ref{multi},  $ \varphi(r)$
is  a multi-Toeplitz operator. By Theorem \ref{Fourier},
$\varphi(r)$  has a unique Fourier representation
$$
\sum_{|\alpha|>0}  B_{(\alpha)} (r)\otimes r^{|\alpha|}S_\alpha^* +
A_{(0)}(r)\otimes  I+ \sum_{|\alpha|>0}  A_{(\alpha)} (r)\otimes
r^{|\alpha|}S_\alpha,
$$
where $\{A_{(\alpha)}(r)\}_{\alpha\in \FF_n^+}$ and
 $\{B_{(\alpha)}(r)\}_{\alpha\in \FF_n^+\backslash \{g_0\}}$ are some
 sequences  of operators in $B(\cE)$. Applying Theorem \ref{Toeplitz}, we
 deduce that the map $h:[B(\cH)^n]_1\to B(\cE)\otimes_{min} B( \cH)$ defined by
 \begin{equation}
\label{gy}
 h(Z_1,\ldots, Z_n):=\sum_{k=1}^\infty \sum_{|\alpha|=k}
  B_{(\alpha)} (r)\otimes r^{|\alpha|}Z_\alpha^* + A_{(0)}(r)
 \otimes I+ \sum_{k=1}^\infty \sum_{|\alpha|=k} A_{(\alpha)} (r)\otimes r^{|\alpha|}Z_\alpha
 \end{equation}
 is a free pluriharmonic function,
 where the series are convergent in the operator norm topology.
Choose a sequence  of polynomials  $\{p_m(S_1,\ldots,
S_n)\}_{m=1}^\infty$  in $\cP_n(\cE)^*+ \cP_n(\cE)$, such that
$\|p_m(rS_1,\ldots, rS_n)- \varphi(r)\|\to 0$, as $m\to \infty$.
Applying  again Theorem \ref{Toeplitz} to the multi-Toeplitz
operator $A:= \varphi(r)- p_m(rS_1,\ldots, rS_n)$, we deduce that
\begin{equation}
\label{ht}
 \|h(tS_1,\ldots, tS_n)-p_m(rtS_1,\ldots, rtS_n)\|\leq
\|\varphi(r)- p_m(rS_1,\ldots, rS_n)\|
\end{equation}
 for any $t\in [0,1)$. If $Y:=(Y_1,\ldots, Y_n)\in [B(\cH)^n]_r$,
 there exists $t_0\in (0,1)$ such that $\frac{1}{t_0 r} Y\in
 [B(\cH)^n]_1$.
 Due to the noncommutative von Neumann inequality, we have
 $$
 \left\|h\left(\frac{1}{r}
 Y_1,\ldots,\frac{1}{r}Y_n\right)-p_m(Y_1,\ldots, Y_n)\right\|\leq
 \|h(t_0S_1,\ldots, tS_n)-p_m(t_0rS_1,\ldots, t_0rS_n)\|.
 $$
 Hence and using \eqref{ht}, we obtain
 $$
\left\|h\left(\frac{1}{r}
 Y_1,\ldots,\frac{1}{r}Y_n\right)-p_m(Y_1,\ldots, Y_n)\right\|\leq
\|\varphi(r)- p_m(rS_1,\ldots, rS_n)\|
$$
for any $(Y_1,\ldots, Y_n)\in [B(\cH)^n]_r$. This implies that $
p_m(Y_1,\ldots, Y_n)$ converges in the norm topology to
$h\left(\frac{1}{r}
 Y_1,\ldots,\frac{1}{r}Y_n\right)$, as
$m\to \infty$.
 Since $ p_m(Y_1,\ldots,
Y_n)={\bf  P}_{\frac{1}{r} Y}[p_m(rS_1,\ldots, rS_n)] $ and taking
the limit in the operator norm, as $m\to \infty$, we obtain
$
 h\left(\frac{1}{r}
 Y_1,\ldots,\frac{1}{r}Y_n\right)={\bf  P}_{\frac{1}{r} Y}[ \varphi(r)]
$
for any $(Y_1,\ldots, Y_n)\in [B(\cH)^n]_r$.
 Hence,  and using relation \eqref{gy}, we deduce
that
\begin{equation}
\label{PYF} {\bf   P}_{\frac{1}{r} Y}[
\varphi(r)]=\sum_{k=1}^\infty \sum_{|\alpha|=k}  B_{(\alpha)}
(r)\otimes Y_\alpha^* + A_{(0)}(r)
 \otimes I+ \sum_{k=1}^\infty \sum_{|\alpha|=k}  A_{(\alpha)} (r)\otimes Y_\alpha
 \end{equation}
for any $(Y_1,\ldots, Y_n)\in [B(\cH)^n]_r$. If $r<t<\gamma$, then,
as above, one can show that
\begin{equation}
\label{PYF2}
 {\bf  P}_{\frac{1}{t} Z}[ \varphi(t)]=\sum_{k=1}^\infty
\sum_{|\alpha|=k}  B_{(\alpha)} (t)\otimes Z_\alpha^* + A_{(0)}(t)
 \otimes I+ \sum_{k=1}^\infty \sum_{|\alpha|=k}  A_{(\alpha)} (t)\otimes Z_\alpha
 \end{equation}
for any $(Z_1,\ldots, Z_n)\in [B(\cH)^n]_t$.
On the other hand, due to relation \eqref{fipifi}, Lemma \ref{PPP},
and the fact that $\varphi(r)$ and $ \varphi(t)$ are in $
B(\cE)\otimes_{min} B(F^2(H_n))$, we deduce that
$$
{\bf  P}_{\frac{1}{r} Y} [ \varphi(r)]={\bf P}_{\frac{1}{r} Y}
\left({\bf  P}_{\frac{r}{t} S} [\varphi(t)]\right)={\bf
P}_{\frac{1}{t} Y}[ \varphi(t)]
$$
for  $Y:=(Y_1,\ldots, Y_n)\in [B(\cH)^n]_r$.
Hence, using relations \eqref{PYF},  \eqref{PYF2}, and the
uniqueness of free pluriharmonic functions,  we deduce that
  $ B_{(\alpha)} (r)= B_{(\alpha)} (t)$ for  $\beta\in
  \FF_n^+\backslash\{g_0\}$,
and $ A_{(\alpha)} (r)= A_{(\alpha)} (t)$ for  $\alpha\in\FF_n^+$.
Therefore, the coefficients do not depend on $r\in (0,\gamma)$, so
we may set
   $A_{(\alpha)}:=A_{(\alpha)}(r)$, $\alpha\in\FF_n^+$, and
$B_{(\alpha)}:=B_{(\alpha)}(r)$, $\beta\in
\FF_n^+\backslash\{g_0\}$. Now, it is  clear that the map
$v:[B(\cH)^n]_\gamma\to B(\cE)\otimes_{min} B(\cH)$ given by
\eqref{v-def} is well-defined and
$$
v(X_1,\ldots, X_n)=\sum_{k=1}^\infty \sum_{|\alpha|=k}  B_{(\alpha)}
\otimes X_\alpha^* + A_{(0)}\otimes
 I+ \sum_{k=1}^\infty \sum_{|\alpha|=k}  A_{(\alpha)}\otimes  X_\alpha
 $$
for any $(X_1,\ldots, X_n)\in [B(\cH)^n]_\gamma$, where the series
are convergent in the norm operator topology. Therefore, $v$ is a
free pluriharmonic function.  Due to Theorem \ref{Toeplitz},
$\varphi(r)=\lim_{t\to 1}
 h(tS_1,\ldots, tS_n)$. Since $v$ is continuous on
 $[B(\cH)^n]_\gamma$, we have $\lim_{t\to 1}
 h(tS_1,\ldots, tS_n)=\lim_{t\to 1}v(rtS_1,\ldots, rtS_n)=v(rS_1,\ldots, rS_n)$.
 Consequently, we have $\varphi(r)=v(rS_1,\ldots, rS_n)$.
 This completes the proof.
\end{proof}

 Now we obtain  a Weierstrass   type convergence theorem \cite{Co} for the vector space
 $Har(B(\cH)^n_\gamma)$, $\gamma>0$,
    of  all free pluriharmonic functions on the
    open   unit ball  $[B(\cH)^n]_\gamma$ with  coefficients  in $B(\cE)$.
This enables us to introduce a metric on $Har(B(\cH)^n_\gamma)$ with
respect to which it becomes a complete metric space.

Given $r\in [0,\gamma)$,  denote by  $ [B(\cH)^n]_{ r}^{-}$ the
noncommutative closed ball
$$ [B(\cH)^n]_{ r}^{-}:=\{ (X_1,\ldots, X_n)\in B(\cH)^n:\
 \|X_1X_1^*+\cdots+X_nX_n^*\|^{1/2}\leq r\}.
$$
Assume now that $\cH$ is an infinite dimensional  Hilbert space and
recall  that a free pluriharmonic function is  uniquely
determined by its  representation   $\cH$.

 Here
 is our  version of Weierstrass theorem   for
  free pluriharmonic functions.

\begin{theorem}\label{Weierstrass}
 Let $\{u_m\}_{m=1}^\infty\subset Har(B(\cH)^n_\gamma)$, $\gamma>0$,  be a
 sequence of free
  pluriharmonic  functions   which is uniformly convergent on any closed
  ball  $[B(\cH)^n]_r^-$, $r\in [0,\gamma)$. Then there is a free pluriharmonic
  function
     $u\in Har(B(\cH)^n_\gamma)$   such that
$u_m$ converges to $u$  on any closed  ball $[B(\cH)^n]_r^-$.
\end{theorem}

\begin{proof}
Since $\cH$ is infinite dimensional and due to the noncommutative
von Neumann inequality, one can see that a sequence
$\{u_m\}_{m=1}^\infty\subset Har(B(\cH)^n_\gamma)$ of free
pluriharmonic functions   converges uniformly on $[B(\cH)^n]_r^-$
  if and only if the sequence
  $\{u_m(rS_1,\ldots, rS_n)\}_{m=1}^\infty$ is convergent
  in the operator norm topology of $B(\cE\otimes F^2(H_n))$.
For each $m\in \NN$, $u_m(rS_1,\ldots, rS_n)$ is in $\cA_n(\cE)^*+
\cA_n(\cE)$ and $\varphi(r):=\lim\limits_{m\to\infty}
u_m(rS_1,\ldots, rS_n)$ is in $\overline{\cP_n(\cE)^*
+\cP_n(\cE)}^{\|\cdot\|}$. Since $u_m$ is free pluriharmonic and using
 the properties
of the Poisson transform,  we have $ u_m(rS_1,\ldots, rS_n)={\bf
P}_{\frac{r}{t} S}[u_m(t S_1,\ldots, t S_n)] $ for $0\leq
r<t<\gamma$. Taking the limit as $m\to \infty$, we obtain
$
 \varphi(r)={\bf  P}_{\frac{r}{t} S} [\varphi(t)].
$
On the other hand, due to Theorem \ref{MVP}, we have
\begin{equation}
\label{umP} u_m(X_1,\ldots, X_n)={\bf
P}_{\frac{1}{r}X}[u_m(rS_1,\ldots, rS_n)]
\end{equation}
for any $X:=(X_1,\ldots,X_n)\in [B(\cH)^n]_r$ and $r\in (0,\gamma)$.
Since $u_m$ converges uniformly on $[B(\cH)^n]_r^-$, there exists
$v(X):=\lim_{m\to\infty} u_m(X)$. Now, relation \eqref{umP} implies
$v(X)={\bf P}_{\frac{1}{r}X} [\varphi(r)]$ for any $X\in
[B(\cH)^n]_r$ and $r\in (0,\gamma)$. Applying again Theorem
\ref{MVP}, we deduce that $v$ is a free pluriharmonic function on
$[B(\cH)^n]_\gamma$. The proof is complete.
 \end{proof}

Let  $C(B(\cH)^n_\gamma, B(\cE\otimes \cH))$ be the vector space of
all continuous functions from
 the open  noncommutative  ball
$[B(\cH)^n]_\gamma$ to $B(\cE\otimes \cH)$. If $f,g\in
C(B(\cH)^n_\gamma, B(\cE\otimes \cH))$ and $0<r<\gamma$, we define
$$
\rho_r(f,g):=\sup_{(X_1,\ldots, X_n)\in [B(\cH)^n]_r^-}
\|f(X_1,\ldots, X_n)-g(X_1,\ldots, X_n)\|.
$$
 Let $0<r_m<\gamma$ be such that $\{r_m\}_{m=1}^\infty$ is an increasing
  sequence convergent  to $\gamma$.
For any  functions $f,g\in C(B(\cH)^n_\gamma, B(\cE\otimes\cH))$, we define
$$
\rho (f,g):=\sum_{m=1}^\infty \left(\frac{1}{2}\right)^m
\frac{\rho_{r_m}(f,g)}{1+\rho_{r_m}(f,g)}.
$$
 As in \cite{Po-holomorphic}, in the particular case when $\cE=\CC$,    one can   prove  that
 if
  $\{f_k\}_{k=1}^\infty$  and $f$  are functions in  $C(B(\cH)^n_\gamma, B(\cE\otimes\cH))$,
  then $f_k$ is convergent  to  $f$  in the metric $\rho$ if and only
  if $f_k\to f$ uniformly on any closed ball $[B(\cH)^n]_{r_m}^-$,
   \ $m=1,2,\ldots$.
Moreover, one can show that   $\left( C(B(\cH)^n_\gamma,
B(\cE\otimes\cH)), \rho\right)$ is a complete metric space.

\begin{theorem}\label{complete-metric}
$\left(Har(B(\cH)^n_\gamma), \rho\right)$  is a complete metric
space.
 \end{theorem}

\begin{proof}
Since $Har(B(\cH)^n_\gamma)\subset C(B(\cH)^n_\gamma,
B(\cE\otimes\cH))$ and $\left(C(B(\cH)^n_\gamma, B(\cE\otimes\cH)),
\rho\right)$ is a complete metric space, it is enough to show that
  $\left(Har(B(\cH)^n_\gamma), \rho\right)$ is
closed in $\left( C(B(\cH)^n_\gamma, B(\cE\otimes\cH)),
\rho\right)$. Let $\{u_k\}_{k=1}^\infty\subset Har(B(\cH)^n_\gamma)$
and $u\in C(B(\cH)^n_\gamma, B(\cE\otimes\cH)$ be  such that
$\rho(u_k,u)\to 0$, as $k\to\infty$. Consequently, $u_k\to u$
uniformly on any closed ball $[B(\cH)^n]_{r_m}^-$, \ $m=1,2,\ldots$.
Applying now Theorem \ref{Weierstrass}, we deduce that $u\in
Har(B(\cH)^n_\gamma)$.
  This completes the proof of the theorem.
 \end{proof}

We say  that a free pluriharmonic function $u$  is positive if any
representation on a Hilbert space is positive, i.e., $u(X_1,\ldots, X_n)\geq 0$
for any $(X_1,\ldots, X_n)\in  [B(\cK)^n]_\gamma$ and  any Hilbert space $\cK$.
Our next result is    a Harnack type inequality for free
pluriharmonic functions.

\begin{theorem}
\label{Harnack1} If $u$ is a  positive free pluriharmonic function
on $[B(\cH)^n]_\gamma$  with operator-valued  coefficients and $0<r<\gamma$, then
$$
\|u(X_1,\ldots, X_n)\|\leq \|u(0)\|\,\frac{\gamma+r}{\gamma-r}\quad
\text{ for any } \ (X_1,\ldots, X_n)\in [B(\cH)^n]_r^-.
$$
\end{theorem}

\begin{proof}
Any self-adjoint free pluriharmonic function on the noncommutative
ball $[B(\cH)^n]_\gamma$ has a  representation
$$
u(Y_1,\ldots, Y_n)=\sum_{k=1}^\infty \sum_{|\alpha|=k}
A_{(\alpha)}^*\otimes  Y_\alpha^*+ A_{(0)}\otimes  I+
\sum_{k=1}^\infty \sum_{|\alpha|=k} A_{(\alpha)}\otimes  Y_\alpha
$$
for $(Y_1,\ldots, Y_n)\in[B(\cH)^n]_\gamma$. Due to Theorem
\ref{MVP}, if $0<r<\gamma$ and $ X:=(X_1,\ldots,X_n)\in
[B(\cH)^n]_r$, then
\begin{equation}
\label{uPu} u(X_1,\ldots, X_n)= {\bf
P}_{\frac{1}{r}X}[u(rS_1,\ldots, rS_n)].
\end{equation}
Notice that the map $h:[B(\cH)^n]_1\to B(\cH)$ defined by $
h(Z_1,\ldots, Z_n):=u(\gamma Z_1,\ldots, \gamma Z_n)$ is a positive
free pluriharmonic function. Due to Theorem 3.1 from \cite{Po-Bohr}
(see also Lemma \ref{plu-pro}), we have
\begin{equation*}
 \left\| \sum_{|\alpha|=k} A_{(\alpha)}^* A_{(\alpha)}\right\|^{1/2}\leq
\frac{\|A_{(0)}\|}{\gamma^k}\quad \text{ for any }\  k\in \NN.
\end{equation*}
Using this inequality and relation \eqref{uPu}, we deduce that, for
any $ (X_1,\ldots, X_n)\in [B(\cH)^n]_r$,
\begin{equation*}
\begin{split}
\|u(X_1,\ldots, X_n)\|&\leq \|u(rS_1,\ldots, rS_n)\| \leq
\|A_{(0)}\| +2\sum_{k=1}^\infty\left\|\sum_{|\alpha|=k}
A_{(\alpha)}\otimes  r^{|\alpha|}S_\alpha\right\|\\
&= A_{(0)} +2\sum_{k=1}^\infty  r^k \left\|\sum_{|\alpha|=k}
A_{(\alpha)}^* A_{(\alpha)}\right\|^{1/2}
\leq \|A_{(0)}\| +2\sum_{k=1}^\infty  r^k \frac{\|A_{(0)}\|}{\gamma^k}\\
&= \|A_{(0)}\|\left(
1+\frac{2\frac{r}{\gamma}}{1-\frac{r}{\gamma}}\right) =\|A_{(0)}\|\,
\frac{\gamma+r}{\gamma-r},
\end{split}
\end{equation*}
which completes the proof.
\end{proof}

Now, we can obtain a Harnack type convergence theorem
 for free pluriharmonic functions.

\begin{theorem}
\label{Harnack2} Let $\{u_m\}_{m=1}^\infty$ be a sequence of
   free pluriharmonic functions on $[B(\cH)^n]_\gamma$
with operator-valued coefficients
such that  $\{u_m(0)\}_{m=1}^\infty$ is a  convergent  sequence in
the operator norm and
$$
 u_1\leq u_2\leq \cdots.
$$
 Then $u_m$ converges to a free
pluriharmonic function on $[B(\cH)^n]_\gamma$.
\end{theorem}
\begin{proof}
We may assume that $u_1\geq 0$ (if not, consider the sequence
$\{u_m-u_1\}_{m=1}^\infty$).  If $m>k$, then, applying the Harnack
type inequality of Theorem \ref{Harnack1} to the positive free
pluriharmonic function $u_m-u_k$, we obtain
$$
\|u_m(X)-u_k(X)\|\leq \|u_m(0)-u_k(0)\|\frac{\gamma+r}{\gamma-r}
$$
for any $X\in [B(\cH)^n]_r^-$.  Since $\{u_m(0)\}$ is convergent
in the operator norm, we
deduce that   $\{u_m\}_{m=1}^\infty$ is a uniformly Cauchy sequence
on $[B(\cH)^n]_r^-$. Applying Theorem \ref{Weierstrass}, we deduce
that $u_m$ converges to a free pluriharmonic function on
$[B(\cH)^n]_\gamma$. This completes the proof.
\end{proof}

We remark that if $\cE=\CC$ in Theorem \ref{Harnack2},  then  it is
enough to assume that the sequence $\{u_m(0)\}$ is bounded.

 The following result can be
seen as a maximum (resp. minimum) principle for free  pluriharmonic
functions.

\begin{theorem} \label{max-min} Let $u:[B(\cH)^n]_\gamma\to B(\cE)\otimes_{min}B(\cH)$ be a
self-adjoint  free pluriharmonic function with operator-valued  coefficients
satisfying either one of
the following conditions:
\begin{enumerate}
\item[(i)] $u(X_1,\ldots, X_n)\leq u(0)$ for any $(X_1,\ldots,
X_n)\in [B(\cH)^n]_\gamma$;
\item[(ii)] $u(X_1,\ldots, X_n)\geq u(0)$ for any $(X_1,\ldots,
X_n)\in [B(\cH)^n]_\gamma$;
\item[(iii)] $u(rS_1,\ldots, rS_n)\leq u(t S_1,\ldots, t
S_n)$ for some $r, t\in [0,\gamma)$, $r\neq t$.
\end{enumerate}
Then $u=u(0)$.
\end{theorem}
\begin{proof}
Assume that  condition (i) holds. Since $v:=u-u(0)$ is a positive
free pluriharmonic function on $[B(\cH)^n]_\gamma$, we can apply
Theorem \ref{Harnack1} and deduce that $v=0$. If (ii) holds, the
proof is similar.

Finally, if we assume that (iii) holds, then $w(X_1,\ldots,X_n):=u(t
X_1,\ldots, t X_n)-u(r X_1,\ldots, r X_n)$ is a positive free
pluriharmonic function on $[B(\cH)^n]_1$ with $w(0)=0$. Applying
again Theorem \ref{Harnack1}, we deduce that $w(X_1,\ldots,X_n)=0$
for $(X_1,\ldots,X_n)\in [B(\cH)^n]_1$. Hence $w(\gamma S_1,\ldots,
\gamma S_n)=0$ for $0<\gamma<1$, which implies $u=u(0)$.
\end{proof}

 In  \cite{Po-holomorphic},
 we  developed  a  {\it free holomorphic functional calculus}. Using those
 ideas,
   we can similarly prove   that
 if $$f=\sum_{k=1}^\infty \sum_{|\alpha|=k} B_{(\alpha)}\otimes  Z_\alpha^* +
  A_{(0)}\otimes  I+\sum_{k=1}^\infty \sum_{|\alpha|=k} A_{(\alpha)}\otimes  Z_\alpha$$
  is a free  pluriharmonic  function on  $[B(\cH)^n]_\gamma$ with
  coefficients in $B(\cE)$
   and  $(T_1,\ldots, T_n)\in B(\cH)^n$ is any $n$-tuple of operators
   with joint spectral radius
$r(T_1,\ldots, T_n)<\gamma$, then $ f(T_1,\ldots, T_n) $ is  a
bounded linear operator, where the corresponding series converge in
norm. This provides a {\it free
 pluriharmonic  functional calculus}, which turns out to be continuous
as a map from the complete metric space $(Har(B(\cH)^n_\gamma),
\rho)$ to $B(\cE\otimes\cH)$ with the operator norm topology. Since
the proof of the next result is similar to  the proof of Theorem 5.8
from \cite{Po-holomorphic}, we shall omit it. We denote by
$Har_\CC(B(\cH)^n_\gamma)$ the set of all free pluriharmonic
functions on $[B(\cH)^n]_1$ with scalar coefficients.

\begin{theorem}
 If $T:=(T_1,\ldots, T_n)\in B(\cH)^n$ is any $n$-tuple of operators with
   joint spectral radius
$r(T_1,\ldots, T_n)<\gamma$ then the mapping $\Phi_T:
Har_\CC(B(\cH)^n_\gamma) \to B(\cH)$ defined by
 $$
\Phi_T(u):=u(T_1,\ldots, T_n)
$$
is a continuous   linear map such that $\Phi_T(X_\alpha)= T_\alpha$
and   $\Phi_T(X_\alpha^*)= T_\alpha^*$ for any $\alpha\in \FF_n^+$.
Moreover,  $\Phi_T$ is uniquely determined by these conditions.
\end{theorem}

\section{Bounded free pluriharmonic functions}

In this section we characterize the set of all bounded free
 pluriharmonic functions on the noncommutative unit ball $[B(\cH)^n]_1$ and
  obtain  a Fatou
 type result concerning their boundary  behavior.

 A
function $u:[B(\cH)^n]_1\to B(\cE\otimes \cH)$ is called bounded if
$$\|u\|:=\sup_{(X_1,\ldots, X_n)\in [B(\cH)^n]_1}||u(X_1,\ldots,
X_n)\|<\infty.
$$
We say that a free pluriharmonic function is bounded if its
representation on any Hilbert space is bounded. As we will see  in
the next result, it is enough to assume that the Hilbert space is
separable and infinite dimensional.

 Let $\cH$ be an  infinite dimensional Hilbert space
and denote by  $Har^\infty(B(\cH)^n_1)$  the set of all bounded free
pluriharmonic functions on  $[B(\cH)^n]_1$ with coefficients in
$B(\cE)$.
 For each $m=1,2,\ldots$,
we define the norms $\|\cdot
\|_m:M_m\left(Har^\infty(B(\cH)^n_1)\right)\to [0,\infty)$ by
setting
$$
\|[u_{ij}]_m\|_m:= \sup \|[u_{ij}(X_1,\ldots, X_n)]_m\|,
$$
where the supremum is taken over all $n$-tuples $(X_1,\ldots,
X_n)\in [B(\cH)^n]_1$. It is easy to see that the norms
$\|\cdot\|_m$, $m=1,2,\ldots$, determine  an operator space
structure  on $Har^\infty(B(\cH)^n_1)$,
 in the sense of Ruan (see e.g. \cite{ER}).

The main result of this section is the following characterization of
bounded  free pluriharmonic  functions.
 \begin{theorem}\label{bounded}
  If $u:[B(\cH)^n]_1\to B(\cE)\otimes_{min}B(\cH)$, then
 the following statements are equivalent:
\begin{enumerate}
\item[(i)] $u$ is a bounded free pluriharmonic function on
$[B(\cH)^n]_1$;
\item[(ii)]
there exists $f\in \overline{\cA_n(\cE)^*+\cA_n(\cE)}^{SOT}$ such
that \ $u(X)={\bf P}_X[f]$ for $X\in [B(\cH)^n]_1$.
\end{enumerate}
In this case,
  $f=\text{\rm SOT-}\lim\limits_{r\to 1}u(rS_1,\ldots, rS_n).
  $
   Moreover, the map
$$
\Phi:Har^\infty((B(\cH)^n_1)\to \overline{\cA_n(\cE)^*+\cA_n(\cE)}^{SOT}\quad
\text{ defined by } \quad \Phi(u):=f
$$ is a completely   isometric isomorphism of operator spaces, where
$\cA_n(\cE):=B(\cE)\otimes_{min} \cA_n$ and $\cA_n$ is
 the noncommutative disc algebra.
\end{theorem}

\begin{proof}
Assume   that $u$ is a bounded free pluriharmonic function on
 operatorial unit ball    and let
$$u(X_1,\ldots, X_n):=\sum_{k=1}^\infty \sum_{|\alpha|=k}
B_{(\alpha)}\otimes  X_\alpha^* +A_{(0)}\otimes I+ \sum_{k=1}^\infty
 \sum_{|\alpha|=k} {A_{(\alpha)}}\otimes  X_\alpha
$$
be its representation on the Hilbert space $\cH$.
  According to
Proposition \ref{plu-ha},  we deduce  that, for any $r\in[0,1)$, $
u(rS_1,\ldots, rS_n) $
 is in $\cA_n(\cE)^*+\cA_n(\cE)$.
One can show that $u$ is bounded if and only if  $\sup\limits_{0\leq
r<1}\|u(rS_1,\ldots, rS_n)\|<\infty$.
Indeed, if $(X_1,\ldots, X_n)\in [B(\cH)^n]_1$, then there exists
$r\in (0,1)$ such that $(\frac{1}{r}X_1,\ldots, \frac{1}{r}X_n)\in
[B(\cH)^n]_1$. Since $u(rS_1,\ldots, rS_n)\in \cA_n(\cE)^*
+\cA_n(\cE)$, the noncommutative von Neumann inequality
\cite{Po-von} implies $ \|u(X_1,\ldots, X_n)\| \leq \|u(rS_1,\ldots,
rS_n)\|. $ Hence, we deduce that $ \|u\|\leq \sup\limits_{0\leq r<1}
\|u(rS_1,\ldots, rS_n)\|<\infty.$
 Since $\cH$ is infinite dimensional, the reverse inequality
is obvious, therefore,
\begin{equation}\label{norm-u}
\|u\|= \sup\limits_{0\leq r<1} \|u(rS_1,\ldots, rS_n)\|.
\end{equation}
  Now, due to Theorem \ref{Toeplitz},
$ u(S_1,\ldots, S_n) $ is the Fourier representation of a
multi-Toeplitz operator $f\in B(\cE\otimes F^2(H_n))$ and
\begin{equation} \label{soso} f=\text{\rm
SOT-}\lim\limits_{r\to1}u(rS_1,\ldots, rS_n).
\end{equation}
Hence, we deduce that $f\in \overline{\cA_n(\cE)^*+\cA_n(\cE)}^{SOT}$. The
next step is to prove that $u(X)={\bf P}_X[f]$ for $X\in [B(\cH)^n]_1$.
Since $u_r(S_1,\ldots,S_n):=u(rS_1,\ldots, rS_n)$ is in
$\cA_n(\cE)^*+\cA_n(\cE)$, we can use the properties of the noncommutative
Poisson transform and deduce  (first on polynomials of the form
$\sum_{|\alpha|\leq m} C_{(\alpha)}\otimes S_\alpha$) that
\begin{equation}\label{cube}
u_r(X_1,\ldots, X_n)={\bf P}_X[u_r(S_1,\ldots, S_n)]
=(I_\cE \otimes K_X^*)[u_r(S_1,\ldots, S_n)\otimes I_\cH](I_\cE\otimes K_X )
\end{equation}
for any $X:=(X_1,\ldots, X_n)\in [B(\cH)^n]_1$.
 Taking into account that the map $Y\mapsto  Y \otimes I$ is
SOT-continuous on bounded subsets of $B(\cE\otimes F^2(H_n))$ and
$\sup\limits_{0\leq r<1}\|u_r(S_1,\ldots, rS_n)\|<\infty$, we can
use
 relations \eqref{soso} and \eqref{cube}  to obtain
\begin{equation}\label{so-kx}
\text{\rm SOT-}\lim_{r\to 1} u_r(X_1,\ldots, X_n)={\bf P}_X[f].
\end{equation}
On the other hand, taking into account that $u$ is a free
pluriharmonic function on $[B(\cH)^n]_1$, and hence continuous,
  we have
$\text{\rm SOT-}\lim_{r\to1}u_r(X_1,\ldots,X_n) =u(X_1,\ldots,X_n)$.
Therefore, $(i)\implies (ii)$.

To prove that $(ii)\implies (i)$, let $f\in
\overline{\cA_n(\cE)^*+\cA_n(\cE)}^{SOT}$ and define
$u:[B(\cH)^n]_1\to B(\cE)\otimes_{min}B(\cH)$  by  setting  \
$u(X):={\bf P}_X[f]$.  We show first that $u$ is a pluriharmonic
function.  Notice that due to Corollary \ref{multi}, $f$ is a
multi-Toeplitz operator. Let $ \varphi(S_1,\ldots,
S_n):=\sum_{k=1}^\infty \sum_{|\alpha|=k} B_{(\alpha)}\otimes
S_\alpha^*+ A_{(0)}\otimes I+\sum_{k=1}^\infty \sum_{|\alpha|=k}
{A_{(\alpha)}} \otimes S_\alpha $ be its Fourier representation.
According to Theorem \ref{Toeplitz}, $ \varphi(rS_1,\ldots, rS_n) $
is in   $\cA_n(\cE)^*+ \cA_n(\cE)$ for any $r\in [0,1)$,
\begin{equation}\label{sup}
\sup_{0\leq r< 1}\|\varphi(rS_1,\ldots, rS_n)\|= \|f\|,\quad \text{
and } \ f=\text{\rm SOT-}\lim_{r\to 1}\varphi(rS_1,\ldots, rS_n).
\end{equation}
Consequently, due to Proposition \ref{plu-ha}, the map
$g:[B(\cH)^n]_1\to B(\cE)\otimes_{min}B(\cH)$ defined by
$$g(X_1,\ldots, X_n):=\sum_{k=1}^\infty \sum_{|\alpha|=k} B_{(\alpha)}\otimes
X_\alpha^* +A_{(0)}\otimes  I+ \sum_{k=1}^\infty
 \sum_{|\alpha|=k} {A_{(\alpha)}}\otimes  X_\alpha
$$ is a free pluriharmonic function. Now, let us show that $u=g$.
Since  $\varphi_r(S_1,\ldots, S_n)$ is in $\cA_n(\cE)^*+\cA_n(\cE)$ for any
$r\in [0,1)$, we have
$$
\varphi_r(X_1,\ldots, X_n)={\bf P}_X[\varphi_r(S_1,\ldots, S_n)]
=(I_\cE\otimes K_X^*)[\varphi_r(S_1,\ldots, S_n)\otimes I_\cH](I_\cE\otimes K_X).
$$
As above, since the map $Y\mapsto  Y\otimes I $ is SOT-continuous on
bounded subsets of $B(\cE\otimes F^2(H_n))$ and  using relation
\eqref{sup}, we deduce that
\begin{equation*}
\text{\rm SOT-}\lim_{r\to 1} \varphi_r(X_1,\ldots,
X_n)=(I_\cE\otimes K_X^*)[f\otimes I_\cH](I_\cE\otimes K_X)=u(X_1,\ldots, X_n).
\end{equation*}
On the other hand, since $\varphi_r(X_1,\ldots, X_n)=g(rX_1,\ldots,
rX_n)$ and due to the continuity of $g$ on  $[B(\cH)^n]_1$,  we have
$\text{\rm SOT-}\lim_{r\to1}\varphi_r(X_1,\ldots,X_n)
=g(X_1,\ldots,X_n)$. Hence, and using  relation \eqref{so-kx}, we
obtain $u=g$. This completes the proof of the implication
$(ii)\implies (i)$.

To prove the last part of the theorem, notice that if $[u_{ij}]_m\in
M_m(Har^\infty(B(\cH)^n_1)$ then, as above (see relation
\eqref{norm-u}), one can show that $ \|[u_{ij}]_m\|_m=\sup_{0\leq
r<1}\|[u_{ij}(rS_1,\ldots, rS_n)]_m\| $ and  that the operators
$f_{ij}:=\text{\rm SOT-}\lim\limits_{r\to 1}u_{ij}(rS_1,\ldots,
rS_n)$  are multi-Toeplitz   for $i,j=1,\ldots,m$.
According to relation \eqref{PHU}, we have
$$
P_{\cE\otimes F^2(H_n)}[(I_\cE\otimes U^*)(f_{ij}\otimes
I_{\cD_*})(I_\cE\otimes U)]|_{\cE \otimes F^2(H_n)}=u_{ij}(rS_1,\ldots, rS_n),\quad 0\leq r<1.
$$
Hence, we deduce that $ \sup_{0\leq r<1}\|[u_{ij}(rS_1,\ldots,
rS_n)]_m\|\leq \|[f_{ij}]_m\|.$ Moreover, since $[f_{ij}]_m$ is
equal to $\text{\rm SOT-}\lim_{r\to 1}[u_{ij}(rS_1,\ldots,
rS_n)]_m$, we have equality in the above inequality. Therefore, the
map $\Phi$ is a completely isometric isomorphism of operator spaces.
The proof is complete.
\end{proof}

\begin{corollary}\label{bound-equiv}
If $u:[B(\cH)^n]_1\to B(\cE)\otimes_{min} B(\cH)$ then the following statements are
equivalent:
\begin{enumerate}
\item[(i)]
$u$ is a bounded free pluriharmonic function;
\item[(ii)] there is a bounded  function $\varphi:[0,1)\to \overline
{\cP_n(\cE)^*+ \cP_n(\cE)}^{\|\cdot \|}$ such that
$$
\varphi(r)={\bf P}_{\frac{r}{t}S}[\varphi(t)]\quad \text{ for } \ 0\leq
r<t<1, $$
   and
$u(X_1,\ldots, X_n)={\bf P}_{\frac{1}{r} X}[\varphi(r)]$ for any $
X:=(X_1,\ldots, X_n)\in  [B(\cH)^n]_r$  and  $ r\in (0,1)$.
\end{enumerate}
Moreover, $u$ and $\varphi$ uniquely determine each other and
satisfy the equation $u(rS_1,\ldots, rS_n)=\varphi(r)$ for  $r\in
[0,1)$.
\end{corollary}

\begin{proof}
Assume that $u$ is a bounded free pluriharmonic function on
$[B(\cH)^n]_1$. According  to Theorem \ref{bounded}, there exists
 $f\in \overline{\cA_n(\cE)^*+\cA_n(\cE)}^{SOT}$ such that \
$u(X)={\bf P}_X[f]$ for $X\in [B(\cH)^n]_1$ and $\sup\limits_{0\leq
r<1} \|u(rS_1,\ldots, rS_n)\|<\infty$. Setting
$\varphi(r):=u(rS_1,\ldots, rS_n)={\bf P}_{rS}[f]$ for $r\in [0,1)$,
we obtain a function $\varphi:[0,1)\to \overline {\cP_n(\cE)^*+
\cP_n(\cE)}^{\|\cdot \|}$. Notice that $ u(rS_1,\ldots, rS_n) ={\bf
P}_{\frac{r}{t}S}[u(tS_1,\ldots, tS_n)], $
 which implies $
\varphi(r)={\bf P}_{\frac{r}{t}S}[\varphi(t)]$  for $ 0\leq r<t<1$. On the
other hand, since $u$ is free pluriharmonic, Theorem \ref{MVP}
implies $ u(X_1,\ldots, X_n)={\bf P}_{\frac{1}{r} X}[\varphi(r)] $ for any
$X:=(X_1,\ldots, X_n)\in [B(\cH)^n]_r$ and $r\in (0,1)$. Therefore,
(ii) holds.

Conversely, assume that condition (ii) holds. Using again Theorem
\ref{MVP}, we deduce that the map $u:[B(\cH)^n]_1\to B(\cE)\otimes_{min} B(\cH)$
 defined
by $u(X_1,\ldots, X_n):= {\bf P}_{\frac{1}{r}X}[\varphi(r)]$ for any
$X:=(X_1,\ldots, X_n)\in [B(\cH)^n]_r$ and $r\in (0,1)$ is free
pluriharmonic.  Since $\varphi$ is bounded, the relation above
implies $\|u\|\leq \|\varphi\|_\infty$, which completes the proof.
\end{proof}

We can prove now the following Fatou type result concerning the
boundary behavior of    bounded free pluriharmonic functions. This
also extends the $F_n^\infty$-functional calculus for pure row
contractions \cite{Po-funct}. We recall that $(X_1,\ldots, X_n)\in
B(\cH)^n$ is a pure $n$-tuple if $\sum_{|\alpha|=k} X_\alpha
X_\alpha^*\to 0$, as $k\to\infty$, in the strong operator topology.

\begin{theorem}\label{Fatou1}
Let $u$  be a bounded free pluriharmonic function on $[B(\cH)^n]_1$
with coefficients in $B(\cE)$.
 If   $(X_1,\ldots, X_n)\in [B(\cH)^n]_1^-$ is a pure
$n$-tuple of operators, then $
 \text{\rm SOT-}\lim_{r\to 1} u(rX_1,\ldots,
rX_n) $ exists.
\end{theorem}
\begin{proof}
Since $X:=(X_1,\ldots, X_n)\in [B(\cH)^n]_1^-$ is a pure $n$-tuple
of operators, the Poisson kernel $K_X$ is an isometry and
$p(X_1,\ldots, X_n)=(I_\cE \otimes K_X^*)[p(S_1,\ldots,S_n)\otimes
I_\cH] (I_\cE\otimes K_X)$ for any polynomial $p(S_1,\ldots,S_n)\in
\cP_n(\cE)$.  Using the fact that $u(rS_1,\ldots, rS_n)\in
\cA_n(\cE)^*+ \cA_n(\cE)$ for $r\in [0,1)$, we deduce that
\begin{equation}
\label{UKKU}
 u(rX_1,\ldots, rX_n)=(I_\cE \otimes K_X^*)(u(rS_1,\ldots, rS_n)\otimes
I)(I_\cE \otimes K_X).
\end{equation}
 Due to the boundedness  of  $u$, Theorem
\ref{bounded} implies $\text{\rm SOT-}\lim\limits_{r\to 1}
u(rS_1,\ldots, rS_n)=f\in \overline{\cA_n(\cE)^*+\cA_n(\cE)}^{SOT}$
and $\sup_{r\in[0,1)} \|u(rS_1,\ldots, rS_n)\|<\infty$.
Now,
using relation \eqref{UKKU}, we deduce that $\text{\rm
SOT-}\lim\limits_{r\to 1} u(rX_1,\ldots, rX_n)$ exists, which
completes the proof.
\end{proof}

\bigskip

\section{
Dirichlet extension problem for free pluriharmonic functions}

In this section  we solve the Dirichlet extension problem for the
noncommutative ball $[B(\cH)^n]_1$ and obtain a version of the
maximum principle for free  pluriharmonic (resp. $C^*$-harmonic)
functions.

We denote by $Har^c((B(\cH)^n_1)$ the set of all
  free pluriharmonic functions on $[B(\cH)^n]_1$ with operator-valued coefficients, which
 have continuous extensions   (in the operator norm topology) to
the closed ball $[B(\cH)^n]_1^-$. Throughout this section we assume
that $\cH$ is  an  infinite dimensional Hilbert space.

\begin{theorem}\label{Dirichlet}  If $u:[B(\cH)^n]_1\to B(\cE)\otimes_{min} B( \cH)$, then
 the following statements are equivalent:
\begin{enumerate}
\item[(i)] $u$ is a free pluriharmonic function on $[B(\cH)^n]_1$ which
 has a continuous extension  (in the operator norm topology) to
the closed ball $[B(\cH)^n]_1^-$;
\item[(ii)]
there exists $f\in \overline{\cA_n(\cE)^*+\cA_n(\cE)}^{\|\cdot\|}$ such that
$u(X)={\bf P}_X(f)$ for $X\in [B(\cH)^n]_1$;
\item[(iii)] $u$ is a free pluriharmonic function on $[B(\cH)^n]_1$
such that \ $u(rS_1,\ldots, rS_n)$ converges in the operator norm
topology, as $r\to 1$.
\end{enumerate}
In this case, $  f=\lim\limits_{r\to 1}u(rS_1,\ldots, rS_n),$ where
the convergence is in the operator norm. Moreover, the map $
\Phi:Har^c(B(\cH)^n_1)\to
\overline{\cA_n(\cE)^*+\cA_n(\cE)}^{\|\cdot\|}\quad \text{ defined
by } \quad \Phi(u):=f $ is a  completely   isometric isomorphism of
operator spaces, where $\cA_n(\cE):=B(\cE)\otimes_{min} \cA_n$ and
$\cA_n$ is the
 noncommutative disc algebra.
\end{theorem}

\begin{proof}
First we prove that  $(iii)\implies (ii)$. Assume that $u$ is a free
pluriharmonic function on $[B(\cH)^n]_1$ with coefficients in
$B(\cE)$  such that $u(rS_1,\ldots, rS_n)$ converges in the operator
norm as $r\to 1$.  Using Proposition \ref{plu-ha}, we deduce that
$u(rS_1,\ldots, rS_n)\in \cA_n(\cE)^*+\cA_n(\cE)$ and, due to (iii),
there exists $f$ in $ \overline{\cA_n(\cE)^*+\cA_n(\cE)}^{\|\cdot\|}$
such that $u(rS_1,\ldots, rS_n)\to f$ in the operator norm topology
as $r\to 1$. We recall   that the noncommutative Poisson transform
${\bf P}_X$ is defined by
 $
{\bf P}_X[f]:=(I_\cE\otimes K_X^*)(f\otimes I_\cH)(I_\cE\otimes
K_X)$, \ $f\in B(\cE\otimes F^2(H_n)). $ On the other hand, since
$u(rS_1,\ldots, rS_n)\in \cA_n(\cE)^*+\cA_n(\cE)$, one can prove
(first on polynomials) that
$$
{\bf P}_X[u(rS_1,\ldots, rS_n)]=(I_\cE\otimes K_X^*)[u(rS_1,\ldots,
rS_n)\otimes I_\cH](I_\cE\otimes K_X)=u(rX_1,\ldots, rX_n)
$$
for any $r\in [0,1)$ and $X:=(X_1,\ldots, X_n)\in [B(\cH)^n]_1$.
Since $u(rS_1,\ldots, rS_n)\to f$ in the operator norm, we deduce
that $u(rX_1,\ldots, rX_n)\to{\bf P}_X[f]$, as $r\to 1$. Taking into
account that any free pluriharmonic function is continuous, we have
$u(rX_1,\ldots, rX_n)\to u(X_1,\ldots, X_n)$ in norm, as $r\to 1$.
Summing up the results above, we deduce that $u(X_1\ldots, X_n)={\bf
P}_X[f]$ for $X:=(X_1\ldots, X_n)\in [B(\cH)^n]_1$, which proves
(ii).

Now we prove the implication $(ii)\implies (iii)$. Assume that $f\in
\overline{\cA_n(\cE)^*+\cA_n(\cE)}^{\|\cdot\|}$.
 Due to Corollary \ref{multi}, $f$ is a multi-Toeplitz operator.
 Let $\sum_{|\alpha|\geq 1} B_{(\alpha)}\otimes S_\alpha^*
+ A_{(0)}\otimes I +\sum_{|\alpha|\geq 1} A_{(\alpha)}\otimes
S_\alpha
$
 be the  Fourier
representation of $\varphi(S_1,\ldots, S_n)$. By  Theorem
\ref{Toeplitz}, for each $r\in [0,1)$, the operator
 $
\varphi(rS_1,\ldots, rS_n) $ is in $\cA_n(\cE)^*+\cA_n(\cE)$. Now,
we prove that
\begin{equation}
\label{PrS}
 \varphi(rS_1,\ldots, rS_n)={\bf P}_{rS}[f],
 \end{equation}
  where
$rS:=(rS_1,\ldots, rS_n)$ and $0\leq r<1$. Let $\gamma,\beta\in
\FF_n^+$ be fixed and  $q:=\max\{|\beta|, |\gamma|\}$,  and define
$Q_{\gamma, \beta}:=\sum_{1\leq|\sigma|\leq q} B_{(\sigma)}\otimes
S_\sigma^*+A_{(0)}\otimes I+ \sum_{1\leq|\sigma|\leq q}
A_{(\sigma)}\otimes  S_\sigma$. Notice that
\begin{equation*}
\begin{split}
\left< {\bf P}_{rS}[f](x\otimes e_\gamma), y\otimes e_\beta  \right>
&= \left< (f\otimes
I_{F^2(H_n)})(I_\cE\otimes K_{rS})(x\otimes e_\gamma),(I_\cE\otimes  K_{rS}) (y\otimes e_\beta)\right>\\
&= \left< \sum_{\alpha\in \FF_n^+} f(x\otimes e_\alpha)\otimes
\Delta_{rS} r^{|\alpha|} S_\alpha^* e_\gamma, \sum_{\omega\in
\FF_n^+} y\otimes e_\omega\otimes \Delta_{rS} r^{|\omega|}
S_\omega^* e_\beta\right>
\\
&=\sum_{k=0}^\infty \sum_{|\alpha|=k}\left( \sum_{m=0}^\infty
\sum_{|\omega|=m} \left< f(x\otimes e_\alpha),y\otimes
e_\omega\right> \left<\Delta_{rS} r^{|\alpha|} S_\alpha^*
e_\gamma,\Delta_{rS}
r^{|\omega|} S_\omega^* e_\beta\right>\right)\\
&= \sum_{k=0}^q \sum_{|\alpha|=k}\left( \sum_{m=0}^q
\sum_{|\omega|=m} \left< Q_{\gamma, \beta}(x\otimes e_\alpha),
y\otimes e_\omega\right> \left<\Delta_{rS} r^{|\alpha|} S_\alpha^*
e_\gamma,\Delta_{rS}
r^{|\omega|} S_\omega^* e_\beta\right>\right)\\
&= \sum_{k=0}^\infty \sum_{|\alpha|=k}\left( \sum_{m=0}^\infty
\sum_{|\omega|=m} \left< Q_{\gamma, \beta}(x\otimes e_\alpha),
y\otimes e_\omega\right> \left<\Delta_{rS} r^{|\alpha|} S_\alpha^*
e_\gamma,\Delta_{rS}
r^{|\omega|} S_\omega^* e_\beta\right>\right)\\
&=
 \left< (Q_{\gamma, \beta}\otimes
I_{F^2(H_n)})(x\otimes K_{rS}(e_\gamma)), (y\otimes K_{rS} (e_\beta))\right>\\
&= \left< {\bf P}_{rS}[Q_{\gamma, \beta}](x\otimes e_\gamma), (y\otimes e_\beta)  \right>  \\
&= \left<\varphi(rS_1,\ldots, rS_n)(x\otimes  e_\gamma), y\otimes
e_\beta\right>
\end{split}
\end{equation*}
for any $x,y\in \cE$ and $\gamma, \beta\in \FF_n^+$. Consequently,
relation \eqref{PrS} holds.

For any $q\in \cP_n(\cE)^*+\cP_n(\cE)$, we have $q=\lim_{r\to 1}
{\bf P}_{rS}[q]$ in the operator norm topology.
  Since
$f\in \overline{\cA_n(\cE)^*+\cA_n(\cE)}^{\|\cdot\|}$ using a
standard approximation argument and the continuity in  the operator
norm of the noncommutative Poisson transform ${\bf P}_X $, we deduce
that $f=\lim_{r\to 1} {\bf P}_{rS}[f]$. Hence and using \eqref{PrS},
we
 have
$\varphi(rS_1,\ldots, rS_n)\to f$,  as $ r\to 1, $ in the operator
norm topology.
 Now,  define $u(X):={\bf P}_X[f]$ for $X\in
[B(\cH)^n]_1$ and note that
 Theorem \ref{bounded} (see also its proof) implies that
 $u$ is a free pluriharmonic function and $u(rS_1,\ldots,
 rS_n)=\varphi(rS_1,\ldots,
 rS_n)$. Therefore item (iii) follows.

Since $\cH$ is infinite dimensional, the implication $(i)\implies
(iii)$ is obvious. It remains to prove that (ii)$\implies$(i).
Assume that $f\in \overline{\cA_n(\cE)^*+\cA_n(\cE)}^{\|\cdot\|}$
and $u(X)={\bf P}_X[f]$ for $X\in [B(\cH)^n]_1$. Due to Theorem
\ref{bounded}, $u$ is a bounded free pluriharmonic function. As
above, one can show that
   for any $n$-tuple
$Y:=(Y_1,\ldots, Y_n)\in [B(\cH)^n]_1^-$, $ \tilde u(Y_1,\ldots,
Y_n):= \lim_{r\to 1}  {\bf P}_{rY}[f] $ exists in the operator norm
and, since $\|{\bf P}_{rY}[f]\|\leq \|f\|_\infty$, $ \|\tilde
u(Y_1,\ldots, Y_n)\|\leq \|f\|_\infty$   for any $ (Y_1,\ldots,
Y_n)\in [B(\cH)^n]_1^-. $
 Notice also that $\tilde u$ is an extension of the free pluriharmonic
  function $u$ defined by  $u(X):={\bf P}_X[f]$, $X\in  [B(\cH)^n]_1$.
Indeed, if $ (X_1,\ldots, X_n)\in [B(\cH)^n]_1$, then
\begin{equation*}
\begin{split}
\tilde u(X_1,\ldots, X_n)=\lim_{r\to 1} P_{rX}[f] =\lim_{r\to
1}u(rX_1,\ldots, rX_n)=u(X_1,\ldots, X_n).
\end{split}
\end{equation*}
The last equality is due to   the continuity of free pluriharmonic
functions.

Now, let us prove that $\tilde u:[B(\cH)^n]_1^-\to
B(\cE)\otimes_{min} B(\cH)$ is continuous. Since $u(rS_1,\ldots,
rS_n)$ converges in the operator norm to $f$, for any $\epsilon>0$
there exists $r_0\in [0,1)$ such that $\|f-u(r_0S_1,\ldots, r_0
S_n)\|<\epsilon$. Since $f-u(r_0S_1,\ldots, r_0 S_n)\in
\overline{\cA_n(\cE)^*+\cA_n(\cE)}^{\|\cdot\|}$, and using the
properties of the noncommutative Poisson transform, we deduce the
von Neumann type inequality
\begin{equation}
\label{tild-f} \|\tilde u(T_1,\ldots, T_n)-u(r_0T_1,\ldots,
r_0T_n)\|\leq \|f-u(r_0S_1,\ldots, r_0S_n) \|< \frac{\epsilon}{3}
\end{equation}
for any $(T_1,\ldots, T_n)\in [B(\cH)^n]_1^-$.  Since  $u$ is a
continuous function on $[B(\cH)^n]_1$, there exists $\delta>0$ such
that $ \|u(r_0T_1,\ldots, r_0T_n)-u(r_0Y_1,\ldots,
r_0Y_n)\|<\frac{\epsilon}{3} $
 for any $n$-tuples $(T_1,\ldots, T_n)$ and  $(Y_1,\ldots, Y_n)$   in $[B(\cH)^n]_1^-$
     such that  $\|(T_1-Y_1,\ldots, T_n-Y_n)\|<\delta$.
Hence, and using \eqref{tild-f}, we have
\begin{equation*}
\begin{split}
\|\tilde u(T_1,\ldots, T_n)-\tilde u(Y_1,\ldots, Y_n)\| &\leq
\|\tilde u(T_1,\ldots, T_n)-u(r_0T_1,\ldots, r_0T_n)\|\\
&\qquad
 + \| u(r_0T_1,\ldots, r_0T_n)- u(r_0Y_1,\ldots, r_0Y_n)\|\\
&\qquad + \|u(r_0Y_1,\ldots, r_0Y_n)-\tilde u(Y_1,\ldots, Y_n)\|
<\epsilon,
\end{split}
\end{equation*}
whenever $\|(T_1-Y_1,\ldots, T_n-Y_n)\|<\delta$. This proves the
continuity of $\tilde u$ on $[B(\cH)^n]_1^-$. The last part of the
theorem follows from Theorem \ref{bounded}. This completes the
proof.
\end{proof}

The proof of the next result is similar to the proof of Corollary
\ref{bound-equiv} but one has to use Theorem \ref{Dirichlet}. We
shall omit it.

\begin{corollary}
If $u:[B(\cH)^n]_1\to B(\cE)\otimes_{min}B(\cH)$ then the following
statements are equivalent:
\begin{enumerate}
\item[(i)]
$u$ is a   free pluriharmonic function which has continuous
extension to the closed ball $[B(\cH)^n]_1^-$;
\item[(ii)] there exists  a continuous map $\varphi:[0,1]\to
\overline{\cP_n(\cE)^* +\cP_n(\cE)}^{\|\cdot\|}$ such that
$$
\varphi(r)={\bf P}_{\frac{r}{t}S}[\varphi(t)]\quad \text{ for } \
0\leq r<t\leq1, $$
 and
$u(X_1,\ldots, X_n)={\bf P}_{\frac{1}{r} X}[\varphi(r)]$ for any $
X\in  [B(\cH)^n]_r$  and  $ r\in (0,1]$.
\end{enumerate}
Moreover, $u$ and $\varphi$ uniquely determine each other and
satisfy the  equations $u(rS_1,\ldots, rS_n)=\varphi(r)$ if $r\in
[0,1)$ and $\varphi(1)=\lim\limits_{r\to 1} u(rS_1,\ldots, rS_n)$,
where the convergence is in the operator norm topology.
\end{corollary}

In what follows we introduce  the   class of  {\it  $C^*$-harmonic
functions} on the noncommutative ball $[B(\cH)]_1$.
 Let $\varphi:[0,1)\to
B(\cE)\otimes_{min}C^*(S_1,\ldots, S_n)$  be a map with the property
that
\begin{equation}
\label{gen-funct} \varphi(r)={\bf
P}_{\frac{r}{t}S}[\varphi(t)]\qquad \text{ for } \ 0\leq r<t<1,
\end{equation}
and  define the function $u:[B(\cH)^n]_1\to
B(\cE)\otimes_{min}B(\cH)$ by setting
 \begin{equation}
 \label{def-har}
 u(X_1,\ldots, X_n):= {\bf P}_{\frac{1}{r}X}[\varphi(r)]
\end{equation}
 for any  $ X:=(X_1,\ldots, X_n)\in
 [B(\cH)^n]_r$  and  $ r\in (0,1)$.
Notice that $u$ is well-defined. Indeed, if $0< r<t<1$ and $ X \in
 [B(\cH)^n]_r$, then, using relation \eqref{gen-funct} and Lemma
 \ref{PPP}, we  have
 $$
{\bf P}_{\frac{1}{r}X}[\varphi(r)]=({\bf P}_{\frac{1}{r}X}\circ{\bf
P}_{\frac{r}{t}S})[\varphi(t)]={\bf P}_{\frac{1}{t}X}[\varphi(t)],
$$
which proves our assertion. The map  $u$ defined by \eqref{def-har}
is called {\it
  $C^*$-harmonic function} on $[B(\cH)^n]_1$   and $\varphi$ is
 called the
 generating function of $u$.

We remark that, according to Theorem \ref{MVP}, any free
pluriharmonic function is a  $C^*$-harmonic function, while the
converse is not true. For instance, consider the function
$u(X_1,\ldots, X_n):=\sum_{\alpha,\beta\in \Lambda}
A_{\alpha,\beta}\otimes X_\alpha X_\beta^*$,\quad $(X_1,\ldots,
X_n)\in
 [B(\cH)^n]_1
 $
 where $\Lambda$ is any finite subset  of $\FF_n^+$ and
 $A_{\alpha,\beta}\in B(\cE)$.

\begin{proposition}\label{har-prop} Let $\varphi:[0,1)\to
B(\cE)\otimes_{min}C^*(S_1,\ldots, S_n)$  be a map  satisfying
relation \eqref{gen-funct}  and let  $u$ be the $C^*$-harmonic
function generated by $\varphi$. Then the following statements hold:
\begin{enumerate}
\item[(i)]
$\varphi$ is continuous on $[0,1)$ and $\|\varphi(r)\|\leq
\|\varphi(t)\|$ for $0\leq r<t<1$;
\item[(ii)] $u$ is a bounded  $C^*$-harmonic function if and only if its generating function
 $\varphi$ is bounded.
\end{enumerate}
\end{proposition}
\begin{proof}
 Using the continuity of the noncommutative Berezin transform (see
 Theorem \ref{berezin} part (ii)) and the fact that
 $$\|\varphi(r_1)-\varphi(r_2)\|=\left\| {\bf
P}_{\frac{r_1}{t}S}[\varphi(t)]-{\bf
P}_{\frac{r_2}{t}S}[\varphi(t)]\right\|
$$
for $0\leq  r_1<r_2<t$, we deduce that  $\varphi$ is continuous. On
the other hand,   we have $\left\|{\bf
P}_{\frac{r}{t}S}[\varphi(t)]\right\|\leq \|\varphi(t)\|$ for $0\leq
r<t<1$, which proves the second part of $(i)$.

To prove (ii),  assume  that  $\varphi$ is bounded and
$\sup_{r\in[0,1)} \|\varphi(r)\|\leq M$ for some $M>0$. Then, using
relation  \eqref{def-har},    we deduce that
$$\|u(X)\|=\left\|{\bf P}_{\frac{1}{r}X}[\varphi(r)]\right\|\leq
\|\varphi(r)\|\leq M
$$
for any $ X\in
 [B(\cH)^n]_r$  and  $ r\in (0,1)$.  Conversely, assume that  $u$ is a bounded  $C^*$-harmonic
function.  In particular, if $\cH=F^2(H_n)$, $X=(rS_1,\ldots,
rS_n)$, and $r<t<1$,  then, due to relations \eqref{gen-funct} and
\eqref{def-har}, we have $ \varphi(r)={\bf
P}_{\frac{r}{t}S}[\varphi(t)]=u(rS_1,\ldots, rS_n). $ Hence,
$\varphi$ is bounded on the interval $[0,1)$. The proof is complete.
\end{proof}

  The following result is needed to solve the
  Dirichlet extension problem for  $C^*$-harmonic functions.

\begin{theorem} \label{Diri-C*}
Let $f\in B(\cE)\otimes_{min}C^*(S_1,\ldots, S_n)$ and define
$u:[B(\cH)^n]_1\to B(\cE)\otimes_{min}B(\cH)$ by
$$
u(X):=  {\bf P}_{X}[f],\qquad X\in [B(\cH)^n]_1.
$$
Then \begin{enumerate}
\item[(i)]
 $u$ has a continuous extension $\widetilde u$  to the closed ball
$[B(\cH)^n]_1^-$ and the map
$$\Phi:B(\cE)\otimes_{min}C^*(S_1,\ldots, S_n)\to
C\left([B(\cH)^n]_1^-, B(\cE\otimes \cH)\right), \quad
\Phi(f)=\widetilde u,$$ is a complete linear isometry;
\item[(ii)]
$u$ has the  Poisson mean value property, i.e.,
$$
u(X)={\bf P}_{\frac{1}{r} X} [u(rS_1,\ldots, rS_n)]\quad \text{ for
any
 } \ X\in [B(\cH)^n]_r \ \text{ and } \  r\in (0,1),
$$
and $u(rS_1,\ldots, rS_n)\in B(\cE)\otimes_{min}C^*(S_1,\ldots,
S_n)$ for $r\in [0,1)$;
\item[(iii)] $\sup\limits_{\|X\|=r_1}
\|u(X)\|\leq \sup\limits_{\|X\|=r_2} \|u(X)\|  \quad \text{ for } \
0\leq r_1\leq r_2< 1;$
\item[(iv)]  $f= \lim\limits_{r\to 1}
u(rS_1,\ldots, rS_n)$ in the operator norm topology and
 \begin{equation*}
\begin{split}
\sup\limits_{X\in [B(\cH)^n]_1} \|u(X)\|&= \sup\limits_{\|X\|=1}
\|\widetilde u(X)\|=\lim\limits_{r\to 1} \|u(rS_1,\ldots, rS_n)\|\\
&=\sup_{0\leq r< 1}\|u(rS_1,\ldots, rS_n)\|=\|f\|;
\end{split}
\end{equation*}
\item[(v)]
 $u$ is a free pluriharmonic function on $[B(\cH)^n]_1$ if and
only if $f\in \overline{\cA_n(\cE)^*+ \cA_n(\cE)}^{\|\cdot\|}$.
\end{enumerate}
\end{theorem}
\begin{proof}

Similarly  to the relation \eqref{Po-tran}, one  can define
$\widetilde u:[B(\cH)^n]_1^-\to B(\cE)\otimes_{min}B(\cH)$ by
$$
\tilde u(X):=\lim_{r\to 1} {\bf P}_{rX}[f],\qquad X\in
[B(\cH)^n]_1^-,
$$
where the limit exists in the operator norm topology. Moreover, we
have $\tilde u(X)=u(X)$ for $X\in [B(\cH)^n]_1$.
Now, we prove that $\widetilde u$ is continuous on $[B(\cH)^n]_1^-$.
 Let $\epsilon>0$ and $q$ be a polynomial of the form
\begin{equation}\label{pol}
q=q(S_1,\ldots,S_n):=\sum C_{\alpha,\beta}\otimes  S_\alpha
S_\beta^*, \quad C_{\alpha,\beta}\in B(\cE),
\end{equation}
  such that $\|f-q\|< \frac
{\epsilon}{3}$. Since $\widetilde u(X)-q(X)=\lim_{r\to 1} ({\bf
P}_{rX}[f]-{\bf P}_{rX}[q])$, we deduce that
\begin{equation}\label{ine-ine}
\|\widetilde u(X)-q(X)\|\leq \|f-q\|< \frac {\epsilon}{3}
\end{equation}
for any $X\in [B(\cH)^n]_1^-$.  Choose $\delta>0$ such that
$\|q(X)-q(Y)\|<\frac {\epsilon}{3}$ whenever $X,Y\in [B(\cH)^n]_1^-$
with $\|X-Y\|<\delta$. Consequently,
\begin{equation*}
\begin{split}
\|\widetilde u(X)-\widetilde u(Y)\|&\leq\|\widetilde
u(X)-q(X)\|+\|q(X)-q(Y)\| +\|q(Y)-\widetilde u(Y)\|\leq \epsilon
\end{split}
\end{equation*}
for any  $X,Y\in [B(\cH)^n]_1^-$ with $\|X-Y\|<\delta$. Therefore,
$\widetilde u$ is continuous.
 Taking into account that $\cH$ is infinite dimensional and using
 the noncommutative von Neumann
inequality, relation \eqref{ine-ine} implies
$$
\|f\|- \frac {\epsilon}{3}\leq \|q\|=\sup_{\|X\|\leq 1} \|q(X)\|\leq
\frac {\epsilon}{3} +\|\widetilde u(X)\|
$$
for any $X\in [B(\cH)^n]_1^-$. Hence,  $\|f\|\leq \sup_{\|X\|\leq 1}
\|\widetilde u(X)\|$. The reverse inequality is due to the fact that
$ \widetilde u(X):=\lim_{r\to 1} P_{rX}[f]$. Therefore,
\begin{equation}
\label{sup1} \|f\|= \sup_{\|X\|\leq 1} \|\widetilde
u(X)\|=\|\widetilde u\|.
\end{equation}
Similarly, one can prove that $\|[f_{ij}]_m\|=\|[\widetilde
u_{ij}]_m\|$ for any matrix $[f_{ij}]_m\in
M_m(B(\cE)\otimes_{min}C^*(S_1,\ldots, S_n))$, which proves that the
map $\Phi$ is a complete isometry. This completes the proof of part
(i).

Let us prove (ii).   If $r\in [0,1)$ then
$$u(rS_1,\ldots,
rS_n)=(I_\cE \otimes K_{rS})(f\otimes I_{F^2(H_n)})(I_\cE\otimes
K_{rS}), $$
 where $rS:=(rS_1,\ldots, rS_n)$.
Let $\{q_m(S_1,\ldots, S_n)\}_{m=1}^\infty$  be a sequence of
polynomials  of the form \eqref{pol}  such that $q_m(S_1,\ldots,
S_n)\to f$ in the operator norm topology. Since
$$(I_\cE \otimes
K_{rS}^*)[q_m(S_1,\ldots, S_n)\otimes I_{F^2(H_n)}] (I_\cE \otimes
K_{rS})=q_m(rS_1,\ldots, rS_n)$$ is    of the form \eqref{pol}, we
deduce that   $u(rS_1,\ldots, rS_n)\in B(\cE)\otimes_{min}
C^*(S_1,\ldots, S_n)$ for any $r\in [0,1)$, and
\begin{equation}
\label{KfKf} \lim_{r\to 1} (I_\cE\otimes K_{rS}^*)(f\otimes
I_{F^2(H_n)})(I_\cE\otimes K_{rS})=f,
\end{equation}
where the convergence is in the operator norm topology. Due to Lemma
\ref{PPP}, for any $X\in [B(\cH)^n]_r$ and $r\in (0,1)$, we have
$$
u(X)={\bf P}_X[f]=({\bf P}_{\frac{1}{r}X}\circ {\bf P}_{rS})[f]=
{\bf P}_{\frac{1}{r}X}[u(rS_1,\ldots, rS_n)],
$$
which proves (ii).

To prove (iii), let $0\leq r_1<r_2<1$ and set $r:=\frac{r_1}{r_2}$.
Let us prove now that
\begin{equation}
\label{Krs} (I_\cE\otimes K_{rS}^*)[u(r_2S_1,\ldots, r_2S_n)\otimes
I_{F^2(H_n)}] (I_\cE\otimes K_{rS})=u(r_1S_1,\ldots, r_1S_n).
\end{equation}
Indeed, notice that
\begin{equation*}
\begin{split}
&(I_\cE\otimes K_{rS}^*)[u(r_2S_1,\ldots, r_2S_n)\otimes
I_{F^2(H_n)}] (I_\cE\otimes K_{rS}) \\
&\qquad\qquad  = (I_\cE\otimes K_{rS}^*)\left[
(I_\cE\otimes K_{r_2S}^*)(f\otimes I_{F^2(H_n)})(I_\cE\otimes K_{r_2S})\otimes I_{F^2(H_n)}\right](I_\cE\otimes  K_{rS})\\
&\qquad\qquad= \lim_{m\to \infty}(I_\cE\otimes K_{rS}^*)\left[
(I_\cE\otimes K_{r_2S}^*)(q_m(S_1,\ldots, S_n)\otimes I_{F^2(H_n)})
(I_\cE\otimes K_{r_2S})\otimes I_{F^2(H_n)}\right](I_\cE\otimes  K_{rS})\\
&\qquad\qquad= \lim_{m\to \infty} (I_\cE\otimes K_{rS}^*)\left[
 q_m(r_2S_1,\ldots, r_2S_n)\otimes I_{F^2(H_n)}\right] (I_\cE\otimes K_{rS})
 = \lim_{m\to \infty}
 q_m(r_1S_1,\ldots, r_1S_n)\\
 &\qquad\qquad=\lim_{m\to \infty}(I_\cE\otimes K_{r_1S}^*)\left[
 q_m(S_1,\ldots, S_n)\otimes I\right](I_\cE\otimes  K_{r_1S}) =
(I_\cE\otimes K_{r_1S}^*)\left(
  f\otimes I\right)(I_\cE\otimes  K_{r_1S})\\
  &\qquad\qquad=u(r_1S_1,\ldots, r_1S_n),
\end{split}
\end{equation*}
which proves our assertion. Since the Poisson kernel $K_{rS}$ is an
isometry, relation \eqref{Krs} implies
\begin{equation}\label{r1r2}
\|u(r_1S_1,\ldots, r_1S_n)\|\leq \|u(r_2S_1,\ldots, r_2S_n)\|.
\end{equation}
Now, let $X\in B(\cH)^n$ be  such that $0<\|X\|=r<1$. For any
polynomial  of the form \eqref{pol}, we have
$$
(I_\cE\otimes K_X^*)(q\otimes I_\cH)(I_\cE\otimes K_X)=\lim_{t\to 1}
(I_\cE\otimes K^*_{\frac{t}{r} X})\left\{\left[(I_\cE\otimes
K_{rS}^*)(q\otimes I_{F^2(H_n)})(I_\cE\otimes K_{rS})\right]\otimes
I_{\cH}\right\}(I_\cE\otimes  K_{\frac{t}{r} X}).
$$
Since  $K_{\frac{t}{r} X}$ is an isometry, we obtain
$$
\|(I_\cE\otimes K_X^*)(q\otimes I_\cH)(I_\cE\otimes K_X)\|\leq
\|(I_\cE\otimes K_{rS}^*)(q\otimes I_{F^2(H_n)})(I_\cE\otimes
K_{rS})\|.
$$
 An approximation of  $f\in B(\cE)\otimes_{min}C^*(S_1,\ldots, S_n)$ with
polynomials  of the form \eqref{pol} leads to
$$
\|(I_\cE\otimes K_X^*)(f\otimes I_\cH)(I_\cE\otimes K_X)\|\leq
\|(I_\cE\otimes K_{rS}^*)(f\otimes I_{F^2(H_n)})(I_\cE\otimes
K_{rS})\|,
$$
 whence
$\sup\limits_{\|X\|=r} \|u(X)\|\leq \|u(rS_1,\ldots, rS_n)\|$. On
the other hand, since $\cH$ is infinite dimensional and
$\|(rS_1,\ldots, rS_n)\|=r$, it is clear that $\sup\limits_{\|X\|=r}
\|u(X)\|\geq \|u(rS_1,\ldots, rS_n)\|$.  Therefore, $
\sup\limits_{\|X\|=r} \|u(X)\|= \|u(rS_1,\ldots, rS_n)\| $ for any
$r\in [0,1)$. Combining this with  inequality \eqref{r1r2}, we
deduce item (iii).

To prove (iv), notice that,  due to \eqref{KfKf},
\begin{equation}
\label{norm-f} \|f\|=\lim_{r\to 1} \|u(rS_1,\ldots, rS_n)\|.
\end{equation}
Hence, using relations \eqref{sup1},  \eqref{norm-f}, and the fact
that $\cH$ is infinite dimensional,   we deduce that
 $$\lim_{r\to 1}\|u(rS_1,\ldots,
rS_n)\|=\sup_{0\leq r<1}\|u(rS_1,\ldots, rS_n)\|=\sup_{\|X\|\leq 1}
\|\widetilde u(X)\|= \|f\|.
$$
Now, due to  the continuity of $\widetilde u$, we deduce item (iv).
 Item (v) follows from (i) and Theorem \ref{Dirichlet}.
This completes the proof.
 \end{proof}

Now, we can solve the following   Dirichlet extension problem for
 $C^*$-harmonic functions on the noncommutative ball $[B(\cH)^n]_1$.

\begin{theorem}\label{pro-gen}
If  $u:[B(\cH)^n]_1\to B(\cE)\otimes_{min}B(\cH)$, then the
following statements are equivalent:
\begin{enumerate}
\item[(i)] $u$ is a  $C^*$-harmonic function which has a continuous
 extension to $[B(\cH)^n]_1^-$;
\item[(ii)]
there exists $g$ in $ B(\cE)\otimes_{min}C^*(S_1,\ldots, S_n)$ such
that
$$
u(Y)={\bf P}_Y[g] \quad \text{ for any } \ Y\in [B(\cH)^n]_1;
$$
\item[(iii)] there exists a continuous function $\varphi:[0,1]\to
 B(\cE)\otimes_{min}C^*(S_1,\ldots,
S_n)$ such that
$$
\varphi(r)={\bf P}_{\frac{r}{t} S}[\varphi(t)]\quad \text{ for } \
0\leq r<t\leq 1
$$
and $u$ has the  Poisson mean value property with respect to
$\varphi$, i.e.,
$$
u(X)={\bf P}_{\frac{1}{r}X}[\varphi(r)]\quad \text{ for } \ X\in
[B(\cH)^n]_r \ \text{ and } \   r\in (0,1].
$$
\end{enumerate}
\end{theorem}

\begin{proof}
The implication (ii)$\implies$(iii) follows from Theorem
\ref{Diri-C*}
 and Lemma \ref{PPP}, if we define  the mapping
  $\varphi:[0,1]\to B(\cE)\otimes_{min}C^*(S_1,\ldots,
S_n)$ by  setting $\varphi(r):=u(rS_1,\ldots, rS_n)={\bf P}_{rS}[f]$
for
 $r\in
[0,1)$ and $\varphi(1):=\lim_{r\to 1} u(rS_1,\ldots, rS_n)$.
Conversely, assume  that (iii) holds. Setting $g:=\varphi(1)$, we
have $\varphi(r)={\bf P}_{rS}[g]$ for $r\in [0,1)$ and, due to Lemma
\ref{PPP},
$$
u(X)={\bf P}_{\frac{1}{r} X} [\varphi(r)]={\bf P}_{\frac{1}{r}
X}({\bf P}_{rS} [g])={\bf P}_X[g]
$$
for any $X\in [B(\cH)^n]_r$   and  $  r\in (0,1)$. Therefore,
 item (ii) holds.

 The implication $(ii)\implies (i)$ follows from
Theorem \ref{Diri-C*} and the implication $(ii)\implies (iii)$. It
remains to prove that $(i)\implies (ii)$. To this end, assume that
(i) holds. Then  there exists a function $\varphi:[0,1)\to
 B(\cE)\otimes_{min}C^*(S_1,\ldots,
S_n)$ such that $ \varphi(r)={\bf P}_{\frac{r}{t} S}[\varphi(t)]$
for  $ 0\leq r<t < 1, $ and $ u(X):={\bf
P}_{\frac{1}{t}X}[\varphi(t)]$    for  $ X\in [B(\cH)^n]_t $ and $
t\in (0,1). $ Consequently, if $0\leq r<t$  and $X=(rS_1,\ldots,
rS_n)$, we deduce that
$$u(rS_1,\ldots, rS_n)={\bf P}_{\frac{r}{t}
S}[\varphi(t)]=\varphi(r).
$$
Since $u$ has a continuous extension to the closed ball
$[B(F^2(H_n)^n]_1^-$, we  deduce that $ g:=\lim_{r\to 1} \varphi(r)$
exists in the norm topology and  it is  in $
B(\cE)\otimes_{min}C^*(S_1,\ldots, S_n)$. Let $X\in [B(\cH)^n]_1$
and let $t\in (0,1)$ be such that $X\in [B(\cH)^n]_r$. Note that
\begin{equation*}
\begin{split}
\left\|u(X)-{\bf P}_X[g]\right\|&= \left\| {\bf
P}_{\frac{1}{r}X}[\varphi(r)]-{\bf P}_X[g]\right\| \leq \left\| {\bf
P}_{\frac{1}{r}X}[\varphi(r)-g]\right\|+
 \left\| {\bf P}_{\frac{1}{r}X}[g]- {\bf P}_{X}[g]\right\|\\
 &\leq \left\|\varphi(r)-g\right\|+\left\| {\bf P}_{\frac{1}{r}X}[g]- {\bf
 P}_{X}[g]\right\|.
\end{split}
\end{equation*}
Using the continuity of the noncommutative Berezin transform (see
 Theorem \ref{berezin} part (ii)) and taking $r\to 1$, we conclude that
$u(X)={\bf P}_X[g]$, which proves item (ii). The proof is complete.
\end{proof}

A  consequence of Theorem \ref{Diri-C*}  and Theorem \ref{pro-gen}
is the following version of the maximum principle for
$C^*$-harmonic functions.

\begin{corollary}\label{max-mod2}
Let $u$ be a  $C^*$-harmonic  function on $[B(\cH)^n]_1$ with
operator-valued coefficients  and let $r\in [0,1)$. Then
$$ \sup\limits_{\|X\|\leq r} \|u(X)\|=\sup\limits_{\|X\|= r} \|u(X)\|
=\|u(rS_1,\ldots, rS_n)\|.
$$
Moreover, if $0\leq r_1\leq r_2<1$, then \ $\sup\limits_{\|X\|= r_1}
\|u(X)\|\leq \sup\limits_{\|X\|= r_2} \|u(X)\|$.
\end{corollary}

\bigskip

 \section{Noncommutative  transforms: Fantappi\` e, Herglotz,
and Poisson}

In this section, we introduce noncommutative versions of Fantappi\` e, Herglotz,
and Poisson
transforms associated with completely bounded maps on the operator system
$\cR_n^* + \cR_n$ (or $B(F^2(H_n))$), where $\cR_n$ is the noncommutative
 disc algebra
generated by the right creation operators $R_1,\ldots, R_n$ on
 $F^2(H_n)$ and the identity.  These transforms are used to obtain  characterizations
for the set   of all free holomorphic functions with positive real
parts, and to study the geometric structure of the  free pluriharmonic functions  on $[B(\cH)^n]_1$.

Consider the operator space $\cR_n^*+ \cR_n$ and regard
$M_m(\cR_n^*+ \cR_n)$ as a subspace of $M_m(B(F^2(H_n)))$. Let
$M_m(\cR_n^*+ \cR_n)$ have the norm structure that it inherits from
the (unique) norm structure on the $C^*$-algebra $M_m(B(F^2(H_n)))$.
Let $\mu:\cR_n^*+ \cR_n\to B(\cE)$ be a completely  bounded linear
map. Then there exists a completely bounded linear map
$$\widetilde
\mu:=\mu\otimes I:(\cR_n^*+ \cR_n)\otimes_{min} B(\cH)\to
B(\cE)\otimes_{min} B(\cH)
$$
such that $ \widetilde \mu(f\otimes Y)= \mu(f)\otimes Y$ for $f\in
\cR_n^*+ \cR_n$ and  $Y\in B(\cH)$. Moreover, $\|\widetilde
\mu\|_{cb}=\|\mu\|_{cb}$ and, if $\mu$ is completely positive, then
 so is $\widetilde \mu$.

We introduce  the {\it noncommutative Fantappi\` e transform }  of a
completely  bounded linear  map   $\mu:\cR_n^*+ \cR_n\to
B(\cE)$ to be the map $\cF\mu: [B(\cH)^n]_1 \to
B(\cE)\otimes_{min}B(\cH)$ defined by
$$
(\cF\mu)(X_1,\ldots, X_n):= \widetilde\mu[(I-R_1^*\otimes X_1-\cdots
-R_n^*\otimes X_n)^{-1}]
$$
for $(X_1,\ldots, X_n)\in [B(\cH)^n]_1$. Notice that the
noncommutative Fantappi\` e transform  is a   linear map and
$\cF\mu$ is a free holomorphic function in the open unit ball
$[B(\cH)^n]_1$ with coefficients in $B(\cE)$.

If $\mu$ is a completely positive linear map  on the operator system
$\cR_n^* +\cR_n$, we define the {\it noncommutative Herglotz
transform} $H\mu: [B(\cH)^n]_1 \to B(\cE)\otimes_{min}B(\cH)$ by
$$
(H\mu)(X_1,\ldots, X_n):=\widetilde\mu[ 2(I-R_1^*\otimes X_1-\cdots
-R_n^*\otimes X_n)^{-1}-I].
$$

We  introduce  now the {\it noncommutative Poisson transform} of
a completely  bounded linear  map  $\mu:\cR_n^*+\cR_n\to B(\cE)$ to
be the map \ $\cP\mu : [B(\cH)^n]_1\to B(\cE)\otimes_{min}B(\cH)$
defined by
$$
(\cP \mu)(X_1,\ldots, X_n):=\widetilde \mu[P(R, X)],\quad
X:=(X_1,\ldots, X_n)\in [B(\cH)^n]_1,
$$
where the {\it free pluriharmonic Poisson kernel}  $P(R,X)$ is
defined by
\begin{equation}
\label{PRX}
 P(R,X):=\sum_{k=1}^\infty\sum_{|\alpha|=k} R_{\widetilde
\alpha}\otimes X_\alpha^* + I+\sum_{k=1}^\infty\sum_{|\alpha|=k}
R_{\widetilde \alpha}^*\otimes X_\alpha,
\end{equation}
where the convergence is in the operator norm topology of
$B(F^2(H_n)\otimes \cH)$ and $R:=(R_1,\ldots, R_n)$ is the $n$-tuple
of right creation operators.
 Due to the continuity of $\widetilde \mu$ in the operator
norm  and  Proposition \ref{plu-ha}, the Poisson transform $\cP\mu$
is a free pluriharmonic function on $[B(\cH)^n]_1$ with coefficients
in $B(\cE)$.

\begin{proposition}
\label{Poisson-factor} Let \ $\mu:\cR_n^*+\cR_n\to B(\cE)$ be a
completely bounded linear map. The following statements hold.

\begin{enumerate}
\item[(i)]
The map $X\mapsto P(R,X)$ is a positive pluriharmonic function  on
$[B(\cH)^n]_1$ with coefficients in $
 B(F^2(H_n))$ and has the factorization
$ P(X,R)=B_X ^* B_X$,  $ X\in [B(\cH)^n]_1, $ where $$B_X
:=(I\otimes \Delta_X)(I-R_1\otimes X_1^*-\cdots -R_n\otimes
X_n^*)^{-1}
$$
and $\Delta_X:=(I-X_1X_1^*-\cdots -X_nX_n^*)^{1/2}$. \item[(ii)] The
Poisson transform $\cP\mu$ coincides with the Berezin transform
$\cB_\mu(I,\cdot\,)$.
\item[(iii)]
 If $\mu$ is a positive
linear map, then $\cP\mu$ is positive on $[B(\cH)^n]_1$.
\end{enumerate}
\end{proposition}

\begin{proof}
Denote $R_X:=R_1\otimes X_1^*+\cdots +R_n\otimes X_n^*$. Since
$R_i^*R_j=\delta_{ij}I$, $i,j=1,\ldots, n$,  and using \eqref{PRX},
we have
\begin{equation*}
\begin{split}
P(R,X)&= (I-R_X)^{-1}-I+ (I-R_X^*)^{-1}\\
&=(I-R_X^*)^{-1}\left[ I-R_X-(I-R_X^*)(I-R_X)+ I-R_X^*
\right](I-R_X)^{-1}\\
&= (I-R_X^*)^{-1}\left[ I\otimes (I-X_1X_1^*-\cdots -X_nX_n^*)
\right](I-R_X)^{-1} =B_X^* B_X.
\end{split}
\end{equation*}
Consequently, $P(R,X)\geq 0$ for any $X\in [B(\cH)^n]_1$. Using now
the definition of the Berezin transform $\cB_\mu(I,\cdot\,)$, we
have
$$
\cB_\mu(I, X)=\widetilde\mu(B_X^* B_X)=\widetilde \mu(P(R,X))=(\cP
\mu)(X),
$$
which proves our assertion. Part (iii) is now obvious.
\end{proof}

We need a few notations. Assume that $\cH$ is an infinite
dimensional Hilbert space.  We denote by $M^+(B(\cH)^n_1)$ the set
of all free holomorphic functions on the noncommutative ball
$[B(\cH)^n]_1 $ with
coefficients in $B(\cE)$, that can be represented as noncommutative
Herglotz transforms of completely positive linear maps $\mu:\cR_n^*+
\cR_n\to B(\cE)$, up to a constant, that is, an operator of the form \ $i
\text{\rm Im}\, A\otimes I$, $A\in B(\cE)$.
   The set of
all free holomorphic functions $u$  on $[B(\cH)^n]_1$ with positive
real part, i.e.,  Re\,$u(X_1,\ldots,X_n)\geq 0$ for any
$(X_1,\ldots,X_n)\in [B(\cH)^n]_1$,    is denoted by
$Hol^+(B(\cH)^n_1)$.

An operator-valued  positive semidefinite kernel on the free
semigroup $\FF_n^+$ is a map
$ K:\FF_n^+\times\FF_n^+\to B(\cE) $ with the property that for each
$k\in\NN$, for each choice of vectors
 $h_1,\dots,h_k$ in $\cE$, and
$\sigma_1,\dots,\sigma_k$ in $\Sigma$ the inequality
$ \sum\limits_{i,j=1}^k\langle K(\sigma_i,\sigma_j)h_j,h_i\rangle\ge
0 $ holds. Such a kernel is called multi-Toeplitz if it has the
following properties:
$K(g_0,g_0)=I_\cE$, ($g_0$ is the neutral element in $\FF_n^+$)
and
$$
K(\sigma,\omega)=
\begin{cases} K(\sigma\backslash_l\omega,g_0) \ &\text{ if }\sigma>_l\omega  \\
K(g_0,\omega \backslash_l \sigma)\ &\text{ if } \omega>_l\sigma  \\
0 &\text { otherwise } \end{cases}
$$
(see Section 1 for notations).
 We denote by $\cS^+(B(\cH)^n_1)$ the
positive Schur class of free holomorphic functions $\phi$ on
$[B(\cH)^n]_1$ with coefficients in $B(\cE)$  such that the kernel
$K_\phi:\FF_n^+\times \FF_n^+\to B(\cE)$ defined by

\begin{equation}\label{ke1}
         K_{\phi}(\alpha, \beta):=
         \begin{cases}
          A^*_{(\widetilde{\beta\backslash_l \alpha})}
      &\text{ if } \beta>_l\alpha\\
         A_{(0)}+ A^*_{(0)}  &\text{ if } \alpha=\beta\\
          A_{(\widetilde{\alpha\backslash_l \beta})}
      &\text{ if } \alpha>_l\beta\\
          0\quad &\text{ otherwise},
         \end{cases}
         \end{equation}
is positive semi-definite, where $\tilde{\gamma}$ is the reverse of
$\gamma\in \FF_n^+$ and $\phi$ has the representation $
\phi(X_1,\ldots, X_n)=\sum_{k=0}^\infty \sum_{|\alpha|=k}
A_{(\alpha)}\otimes  X_\alpha. $

The main result of this section is the following.
\begin{theorem} \label{H=M=S}
$Hol^+(B(\cH)^n_1)=\cS^+(B(\cH)^n_1)=M^+(B(\cH)^n_1)$.
\end{theorem}
\begin{proof} Let $f$ be in
$Hol(B(\cH)^n_1)$ and have the representation $ f(X_1,\ldots,
X_n):=\sum_{k=0}^\infty \sum_{|\alpha|=k} A_{(\alpha)}\otimes
X_\alpha. $ For each $r\in [0,1)$, define the multi-Toeplitz
operator (with respect to $S_1,\ldots, S_n$)
\begin{equation}
\label{Ar}
 A_r:=\sum_{k=1}^\infty \sum_{|\alpha|=k}
          \frac{1}{2}A^*_{(\alpha)} \otimes r^{|\alpha|}R_\alpha^* +\frac{1}{2}
          (A_{(0)}+A^*_{(0)})\otimes I+
           \sum_{k=1}^\infty \sum_{|\alpha|=k}
          \frac{1}{2}A_{(\alpha)}\otimes  r^{|\alpha|} R_\alpha.
 \end{equation}
Due to the properties of the Poisson transform and the fact that
$\cH$ is infinite dimensional, we can prove  that  $\text{\rm Re}\, f
(X_1,\ldots, X_n)\geq 0$ on $[B(\cH)^n]_1$ if and only if $A_r \geq
0$ for any $r\in [0,1)$. Indeed, for each $X:=(X_1,\ldots,
X_n)\in[B(\cH)^n]_1$, let $r\in [0,1)$ such that $\|X\|<r$. Then,
due to Theorem \ref{MVP}, we have $\text{\rm Re}\,f(X_1,\ldots,
X_n)={\bf P}_{\frac{1}{r}X}[A_r].$ Consequently, if $A_r\geq 0$,
then $\text{\rm Re}\,f(X_1,\ldots, X_n)\geq 0$. The other
implication is obvious due to the fact that $\cH$ is infinite
dimensional.

 Now, we prove that
$Hol^+(B(\cH)^n_1)\subseteq \cS^+(B(\cH)^n_1)$. Assume that $f$ is
in $Hol^+(B(\cH)^n_1)$ and define, for each $r\in [0,1)$,  the
kernel $K_{f,r}:\FF_n^+\times \FF_n^+\to
      B(\cE)$
      by

\begin{equation}\label{ke}
         K_{f,r}(\alpha, \beta):=
         \begin{cases}
         \frac{1}{2}r^{|\beta\backslash_l \alpha|} A^*_{(\widetilde{\beta\backslash_l \alpha})}
      &\text{ if } \beta>_l\alpha\\
         \frac{1}{2}(A_{(0)}+ A^*_{(0)})  &\text{ if } \alpha=\beta\\
          \frac{1}{2}r^{|\alpha\backslash_l \beta|}A_{(\widetilde{\alpha\backslash_l \beta})}
      &\text{ if } \alpha>_l\beta\\
          0\quad &\text{ otherwise}.
         \end{cases}
         \end{equation}

       Note that if $\{h_\beta\}_{|\beta|\leq q}\subset \cE$, then
          \begin{equation*}
          \begin{split}
          &\left< \left(\sum_{k=1}^\infty \sum_{|\alpha|=k}
          A_{(\alpha)}\otimes  r^{|\alpha|}R_\alpha\right)
            \left(\sum_{|\beta|\leq q}  h_\beta\otimes e_\beta
      \right),
         \sum_{|\gamma|\leq q}  h_\gamma\otimes e_\gamma\right>\\
         &=
           \sum_{k=1}^\infty \sum_{|\alpha|=k}\left<\sum_{|\beta|\leq q}
          A_{(\alpha)}h_\beta \otimes r^{|\alpha|}R_\alpha e_\beta,
          \sum_{|\gamma|\leq q}  h_\gamma\otimes e_\gamma\right>
          =\sum_{\alpha\in \FF_n^+\backslash\{g_0\}}\sum_{|\beta|, |\gamma|\leq q}
           r^{|\alpha|}
          \left<A_{(\alpha)} h_\beta, {h}_\gamma\right>
          \left< e_{\beta \tilde{\alpha}}, e_\gamma\right>\\
          &=
          \sum_{ \gamma>\beta; ~|\beta|, |\gamma|\leq q}
          r^{|\gamma\backslash \beta|} \left<A_{(\widetilde{\gamma\backslash
      \beta})}
      h_\beta, {h}_\gamma\right> =
           \sum_{\gamma>\beta; ~|\beta|, |\gamma|\leq q} 2\left<K_{f,r}
            (\gamma, \beta)h_\beta, {h}_\gamma\right>.
          \end{split}
          \end{equation*}
Similar calculations reveal that
$$
\left< \left(\sum_{k=1}^\infty \sum_{|\alpha|=k}
          A^*_{(\alpha)}\otimes r^{|\alpha|} R_\alpha^*\right)
           \right. \left.\left(\sum_{|\beta|\leq q}  h_\beta\otimes e_\beta
      \right),
         \sum_{|\gamma|\leq q}  h_\gamma\otimes e_\gamma \right>
         =
         \sum_{\beta>\gamma; ~|\beta|, |\gamma|\leq q}2\left< K_{f,r}
            (\gamma, \beta)h_\beta, {h}_\gamma\right>.
            $$
Therefore, we obtain
$$
\left< A_r\left(\sum_{|\beta|\leq q}  h_\beta\otimes e_\beta
      \right),
         \sum_{|\gamma|\leq q}  h_\gamma\otimes e_\gamma \right>
         =
         \sum_{|\beta|, |\gamma|\leq q} \left<K_{f,r}
            (\gamma, \beta)h_\beta, {h}_\gamma\right>,
            $$
where the operator $A_r$ is defined by relation \eqref{Ar} and
$K_{f,r}$ is defined by \eqref{ke}. Since $A_r\geq 0 $ for
$r\in[0,1)$,
  we deduce that
             $\left[K_{f,r}
      (\alpha, \beta)\right]_{|\alpha|,|\beta|\leq q}\geq 0$ for
      any
              $r\in [0,1)$.
          Taking $r\to 1$, we obtain
          $\left[K_{f,1}
      (\alpha, \beta)\right]_{|\alpha|,|\beta|\leq q}\geq 0$.
      Therefore, $f\in \cS^+(B(\cH)^n_1)$.

      Now, we prove that $\cS^+(B(\cH)^n_1)\subseteq
      M^+(B(\cH)^n_1)$. Assume that $f\in \cS^+(B(\cH)^n_1)$. Since
      the kernel $K_{f,1}$ is positive semidefinite, so is the kernel $K_{\epsilon,
      f,1}$, $\epsilon>0$,  defined by
      $$
      K_{\epsilon, f,1}(\alpha, \beta):=(D_0+\epsilon I)^{-1/2}
      \left[ K_{f,1}(\alpha,\beta)+ \epsilon \delta_{\alpha \beta}  I
      \right]
      (D_0+\epsilon I)^{-1/2}, \quad  \alpha, \beta\in \FF_n^+,
      $$
      where $D_0:=\frac{1}{2}(A_{(0)}+  A^*_{(0)})$.
Now, since $K_{\epsilon, f,1}$ is a positive semidefinite multi-Toeplitz
kernel which is normalized, i.e., $K_{\epsilon, f,1}(g_0,g_0)=I$,
       we can
      apply
 Theorem 3.1 of \cite{Po-posi}  and   deduce that
there is a completely
      positive
      linear  map
          $\mu_\epsilon:C^*(R_1,\ldots, R_n)\to B(\cE)$ such that
          $
          \mu_\epsilon (R_\alpha)=K_{\epsilon,f,1}(g_0, \alpha),\quad
           \alpha\in \FF_n^+.
          $
The  linear  map  $\nu_\epsilon :C^*(R_1,\ldots, R_n)\to B(\cE)$
defined by $ \nu_\epsilon(g)=(D_0+\epsilon I)^{1/2} \mu_\epsilon(g)
(D_0+\epsilon I)^{1/2}$ is  completely positive and has the property
 that
 $$
\nu_\epsilon (R_\alpha)=\frac{1}{2}  A^*_{(\widetilde \alpha)} \
\text{ if } \ |\alpha|\geq 1 \ \text{ and } \
\nu_\epsilon(I)=\frac{1}{2}(A_{(0)}+  A^*_{(0)}) +\epsilon I. $$
 Setting
$\mu(R_\alpha):=\frac{1}{2}  A^*_{(\widetilde \alpha)}$ if
$|\alpha|\geq 1,$ and $\mu(I)=\frac{1}{2}(A_{(0)}+ A^*_{(0)})$, one
can easily see that
  $$\nu_\epsilon(g)=\mu(g)+ \epsilon \left< g(1),
1\right> I\quad \text{  for  } \ g\in C^*(R_1,\ldots, R_n). $$ Since
$\nu_\epsilon(g)\to \mu(g)$, as $\epsilon \to 0$, we deduce that
$\mu$ is a completely positive linear map.
Now,  using the definition of the noncommutative  Herglotz
transform, we have
\begin{equation*}
\begin{split}
(H\mu)(X_1,\ldots, X_n)&=2(\cF\mu)(X_1,\ldots, X_n)-\mu(I)\otimes I
\\&=\sum_{k=0}^\infty \sum_{|\alpha|=k} A_{(\alpha)}\otimes  X_\alpha
+\left(\frac{ A^*_{(0)}-A_{(0)}}{2}\right)I\\
&=f(X_1,\ldots, X_n)-i(\text{\rm Im}\,  f(0))\otimes I.
\end{split}
\end{equation*}
 Therefore, $f\in M^+(B(\cH)^n_1)$.

Let us prove now that $M^+(B(\cH)^n_1)\subseteq Hol^+(B(\cH)^n_1)$.
To this end, assume that  $\varphi=H\mu$ for some completely
positive linear  map $\mu:C^*(R_1,\ldots, R_n)\to B(\cE)$. Notice
that
\begin{equation*}
\begin{split}
&\frac{1}{2}(\varphi(X_1,\ldots, X_n)+\varphi(X_1,\ldots, X_n)^*)\\
&= \widetilde \mu\left[(I-R_1\otimes X_1^*-\cdots -R_n\otimes
X_n^*)^{-1}-I+(I-R_1^*\otimes X_1-\cdots -R_n^*\otimes
X_n)^{-1}\right] =\widetilde\mu[P(R,X)],
\end{split}
\end{equation*}
where $P(R,X)$ is defined by \eqref{PRX}.  Applying Proposition
\ref{Poisson-factor}, we deduce that $\text{\rm Re}\, \varphi\geq
0$. This completes the proof.
\end{proof}

As a consequence of Theorem \ref{H=M=S} (see also the proof), we
obtain the following noncommutative Herglotz-Riesz representation
for free holomorphic functions  with positive real parts on the
noncommutative ball $[B(\cH)^n]_1$.

\begin{theorem}\label{H-R}
Let $f:[B(\cH)^n]_1 \to B(\cE)\otimes_{min} B(\cH)$ be a free
holomorphic function with $\text{\rm Re}\, f\geq 0$ on
$[B(\cH)^n]_1$. Then
$$
f(X_1,\ldots, X_n)=\widetilde\mu[ 2(I-R_1^*\otimes X_1-\cdots
-R_n^*\otimes X_n)^{-1}-I]+i(\text{\rm Im}\,  f(0))\otimes I
$$
for some completely positive linear  map $\mu$ on the Cuntz-Toeplitz
algebra $C^*(R_1,\ldots, R_n)$.
\end{theorem}
\begin{proof}
 Use the proof of   Theorem \ref{H=M=S} and  apply  Arveson's extension theorem
\cite{Ar}.
\end{proof}

The following result  is a Naimark (\cite{N}) type theorem
concerning the geometric structure of $Hol^+(B(\cH)^n_1)$.

\begin{theorem} \label{geo-struct}
A free holomorphic function $f$ on $[B(\cH)^n]_1$  with coefficients
in $B(\cE)$ has positive real part if and only if there exists an
$n$-tuple of isometries  $(V_1,\ldots, V_n)$ on a Hilbert space $\cK$,
 with orthogonal ranges,
  and a  bounded operator $W:\cE\to \cK$
such that
$$
f(X_1,\ldots, X_n)=  (W^*\otimes I) \left[ 2(I-V_1^*\otimes
X_1-\cdots -V_n^* \otimes X_n)^{-1}- I\right](W\otimes I) +
i(\text{\rm Im}\, f(0))\otimes I.
$$
\end{theorem}
\begin{proof}
According to Theorem \ref{H=M=S},   $f$ is in $Hol^+(B(\cH)^n_1)$ if
and only if  it has the representation
 \begin{equation}
 \label{HR}
f(X_1,\ldots, X_n)=\widetilde\mu[ 2(I-R_1^*\otimes X_1-\cdots
-R_n^*\otimes X_n)^{-1}-I]+i(\text{\rm Im}  f(0))\otimes I
\end{equation}
for some  completely positive linear  map  $\mu$ on the
Cuntz-Toeplitz algebra $C^*(R_1,\ldots, R_n)$ with values in
$B(\cE)$. On the other hand, due to Stinespring's representation
theorem (see \cite{St}), $\mu$ is a  completely positive linear
 map  on $C^*(R_1,\ldots, R_n)$  if and only if there is a
Hilbert space $\cK$,
         a $*$-representation $\pi: C^*(R_1,\ldots, R_n)\to B(\cK)$, and
         a   bounded operator $W:\cE\to \cK$ with $\|\mu(I)\|=\|W\|^2$    such that
         $
         \mu(g)= W^*\pi(g)W$,  \   $g\in C^*(R_1,\ldots, R_n).
        $
        Notice that  $V_i:=\pi(R_i)$, $i=1,\ldots, n$, are isometries
        with orthogonal ranges, and any $*$-representation of $C^*(R_1,\ldots, R_n)$
        is generated by $n$ isometries with orthogonal ranges.
Hence and using \eqref{HR}, we find  the  required form for $f$.
This completes the proof.
\end{proof}

\begin{corollary}\label{positive-pluri}
The map $\mu\mapsto \cP\mu$ is a  linear and one-to-one
correspondence between the space of all completely positive linear
maps on the operator system $\cR_n^* + \cR_n$ and the space of all
positive free pluriharmonic functions on the open noncommutative
ball $[B(\cH)^n]_1$. In particular, any positive free pluriharmonic
function on $[B(\cH)^n]_1$  is the Poisson transform of   a
completely positive linear  map  on the Cuntz-Toeplitz algebra
$C^*(R_1,\ldots, R_n)$.
\end{corollary}
\begin{proof}
The first part follows from Theorem \ref{H=M=S} (see the proof) and
the fact that any positive  free  pluriharmonic function has the
form Re\,$f$ for some free holomorphic function $f$. The second part
is  also due to Theorem  \ref{H=M=S} and Arveson's extension
theorem.
\end{proof}

\begin{corollary}\label{pluri-struct}
A free pluriharmonic function $h$ on $[B(\cH)^n]_1$   with
coefficients in $B(\cE)$ is positive if and only if there exists an
$n$-tuple of isometries $(V_1,\ldots, V_n)$ on a Hilbert space
$\cK$, with orthogonal ranges,  and  a  bounded operator $W:\cE\to
\cK$ such that
$$
h(X_1,\ldots, X_n)=  (W^*\otimes I)\left[ B_X(V_1,\ldots, V_n)^*
B_X(V_1,\ldots, V_n) \right](W\otimes I),
$$
where $B_X(V_1,\ldots, V_n):=(I\otimes \Delta_X)(I-V_1\otimes
X_1^*-\cdots -V_n\otimes X_n^*)^{-1}. $
\end{corollary}
\begin{proof} By Corollary \ref{positive-pluri}, $h$ is a positive
free pluriharmonic function  if and only if there is
 a completely positive map $\mu: C^*(R_1,\ldots, R_n)\to B(\cE)$
 such that $h=\cP\mu$. Using Proposition \ref{Poisson-factor}, we
 deduce that
 $$h(X_1,\ldots, X_n)=\widetilde \mu(B_X^* B_X),
  \quad X:=(X_1,\ldots, X_n)\in [B(\cH)^n]_1,
 $$
  where $B_X$
 is defined by \eqref{KTR}.
 Now, the proof is similar to  the proof of Theorem
 \ref{geo-struct}.
\end{proof}

 Using Theorem \ref{H=M=S},  we can recast some of the results
from \cite{Po-analytic} and \cite{Po-unitary}    to our setting.
More precisely, we can deduce the following Fej\' er type factorization
result and
  Fej\' er and Egerv\' ary-Sz\' azs type inequalities (see \cite{ES}, \cite{Fe}) for the
 moments  of a positive linear functional on the Cuntz-Toeplitz
 algebra $C^*(R_1,\ldots, R_n)$.
\begin{theorem}\label{Fejer}
Let $\mu:C^*(R_1,\ldots, R_n)\to \CC$ be a positive linear
functional with $\mu(R_\alpha)= 0$ for any $\alpha\in \FF_n^+$,
$|\alpha|\geq m$. Then
\begin{enumerate}
\item[(i)]
 there exists a polynomial $p(S_1,\ldots,
S_n)$ in the noncommutative disc algebra $\cA_n$ such that
$$
(\cP \mu) (X)=P_X[p(S_1,\ldots, S_n)^* p(S_1,\ldots, S_n)\otimes
I_\cH], \quad X\in [B(\cH)^n]_1.
$$
\item[(ii)]
$$ \left(\sum_{|\alpha|=k}|\mu(R_\alpha)|^2\right)^{1/2}\leq
\mu(I)\cos \frac{\pi}{\left[\frac{m-1}{k}\right]+2} $$ for $1\leq
k\leq m-1$, where $[x]$ denotes the integer part of $x$.
\end{enumerate}
\end{theorem}
\begin{proof}
Since $\mu$ is a  positive linear functional, Corollary
\ref{positive-pluri} shows that the Poisson transform $\cP\mu$ is a
positive  free pluriharmonic function on $[B(\cH)^n]_1$. Taking into
account that $\mu(R_\alpha)= 0$ for any $\alpha\in \FF_n^+$,
$|\alpha|\geq m$, we have
$$
(\cP\mu)(X_1,\ldots, X_n)=\sum_{1\leq |\alpha|\leq m-1} \bar
a_\alpha X_\alpha^*+ a_0+\sum_{1\leq |\alpha|\leq m-1} a_\alpha
X_\alpha,
$$
where $a_0=\mu(I)$ and $ a_{\alpha}=\mu(R^*_{\widetilde \alpha})$
for $1\leq |\alpha|\leq m-1$. Notice  that
$q(S^*,S):=(\cP\mu)(S_1,\ldots, S_n)$ is a positive multi-Toeplitz
operator with respect to the right creation operators. The Fej\' er
type factorization theorem  from \cite{Po-analytic} implies the
existence of a polynomial $p(S_1,\ldots, S_n)$ such that
$$q(S^*,S)=p(S_1,\ldots, S_n)^*p(S_1,\ldots, S_n).
$$
 On the other hand, using the noncommutative Poisson transform $P_X$, we have
$$(\cP\mu)(X_1,\ldots, X_n)=P_X[q(S^*,S)\otimes I_\cH],
\quad X:=(X_1,\ldots, X_n)\in [B(\cH)^n]_1. $$
Combining these equalities we deduce part (i). Part (ii) follows
from Theorem 8.3 from \cite{Po-unitary} applied to the positive
multi-Toeplitz operator $q(S^*,S)$. The proof is complete.
\end{proof}

\bigskip

\section{The Banach space $Har^1(B(\cH)^n_1)$}

  In this section  we  characterize those  free pluriharmonic function   which
are noncommutative Poisson transforms of completely bounded linear
maps on the operator system $\cR_n^*+\cR_n$, and  those self-adjoint
free pluriharmonic functions which admit Jordan type decompositions.

Throughout this section, we assume that $\cE$ is a separable Hilbert
space. Let $h$ be a free pluriharmonic function on the
noncommutative ball $[B(\cH)^n]_1$ with operator-valued
coefficients in $B(\cE)$, and let $\tau:B(F^2(H_n))\to \CC$ be the
bounded  linear functional defined by $\tau(f)=\left<f(1),1\right>$.
 The radial function associated with $h$,
$$[0,1)\ni r\mapsto h(rR_1,\ldots, rR_n)\in \cR_n(\cE)^*+\cR_n(\cE),$$
generates a family  $\{\nu_{h,r}\}_{r\in [0, 1)}$ of completely  bounded
linear maps
$\nu_{h,r}: \cR_n^*+\cR_n\to  B(\cE)$
uniquely determined by the equations
\begin{equation*}
\begin{split}
\nu_{h,r}(R_{\widetilde\alpha}^*)&:= (\text{\rm id}\otimes
\tau)\left[(I\otimes R_\alpha^*)
 h(rR_1,\ldots, rR_n)\right],
\ \alpha\in \FF_n^+,\\
 \nu_{h,r}(R_{\widetilde\alpha})&:= (\text{\rm id}\otimes \tau)\left[ h(rR_1,\ldots,
rR_n)(I\otimes R_\alpha) \right], \ \alpha\in \FF_n^+\backslash
\{g_0\}.
\end{split}
\end{equation*}
Indeed, let $q$ be a polynomial of the form
\begin{equation}
\label{*p*}
 q:=\sum_{1\leq|\alpha|\leq m} c_\alpha R_{\widetilde
\alpha}^*+ \sum_{|\alpha|\leq m } d_\alpha R_{\widetilde \alpha},
\qquad c_\alpha, d_\alpha\in \CC,
\end{equation}
and notice that
\begin{equation*}
\begin{split}
\nu_{h,r}(q)&=(\text{\rm id}\otimes \tau)\left[\left(I\otimes
\sum_{1\leq |\alpha|\leq m} c_\alpha R_\alpha^*\right)
 h(rR_1,\ldots, rR_n)\right] +
 (\text{\rm id}\otimes \tau)\left[ h(rR_1,\ldots,
rR_n)\left(I\otimes \sum_{|\alpha|\leq m} d_\alpha R_\alpha\right)
\right]
\\
&= (\text{\rm id}\otimes \tau)\left[(I\otimes q)h(rR_1,\ldots,
rR_n)\right] - (\text{\rm id}\otimes \tau)[ d_0 h(rR_1,\ldots,
rR_n)]
  + (\text{\rm id}\otimes \tau)\left[ h(rR_1,\ldots, rR_n)
(I\otimes q)\right].
\end{split}
\end{equation*}
Hence, we deduce that
$$
\|\nu_{h,r}(q)\|\leq (2\|q\|+|d_0|) \|h(rR_1,\ldots, rR_n)\|
\leq 3\|q\|\|h(rR_1,\ldots, rR_n)\|.
$$
Since $\cR_n^*+ \cR_n$ is the norm closure of  polynomials of the
form \eqref{*p*} and passing to matrices over $\cR_n^*+\cR_n$, one
can deduce  that
$$\|[\nu_{h,r}(f_{ij})]_k\|\leq 3\|[f_{ij}]_k\|\|h(rR_1,\ldots, rR_n)\|\quad \text{ for any } \
[f_{ij}]_k\in M_{k}(\cR_n^*+ \cR_n), k\in\NN.
$$
Consequently, $\nu_{h,r}$ is a completely bounded  map for each $r\in [0,1)$.

\begin{lemma}\label{mu-r}
Let $\mu: \cR_n^*+ \cR_n\to B(\cE)$ be a completely bounded linear map. For each
$r\in [0,1]$, define  the linear map $\mu_r: \cR_n^*+ \cR_n\to B(\cE)$
by
$$
\mu_r(R_\alpha):= r^{|\alpha|}\mu(R_\alpha), \ \alpha\in \FF_n^+,
\quad \text{ and } \quad \mu_r(R_\alpha^*):=
r^{|\alpha|}\mu(R_\alpha^*), \ \alpha\in \FF_n^+\backslash \{g_0\}.
$$
Then
\begin{enumerate}
\item[(i)]
$\mu_r$ is a completely bounded linear map;
\item[(ii)] $\mu_r(f)\to \mu(f)$ in
the operator topology, as $r\to 1$;
\item[(iii)]
$ \|\mu\|_{cb}=\sup\limits_{0\leq r<1}
\|\mu_r\|_{cb}=\lim\limits_{r\to 1} \|\mu_r\|_{cb}. $
\end{enumerate}
\end{lemma}

\begin{proof}
Using the noncommutative von Neumann inequality, one can prove that
 $
\|\mu_{r_1}\|\leq \|\mu_{r_2}\|
 $
for  $0\leq r_1<r_2\leq 1$.
Indeed, if $p(R_1,\ldots, R_n)$ and $q(R_1,\ldots, R_n)$ are
polynomials in $\cR_n$, then we have
\begin{equation*}
\begin{split}
\|\mu_{r_1}(q(R_1,\ldots, R_n)^*+p(R_1,\ldots, R_n))\|&=
\|\mu(q(r_1R_1,\ldots, r_1R_n)^*+p(r_1R_1,\ldots, r_1R_n))\|\\
&=\left\|\mu_{r_2}\left(q\left(\frac{r_1}{r_2}R_1,\ldots,
\frac{r_1}{r_2}R_n\right)^*+p\left(\frac{r_1}{r_2}R_1,\ldots,
\frac{r_1}{r_2}R_n\right)\right)\right\|\\
&\leq \|\mu_{r_2}\| \|q(R_1,\ldots, R_n)^*+p(R_1,\ldots, R_n))\|,
\end{split}
\end{equation*}
which proves our assertion.
 In particular, we have $\|\mu_r\|\leq
\|\mu\|$ for any $r\in [0,1)$. Similarly, passing to matrices over
$\cR_n^*+\cR_n$, one can show that $\|\mu_{r_1}\|_{cb}\leq
\|\mu_{r_2}\|_{cb}$ if $0\leq r_1<r_2\leq 1$, and $\|\mu_r\|_{cb}\leq
\|\mu\|_{cb}$ for any $r\in [0,1)$. An approximation argument shows
that $\mu_r(A)\to\mu(A)$ in the operator norm topology, as $r\to 1$,
for any $A\in \cR_n^*+ \cR_n$. Now,
 one can easily see that
$\|\mu\|_{cb}=\sup_{0\leq r<1} \|\mu_r\|_{cb}$. Hence and using the fact that
 the function
$r\mapsto \|\mu_r\|_{cb}$ is increasing, we deduce that
 the limit  $\lim_{r\to 1} \|\mu_r\|_{cb}$
exists and it is equal to $\|\mu\|_{cb}$. This completes the proof.
\end{proof}

\begin{theorem}\label{pluri-measure}
Let $h$ be a free pluriharmonic function  on $[B(\cH)^n]_1$ with
operator-valued coefficients in $B(\cE)$. Then the following statements are equivalent:
\begin{enumerate}
\item[(i)]
there exists a completely  bounded linear map $\mu:C^*(R_1,\ldots,
R_n)\to B(\cE)$ such that $ h=\cP\mu; $
\item[(ii)]
the completely  bounded linear maps \ $\{\nu_{h,r}\}_{r\in[0,1)}$,
associate with the radial function of $h$,
are uniformly  bounded, i.e.,   $\sup\limits_{0\leq r<1}
\|\nu_{h,r}\|_{cb}<\infty$;
\item[(iii)]
there exist positive free pluriharmonic functions $h_1,h_2,h_3,h_4$
on $[B(\cH)^n]_1$ with coefficients in $B(\cE)$
 such that
$$h=(h_1-h_2)+i(h_3-h_4).
$$
\end{enumerate}
\end{theorem}

\begin{proof} Assume that (i) holds. Since $h$ and $\cP \mu$ are
free holomorphic functions  on $[B(\cH)^n]_1$, we deduce that $h$
has the representation
$$h(X_1,\ldots,
X_n)=\sum_{k=0}^\infty\sum_{|\alpha|=k} B_{(\alpha)}\otimes
X_\alpha^* + A_{(0)}\otimes  I+ \sum_{k=0}^\infty\sum_{|\alpha|=k}
A_{(\alpha)}\otimes  X_\alpha, $$
 where
 $A_{(\alpha)}:=\mu(R_{(\widetilde\alpha)}^*)$, \ $\alpha\in \FF_n^+$,  and
 $B_{(\alpha)}:=\mu(R_{(\widetilde\alpha)})$, \
 $\alpha\in \FF_n^+\backslash \{g_0\}$.
Notice that $\nu_{h, r}(R_\alpha)=r^{|\alpha|} \mu(R_\alpha)$,
$\alpha\in\FF_n^+$, and
  $\nu_{h, r}(R_\alpha^*)=r^{|\alpha|}\mu(R^*_\alpha)$,
  $\alpha\in \FF_n^+\backslash \{g_0\}$. Applying
Lemma \ref{mu-r} to $\{\nu_{h, r}\}$, we deduce item (ii).

Now, we prove the implication (ii)$\implies $ (i). To this end,
assume that $h$ is a free pluriharmonic function on $[B(\cH)^n]_1$ with
 coefficients
in $B(\cE)$ and condition  (ii) holds. Let $\{f_j\}$ be a countable
dense subset of $\cR_n^*+ \cR_n$ (for instance, consider all
``noncommutative trigonometric polynomials'' of the form
$\sum_{|\alpha|\leq m} c_\alpha R_\alpha^* + \sum_{|\alpha|\leq m}
d_\alpha R_\alpha$,  whose coefficients lie in some countable dense
subset of the complex plane). For each $j$, we have $\|\nu_{h,
r}(f_j)\|\leq M\|f_j\|$ for any $r\in[0,1)$, where $M:=\sup_{0\leq
r<1} \|\nu_{h, r}\|_{cb}$.

Due to Banach-Alaoglu theorem,  the ball $[B(\cE)]_M^-$  is compact
 in the $w^*$-topology. Since $\cE$ is a separable Hilbert space,
  $[B(\cE)]_M^-$ is  a metric space in the $w^*$-topology which coincides with
   the weak operator topology on $[B(\cE)]_M^-$. Consequently,
 the diagonal process guarantees the
existence  of a sequence $\{r_m\}_{m=1}^\infty$ such that $r_m\to 1$
and WOT-$\lim_{m\to 1} \nu_{h, r_m}(f_j)$ exists for each $f_j$. Fix
$f\in \cR_n^*+ \cR_n$ and $x,y\in \cE$ and let us
 prove that $\{\left<\nu_{h,r_m}(f)x,y\right>\}_{m=1}^\infty$ is a Cauchy sequence. Let
 $\epsilon>0$ and choose  $f_j$ so that
 $\|f_j-f\|<\frac{\epsilon}{3M\|x\|\|y\|}$. Now, we choose $N$ so that
 $
 |\left<(\nu_{h, r_m}(f_j)-\nu_{h, r_k}(f_j))x, y\right>|
 <\frac{\epsilon}{3}\quad \text{ for
 any } \ m,k>N.
 $
Due to the fact that
\begin{equation*}
\begin{split}
|\left<(\nu_{h,r_m}(f)-\nu_{h,r_k}(f))x,y\right>| &\leq
|\left<(\nu_{h,r_m}(f-f_j)x,y\right>|+
|\left<(\nu_{h,r_m}(f_j)-\nu_{h,r_k}(f_j))x,y\right>| \\
&\qquad  + |\left<\nu_{h,r_k}(f_j-f)x,y\right>|\\
&\leq 2 M\|x\|\|y\| \|f-f_j\|+ |\left<(\nu_{h,
r_m}(f_j)-\nu_{h,r_k}(f_j))x,y\right>|
\end{split}
\end{equation*}
we deduce that $|\left<(\nu_{h,r_m}(f)-\nu_{h,r_k}(f))x,
y\right>|<\epsilon$ \ for $m,k>N$. Therefore, we deduce that
$b(x,y):=\lim_{m\to\infty} \left<\nu_{h,r_m}(f)x,y\right>$ exists
for any $x,y\in \cE$ and defines a  functional  $b:\cE\times \cE\to
\CC$ which is linear in the first variable and conjugate linear in
the second.
 Moreover, we have $|b(x,y)|\leq M \|f\|\|x\|\|y\|$ for any $x,y\in \cE$.
 Due to Riesz representation theorem, there exists a unique bounded linear
 operator $B(\cE)$, which we denote by $\nu(f)$, such that
 $b(x,y)=\left<\nu(f)x,y\right>$ for $x,y\in \cE$. Therefore,
 $
 \nu(f)=\text{\rm WOT-}\lim_{r_m\to 1} \nu_{h, r_m}(f)
 $
 for any $f\in \cR_n^*+\cR_n$, and $\|\nu(f)\|\leq M\|f\|$.
Notice that $\nu:\cR_n^*+\cR_n\to B(\cE)$ is a bounded linear map
with $\|\nu\|\leq M$.
 Moreover, $\nu$ is a completely bounded map. Indeed, if $[f_{ij}]_{m}$ is
 an $m\times m$ matrix over $\cR_n^*+\cR_n$, then
 $[\nu(f_{ij})]_{m}=\text{\rm WOT-}\lim_{r_k\to 1}
 [\nu_{h,r_k}(f_{ij})]_{m}.
 $
Hence, $\left\|[\nu(f_{ij})]_{m}\right\|\leq M
\left\|[f_{ij}]_{m}\right\|$ for all $m$, and so $\|\nu\|_{cb}\leq
M$. Notice that, in particular, we have $
\nu(R_{\widetilde\alpha}^*)= A_{(\alpha)}$, \ $\alpha\in \FF_n^+$,
and $\nu(R_{\widetilde\alpha})=  B_{(\alpha)}$, \ $\alpha\in
\FF_n^+\backslash \{g_0\}$, where $\{A_{(\alpha)}\}$ and
$\{B_{(\alpha)}\}$ are the coefficients of $h$. According to
Wittstok's extension
 theorem \cite{W2}, there exists a completely bounded linear map
$\mu:C^*(R_1,\ldots, R_n)\to \CC$ such that $
\mu(R_{\widetilde\alpha}^*)= A_{(\alpha)}$, \ $\alpha\in \FF_n^+$,
and $\mu(R_{\widetilde\alpha})=  B_{(\alpha)}$, \ $\alpha\in
\FF_n^+\backslash \{g_0\}$ and such that $\|\mu\|=\|\nu\|$.
Consequently,  $h=\cP\mu$ and   item (ii) holds.

To prove the implication (i)$\implies$(iii), we  apply
  Wittstock's decomposition theorem \cite{W1}
 to the completely
 bounded linear map $\mu:C^*(R_1,\ldots,
R_n)\to B(\cE)$. Thus $\mu$  has a decomposition of the form
$$
\mu=(\mu_1-\mu_2)+i(\mu_3-\mu_4)$$ where $\mu_1,\mu_2, \mu_3, \mu_4$
are completely positive linear maps on  $C^*(R_1,\ldots, R_n)$ with
values in $B(\cE)$. Due to the linearity of the noncommutative
Poisson transform, we have $
h=(\cP\mu_1-\cP\mu_2)+i(\cP\mu_3-\cP\mu_4)$. Since $(\cP
\mu_j)(X)=\widetilde \mu_j[P(R,X)]$, $j=1,\ldots, 4$ and, due to
Proposition \ref{Poisson-factor}, $P(R,X)\geq 0$ for $X\in
[B(\cH)^n]_1$, we deduce that $\cP \mu_j$, $j=1,\ldots, 4$, are
positive free pluriharmonic functions. Hence  we deduce (iii).

It remains to show that (iii)$\implies$(i). To this end, we assume
that (iii) holds.  Applying  Corollary \ref{positive-pluri} to the
positive free holomorphic functions $h_1, h_2, h_3, h_4$, we find
$\mu_1,\mu_2, \mu_3, \mu_4$, some completely positive linear maps on
$C^*(R_1,\ldots, R_n)$ with values in $B(\cE)$, such that $h_s=\cP
\mu_s$ for $s=1,\ldots,4$. Setting
$\mu:=(\mu_1-\mu_2)+i(\mu_3-\mu_4)$ and using item (iii), we deduce
that $h=\cP\mu$. This completes the proof.
\end{proof}

We remark that, due to Theorem \ref{pluri-measure}, the map
$\mu\mapsto \cP\mu$ is a linear one-to-one correspondence between
the set of all completely bounded linear maps on $\cR_n^*+ \cR_n$
and the set
$$
\{(u_1-u_2)+i(u_3-u_4):\ u_j\geq 0,\  u_j\in Har(B(\cH)^n_1)\}.
$$

Using  again Theorem  \ref{pluri-measure} and the Jordan type
decomposition for selfadjoint  completely bounded linear maps  on
$C^*$-algebras (see \cite{P}), one can easily   deduce the following
result.

\begin{corollary}\label{Jordan}
Let \ $u$ be a selfadjoint free pluriharmonic
function on $[B(\cH)^n]_1$ with coefficients in $B(\cE)$. Then the following statements are equivalent:
\begin{enumerate}
\item[(i)]
$u$ admits a Jordan decomposition $u=u_+-u_-$, where $u_+$ and $u_-$
are positive free pluriharmonic functions on $[B(\cH)^n]_1$;

\item[(ii)]
the  selfadjoint    completely  bounded linear maps \ $\{\nu_{h,r}\}_{r\in[0,1)}$,
associate with the radial function of $h$,
are uniformly  bounded, i.e.,   $\sup\limits_{0\leq r<1}
\|\nu_{h,r}\|_{cb}<\infty$;
\item[(iii)]
there exists a  selfadjoint completely bounded linear map
$\mu:C^*(R_1,\ldots, R_n)\to B(\cE)$ such that $ h=\cP\mu. $
\end{enumerate}

Moreover, one can choose     $u_+=\cP\mu_+$ and $u_-=\cP\mu_-$,
where $\mu=\mu_+-\mu_-$ is the Jordan decomposition of $\mu$, i.e.,
$\mu_+, \mu_-\geq 0$ and $\|\mu\|=\|\mu_+\|+\|\mu_-\|$.
\end{corollary}

 The map $\mu\mapsto \cP\mu$ is a linear one-to-one
correspondence between the set of all  selfadjoint completely
bounded linear maps on $\cR_n^*+ \cR_n$ and the set $ \{u_1-u_2 :\
u_j\geq 0,\  u_j\in Har(B(\cH)^n_1)\}. $

We introduce  now the space $Har^1(B(\cH)^n_1)$  of all    free
pluriharmonic functions  $h$ on $[B(\cH)^n]_1$ with coefficients in $B(\cE)$
 such that
$\sup_{0\leq r<1}\|\nu_{h,r}\|<\infty$ and define
$\|h\|_1:=\sup_{0\leq r<1}\|\nu_{h,r}\|. $
 It is easy to see that
$\|\cdot\|_1$ is a norm on $Har^1(B(\cH)^n_1)$. Denote by CB\,$(
{\cR_n^*+ \cR_n}, B(\cE))$
 the space of all completely bounded linear maps
 from $  {\cR_n^*+ \cR_n}$ to $B(\cE)$.

\begin{theorem}\label{Har-1}
$\left(Har^1(B(\cH)^n_1), \|\cdot \|_1\right)$ is a Banach space
which can be identified with the Banach space
 \text{\rm CB}\,$( {\cR_n^*+ \cR_n}, B(\cE))$.
Moreover, the following statements are equivalent:
\begin{enumerate}
\item[(i)]
$h$ is in $Har^1(B(\cH)^n_1)$;
\item[(ii)]
there is a unique completely bounded linear  map $\mu_h: {\cR_n^*+
\cR_n}\to B(\cE)$ such that \ $h=\cP\mu_h$;
\item[(iii)] there exists an
$n$-tuple of isometries $(V_1,\ldots, V_n)$
on a Hilbert space $\cK$, with orthogonal ranges,  and    bounded operators $W_i:\cE\to \cK$,
$i=1,2$, such that
$$
h(X_1,\ldots, X_n)=  (W_1^*\otimes I)\left[ B_X(V_1,\ldots, V_n)^*
B_X(V_1,\ldots, V_n) \right](W_2\otimes I),
$$
where $B_X(V_1,\ldots, V_n):=(I\otimes \Delta_X)(I-V_1\otimes
X_1^*-\cdots -V_n\otimes X_n^*)^{-1}. $
\end{enumerate}
\end{theorem}

\begin{proof}
Define the map $\Psi:\text{\rm CB}\,( {\cR_n^*+ \cR_n}, B(\cE))\to
Har^1(B(\cH)^n_1)$   by $\Psi(\mu):=\cP\mu$. To prove injectivity of
$\Psi$, let $\mu_1,\mu_2$ be in $\text{\rm CB}\,( {\cR_n^*+ \cR_n},
B(\cE))$ such that $\Psi(\mu_1)=\Psi(\mu_2)$. Then, due to the
uniqueness  of the representation of a free pluriharmonic function
and the definition of the noncommutative Poisson transform of a
completely bounded map on $\cR_n^*+\cR_n$, we deduce that
$\mu_1(R_\alpha)=\mu_2(R_\alpha)$, $\alpha\in \FF_n^+$, and
$\mu_1(R_\alpha^*)=\mu_2(R_\alpha^*)$, $\alpha\in \FF_n^+\backslash
\{g_0\}$. Hence, we have  $\mu_1=\mu_2$. The surjectivity of the map
$\Psi$ is due to Theorem \ref{pluri-measure}. The same theorem  (see
the proof) implies $\|\cP\mu\|_1=\|\mu\|$ for any $\mu$ in \text{\rm
CB}\,$( {\cR_n^*+ \cR_n}, B(\cE))$. This completes the proof of the
equivalence of (i) with (ii) and the identification of
$Har^1(B(\cH)^n_1)$ with \text{\rm CB}\,$( {\cR_n^*+ \cR_n},
B(\cE))$.

To prove (iii), notice that part (i) and Proposition
\ref{Poisson-factor} imply
\begin{equation}
\label{HBP} h(X_1,\ldots, X_n)=(\cP\mu_h)(X_1,\ldots,
X_n)=\widetilde \mu_h(B_X^* B_X),
\end{equation}
where $B_X:=(I\otimes \Delta_X)(I-R_1\otimes X_1^*-\cdots -R_n\otimes
X_n^*)$. On the other  hand,  by Wittstock's extension theorem
\cite{W2}, there exists a completely bounded map
$\phi:C^*(R_1,\ldots, R_n)\to B(\cE)$ that extends $\mu_h$ with
$\|\mu_h\|_{cb}=\|\phi\|_{cb}$. According to Theorem 8.4 from
\cite{P}, which is a generalization of Stinespring's representation
theorem \cite{St}, there exists a Hilbert space $\cK$, a
$*$-representation $\pi:C^*(R_1,\ldots, R_n)\to B(\cK)$, and bounded
operators $W_j:\cE\to \cK$, $j=1,2$, with $\|\phi\|=\|W_1\|\|W_2\|$
such that
$$
\phi(f)=W_1^*\pi(f) W_2,\qquad f\in C^*(R_1,\ldots, R_n).
$$
Notice that $V_i:=\pi(R_i)$, $i=1,\ldots, n$, are isometries with
orthogonal ranges and any $*$-representation of $C^*(R_1,\ldots,
R_n)$ is generated by $n$ isometries with orthogonal ranges. Using
now  relation \eqref{HBP}, one can  complete the proof of part (iii).
\end{proof}

Consider now the space of free holomorphic functions
$H^1(B(\cH)^n_1):=Hol(B(\cH)^n_1)\bigcap Har^1(B(\cH)^n_1)$ together
with the norm $\|\cdot\|_1$. The following result is a consequence
of Theorem \ref{Har-1} and  a weak version of the F. and M.~Riesz
theorem \cite{H}, in our setting.

\begin{corollary}\label{FMR}
$\left(H^1(B(\cH)^n_1), \|\cdot \|_1\right)$ is a Banach space which
can be identified with the annihilator of $\cR_n$ in \text{\rm
CB}\,$( {\cR_n^*+ \cR_n}, B(\cE))$, i.e.,
$$
(\cR_n)^{\perp}:=\{\mu\in \text{\rm CB}\,( {\cR_n^*+ \cR_n},
B(\cE)): \ \mu(R_\alpha)=0 \text{ for all } |\alpha|\geq 1\}.
$$
In particular, for each $f\in H^1(B(\cH)^n_1)$, there is a unique
completely bounded linear map $\mu_f\in (\cR_n)^{\perp}$ such that \
$f=\cP\mu_f$.
\end{corollary}

Let $H^2(B(\cH)^n_1)$ be the set of all free holomorphic functions
on $[B(\cH)^n]_1$ with operator-valued coefficients in $B(\cE)$,
 of the form
$\varphi(X_1,\ldots, X_n)=\sum_{k=0}^\infty \sum_{|\alpha|=k}
 A_{(\alpha)}\otimes X_\alpha$
such that $ \|\varphi\|_2:=\left\|\sum_{\alpha \in \FF_n^+}
A_{(\alpha)}^* A_{(\alpha)}\right\|^{1/2} $ is finite. It is clear
that $(H^2(B(\cH)^n_1), \|\cdot\|_2)$ is a Banach space. We  recall
\cite{Po-holomorphic} that $H^\infty(B(\cH)^n_1)$ is the set of  all
bounded free holomorphic functions on $[B(\cH)^n]_1$. Due to the
results of Section 3, it is clear that
$$
H^\infty(B(\cH)^n_1)=Hol(B(\cH)^n_1)\bigcap Har^\infty(B(\cH)^n_1).
$$

\begin{proposition} $H^\infty(B(\cH)^n_1)\subset H^2(B(\cH)^n_1)\subset
H^1(B(\cH)^n_1)$ and the inclusions are continuous.
\end{proposition}

\begin{proof}
Let  $\varphi\in H^\infty(B(\cH)^n_1)$ have   the representation
$\varphi(X_1,\ldots, X_n)=\sum_{k=0}^\infty \sum_{|\alpha|=k}
A_{(\alpha)}\otimes X_\alpha$. Then its boundary function has the
Fourier representation $\sum_{k=0}^\infty \sum_{|\alpha|=k}
 A_{(\alpha)}\otimes S_\alpha$. Note that, for any  $x\in \cE$ with $\|x\|=1$, we have
 \begin{equation*}
 \begin{split}
 \|\varphi\|_2\leq \left(\sum_{\alpha\in \FF_n^+}
  \|A_{(\alpha)}x\|^2\right)^{1/2}
  =\left\|\left(\sum_{k=0}^\infty \sum_{|\alpha|=k}
 A_{(\alpha)}\otimes S_\alpha\right)(x\otimes 1)\right\|\leq \|\varphi\|_\infty.
 \end{split}
 \end{equation*}
Therefore,  $\varphi\in H^2(B(\cH)^n_1)$ and $\|\varphi\|_2\leq
\|\varphi\|_\infty$.
   Assume  now that $\varphi\in
H^2(B(\cH)^n_1)$. Define the linear map $\mu_\varphi:\cR_n^*
+\cR_n\to B(\cE)$ by $\mu_\varphi(R_\alpha)=0$ for  $\alpha\in
\FF_n^+\backslash \{g_0\}$, and
$\mu_\varphi(R_{\widetilde\alpha}^*)=A_{(\alpha)}$ for $\alpha\in
\FF_n^+$. Due to Corollary \ref{FMR}, to show that $\varphi\in
H^1(B(\cH)^n_1)$,  it is  enough to prove that $\mu_\varphi\in
(\cR_n)^{\perp}$. For any $m\in\NN$, we have
\begin{equation*}
\begin{split}
\left\|\mu_\varphi\left(\sum_{|\alpha|\leq m} c_\alpha R_\alpha^* +
\sum_{|\alpha|\leq m} d_\alpha R_\alpha\right)\right\|
&=\left\|\sum_{|\alpha|\leq m} c_\alpha
A_{(\widetilde\alpha)}\right\| \leq \left(\sum_{|\alpha|\leq m}
|c_\alpha|^2\right)^{1/2}
\left\|\sum_{|\alpha|\leq m} A_{(\alpha)}^* A_{(\alpha)}\right\|^{1/2}\\
&= \left\|\left(\sum_{|\alpha|\leq m} \bar c_\alpha R_\alpha +
\sum_{|\alpha|\leq m} \bar d_\alpha
R_\alpha^*\right)(1)\right\|\left(\sum_{\alpha\in \FF_n^+}
A_{(\alpha)}^* A_{(\alpha)}\right)^{1/2}\\
&\leq \left\|\sum_{|\alpha|\leq m} c_\alpha R_\alpha^* +
\sum_{|\alpha|\leq m} d_\alpha R_\alpha\right\| \|\varphi\|_2.
\end{split}
\end{equation*}
Hence and using Corollary \ref{FMR}, we deduce that
$\|\varphi\|_1=\|\mu_\varphi\|\leq \|\varphi\|_2$. Similarly, passing to matrices
over $\cR_n^*+\cR_n$, one can show that $\mu_\varphi$ is a
completely bounded map.
 Since $\varphi$ is a free
pluriharmonic function and $\varphi=\cP \mu_\varphi$,  it is clear
that $\mu_\varphi$ is the only completely bounded map  on $\cR_n^*+
\cR_n$ with this property.
\end{proof}

\begin{remark}
If $\varphi\in H_\CC^2(B(\cH)^n_1)$ and $\mu_\varphi$ is the
associated   bounded linear functional on $\cR_n^*+\cR_n$, then
    there are
vectors $\eta, \xi\in F^2(H_n)$ such that
$\mu_\varphi(R_{ \alpha}^*)=
\left<R_{ \alpha}^* \eta, \xi\right>$ for any $\alpha\in \FF_n^+$ and
$$
\varphi(X_1,\ldots, X_n)=\sum_{k=0}^\infty \sum_{|\alpha|=k}
\left<R^*_{\widetilde\alpha} \eta, \xi\right> X_\alpha, \qquad (X_1,\ldots,
X_n)\in [B(\cH)^n]_1.
$$
\end{remark}
\begin{proof}

Denote by $\mu_\varphi^*$ the linear functional on $\cR_n^*+ \cR_n$ defined by
$\mu_\varphi^*(f):=\overline{\mu_\varphi(f^*)}$.
 Let $\cR_{n,0}^k$  be the norm closed linear span of the operators
$R_\alpha$, $|\alpha|\geq k$ and assume that $\varphi$ has the representation
$\varphi(X_1,\ldots, X_n):=\sum_{k=0}^\infty \sum_{|\alpha|=k} a_\alpha X_\alpha$. Notice that, for any $m\geq k$, we
have
\begin{equation*}
\begin{split}
\left|\mu_\varphi^*\left( \sum_{k\leq|\alpha|\leq m} d_\alpha
R_\alpha\right)\right|&=\left| \sum_{k\leq|\alpha|\leq m}  d_\alpha
\bar a_{\widetilde\alpha}\right| \leq \left(\sum_{k\leq|\alpha|\leq
m} |d_\alpha|^2\right)^{1/2}
\left(\sum_{k\leq|\alpha|\leq m} |a_\alpha|^2\right)^{1/2}\\
&\leq \left\|\left( \sum_{k\leq|\alpha|\leq m} d_\alpha
R_\alpha\right)(1)\right\|\left(\sum_{k\leq |\alpha|}
|a_\alpha|^2\right)^{1/2} \leq \left\|\sum_{k\leq|\alpha|\leq m}
d_\alpha R_\alpha\right\| \left(\sum_{k\leq |\alpha|}
|a_\alpha|^2\right)^{1/2}.
\end{split}
\end{equation*}
Therefore, $$\left\|\mu_\varphi^*|_{\cR_{n,0}^k}\right\|\leq
\left(\sum_{k\leq |\alpha|} |a_\alpha|^2\right)^{1/2}\to 0,\quad
\text{ as } \ k\to\infty.
$$
Using Proposition 2.2 from \cite{DLP}, we deduce that
$\mu_\varphi^*|_{\cR_n}$ is an absolutely continuous functional on
$\cR_n$, i.e.,  there are vectors $\xi, \eta\in F^2(H_n)$ such that
$\mu_\varphi^*(A)=\left< A\xi, \eta\right>$ for any $A\in \cR_n$.
Therefore, we have $\mu_\varphi(R_{
\alpha}^*)=\overline{\mu_\varphi^*(R_{\alpha})}= \left<R_{ \alpha}^*
\eta, \xi\right>$. Since $\varphi=\cP \mu_\varphi$, we  complete the
proof.
\end{proof}

The remark above  leads to the following question: can  the set
 $(\cR_n)^{\perp}$  be identified with the set of all  absolutely continuous
  functionals
 on  $\cR_n^*$ ?
If the answer is positive, then it will constitute  a noncommutative
multivariable generalization of F. and M.~Riesz theorem \cite{H}.

Using  Theorem  \ref{pluri-measure}, we can recast     some results
from \cite{PPoS} and \cite{Po-Bohr},
 and  obtain the
following  Wiener and Bohr type inequalities (see \cite{B}) for the
analytic moments of a selfadjoint linear functional on
$C^*(R_1,\ldots, R_n)$. Let $\mu, \tau :C^*(R_1,\ldots, R_n)\to \CC$
be  selfadjoint linear functionals, where $\tau$ is defined by
 $\tau(g):=\left<g(1), 1\right>$.  If
$\mu\leq \tau$ on $\cR_n^*+ \cR_n$, then
\begin{enumerate}
\item[(i)]
$\left(\sum\limits_{|\alpha|=k}|\mu(R_\alpha)|^2\right)^{1/2}\leq
1-|\mu(I)|$;
\item[(ii)]
$\sum\limits_{k=0}^\infty\sum\limits_{|\alpha|=k}|\mu(R_\alpha)|
r_\alpha\leq 1$ \ if $\|(r_1,\ldots, r_n)\|\leq \frac{1}{2}$;
\item[(iii)]
$|\mu(I)|+2\sum\limits_{k=0}^\infty\sum_{|\alpha|=k}|\mu(R_\alpha)|
r_\alpha\leq 1$ \  if $\|(r_1,\ldots, r_n)\|\leq \frac{1}{3}$.
\end{enumerate}

\bigskip

\section{Noncommutative Cayley transforms}
We introduce Cayley type transforms acting on formal power series, contractive
 free holomorphic functions on the noncommutative ball $[B(\cH)^n]_1$,
 and multi-analytic  matrices, respectively. These transforms are needed in the next
  section to solve the Carath\' eodory  interpolation problem for free holomorphic
  functions with positive real parts.

Let $f=\sum_{|\alpha|\geq 1} A_{(\alpha)}\otimes Z_\alpha$ be a
formal power series in noncommutative indeterminates $Z_1,\ldots,
Z_n$, coefficients in $B(\cE)$,   and  constant term $0$. For each
$m\in \NN$, $f^m$  defines a power series $\sum_{|\alpha|\geq m}
C_{(\alpha)}\otimes  Z_\alpha$  and it makes sense to consider the
formal power series $ \varphi=1+f+f^2+\cdots. $ We call
$C_{(\alpha)}$ the $\alpha$-coefficient of $f^m$. Notice that if
$m>k$, then the term $f^m$ has all coefficients of order $\leq k$
equal to $0$. Thus, if $\alpha\in \FF_n^+$ with $|\alpha|\leq k$,
then we may define the $\alpha$-coefficient of $\varphi$ as the
$\alpha$-coefficient of the finite sum $1+f+ f^2+\cdots + f^m$.
Notice that $ \varphi=1+\sum_ {|\alpha|\geq 1} B_{(\alpha)}\otimes
Z_\alpha, $ where
\begin{equation}\label{b-alpha}
 B_{(\alpha)}= \sum_{j=1}^{|\alpha|} \sum_{{\gamma_1\cdots
\gamma_j=\alpha }\atop {|\gamma_1|\geq 1,\ldots, |\gamma_j|\geq 1}}
A_{(\gamma_1)}\cdots A_{(\gamma_j)}   \quad \text{ for } \
|\alpha|\geq 1.
\end{equation}
Since $(1-f)\varphi=\varphi (1-f)=1$, we have $(1-f)^{-1}=\varphi$.
Denote  by $\widetilde \CC_0[Z_1,\ldots, Z_n]$ the algebra of all
formal power series in noncommutative indeterminates $Z_1,\ldots,
Z_n$, coefficients in $B(\cE)$,  and constant term $0$. We introduce
the Cayley transform $\widetilde \cC: \widetilde \CC_0[Z_1,\ldots,
Z_n]\to \widetilde \CC_0[Z_1,\ldots, Z_n]$ by setting
$$
\widetilde \cC(f):=(1-f)^{-1}f,\quad f\in \widetilde
\CC_0[Z_1,\ldots, Z_n].
$$

\begin{proposition}\label{Cayley-series}
The Cayley transform for formal power series is a bijection  and
$$\widetilde \cC^{-1}(f)=f(1+f)^{-1},\quad f\in \widetilde
\CC_0[Z_1,\ldots, Z_n].
$$
\end{proposition}
\begin{proof}
If $f_1,f_2\in \widetilde \CC_0[Z_1,\ldots, Z_n]$ and  $\widetilde
\cC(f_1)=\widetilde \cC (f_2)$, then $f_1(1-f_2)=(1-f_1)f_2$, whence
$f_1=f_2$. To prove  that the Cayley transform is surjective, let
$f\in \widetilde \CC_0[Z_1,\ldots, Z_n]$ and notice that
\begin{equation*}
\begin{split}
 \widetilde \cC [f(1+f)^{-1}]&=
 \left[1-f(1+f)^{-1}\right] f(1+f)^{-1} =\left[(1+f)^{-1} (1+f-f)\right]^{-1} (1+f)^{-1}f =f.
\end{split}
\end{equation*}
This completes the proof.
\end{proof}

Denote by $\cP^{(m)}$ the set of all polynomials of degree $\leq 1$ in
$F^2(H_n)$, i.e.,
$$
\cP^{(m)}:=\text{ \rm span} \{ e_\alpha: \ \alpha\in \FF_n^+,
|\alpha|\leq m\},
$$
and define the nilpotent operators $S_i^{(m)}: \cP^{(m)}\to
\cP^{(m)}$ by $ S_i^{(m)}:=P_{\cP^{(m)}} S_i |_{\cP^{(m)}},\quad
i=1,\ldots, n, $ where $S_1,\ldots, S_n$ are the left creation
operators on the Fock space $F^2(H_n)$ and $P_{\cP^{(m)}}$ is the
orthogonal projection of $F^2(H_n)$ onto $\cP^{(m)}$. Notice that
$S_\alpha^{(m)}=0$ if $|\alpha|\geq m+1$. According to
\cite{Po-unitary}, the $n$-tuple of operators $(S_1^{(m)},\ldots,
S_n^{(m)})$ is the universal model for row contractions
$(T_1,\ldots, T_n)$  with $T_\alpha=0$ for $|\alpha|\geq m+1$, and
the following constrained von Neumann inequality holds:
\begin{equation}
\label{vN-nilp}
 \|p(T_1,\ldots, T_n)\|\leq \|p(S_1^{(m)},\ldots,
S_n^{(m)})\|
\end{equation}
for any noncommutative polynomial $p(X_1,\ldots,
X_n)=\sum_{|\alpha|\leq k} A_{(\alpha)}\otimes X_\alpha$, \ $k\in
\NN$.

\begin{lemma}
\label{nilp} Let $f, g$ be free holomorphic functions on the
noncommutative open  ball $[B(\cH)^n]_1$ with operator-valued
coefficients.
\begin{enumerate}
\item[(i)]
If $f$ has the representation $f(X_1,\ldots, X_n):=\sum_{k=0}^\infty
\sum_{|\alpha|=k} A_{(\alpha)}\otimes  X_\alpha$, then
$$f(S_1^{(m)},\ldots,
S_n^{(m)})=\sum_{k=0}^m \sum_{|\alpha|=k} A_{(\alpha)}\otimes
S_\alpha^{(m)}$$ is a bounded linear operator on $\cE\otimes \cP^{(m)}$
and
$$
\|f(S_1^{(m)},\ldots, S_n^{(m)})\|=\sup_{r\in
[0,1)}\|f(rS_1^{(m)},\ldots, rS_n^{(m)})\|.
$$
\item[(ii)] $f=g$ if and only if \ $f(S_1^{(m)},\ldots,
S_n^{(m)})=g(S_1^{(m)},\ldots, S_n^{(m)})$ for any $m\in \NN$.
\item[(iii)] $f\in H^\infty(B(\cH)^n_1)$ if and only if
\ $\sup_{m\in \NN} \|f(S_1^{(m)},\ldots, S_n^{(m)})\|<\infty$.
Moreover, in this case,
$$\|f\|_\infty=\sup_{m\in \NN} \|f(S_1^{(m)},\ldots, S_n^{(m)})\|.
$$
\end{enumerate}
\end{lemma}

\begin{proof}
Since $S^{(m)}_\alpha=0$ for $\alpha\in \FF_n^+$ with $|\alpha|\geq
m+1$, and using inequality \eqref{vN-nilp}, we have
\begin{equation*}
\begin{split}
 \|f(rS_1^{(m)},\ldots, rS_n^{(m)})\|&=\left\|\sum_{k=0}^m
\sum_{|\alpha|=k} A_{(\alpha)}\otimes r^{|\alpha|}
S_\alpha^{(m)}\right\|\\
&\leq \left\|\sum_{k=0}^m \sum_{|\alpha|=k} A_{(\alpha)}\otimes
  S_\alpha^{(m)}\right\|=\|f(S_1^{(m)},\ldots, S_n^{(m)})\|
  \end{split}
  \end{equation*}
  for any $r\in [0,1)$.
Since $\lim_{r\to 1} f(rS_1^{(m)},\ldots, rS_n^{(m)})=
f(S_1^{(m)},\ldots, S_n^{(m)})$ in the operator norm topology, we
deduce part (i). Part (ii) is obvious, so we prove (iii). According
to \cite{Po-holomorphic}, if $f\in H^\infty(B(\cH)^n_1)$, then
$\|f\|_\infty=\sup_{r\in [0,1)}\|f(rS_1,\ldots, rS_n)\|$.
 Since $f(rS_1^{(m)},\ldots, rS_n^{(m)})=P_{\cE\otimes \cP^{(m)}}f(rS_1,\ldots,
 rS_n)|_{\cE\otimes \cP^{(m)}}$, $r\in [0,1)$, we deduce that
$\|f(S_1^{(m)},\ldots, S_n^{(m)})\|\leq \|f\|_\infty$ for any $m\in
\NN$.

Conversely, assume that $\sup_{m\in \NN} \|f(S_1^{(m)},\ldots,
S_n^{(m)})\|<\infty$ and $f\notin H^\infty(B(\cH)^n_1)$. Then for
any $M>0$ there exists $r_0\in [0,1)$ such that $\|f(r_0S_1,\ldots,
r_0S_n)^*\|>M$. Consequently, we can find a vector
$q=\sum_{|\alpha|\leq k} h_\alpha\otimes e_\alpha$  of norm one such
that $\|f(r_0S_1,\ldots, r_0S_n)^*q\|>M$. Notice that
\begin{equation*}
\begin{split}
\|f(S_1^{(k)},\ldots, S_n^{(k)})\|&\geq  \|f(r_0S_1^{(k)},\ldots,
r_0S_n^{(k)})^*\|\geq \|f(r_0S_1,\ldots, r_0S_n)^*q\|>M,
\end{split}
\end{equation*}
which implies $\sup_{m\in \NN} \|f(S_1^{(m)},\ldots,
S_n^{(m)})\|>M$. Hence, we get a contradiction. Therefore, we must
have  $f\in H^\infty(B(\cH)^n_1)$. Moreover, the considerations
above can be used to deduce that $\|f\|_\infty=\sup_{m\in \NN}
\|f(S_1^{(m)},\ldots, S_n^{(m)})\|$. The proof is complete.
\end{proof}

We remark that a result similar to that of Lemma \ref{nilp} holds
for free pluriharmonic functions. The proof is basically the same
but uses the results of Section 3. Consider now  the sets
\begin{equation*}
\begin{split}
H_0^\infty(B(\cH)^n_1)&:=\left\{ f\in  H^\infty(B(\cH)^n_1): \
f(0)=0\right\}\\
Hol_0^+(B(\cH)^n_1)&:=\left\{ g\in Hol(B(\cH)^n_1):\ g(0)=0, \
g(X)^* +I +g(X)\geq 0 \ \text{  for  } \ X\in [B(\cH)^n]_1 \right\}.
\end{split}
\end{equation*}
We introduce the noncommutative Cayley  transform
$$\cC:\left[H_0^\infty(B(\cH)^n_1)\right]_{\leq 1} \to
Hol_0^+(B(\cH)^n_1)\quad \text{ defined by } \ \cC f:=g,$$
 where $g$ is the
free holomorphic function on $[B(\cH)^n]_1$ uniquely determined by
the formal power series
 $(1-\widetilde f)^{-1}\widetilde f$, where
$\widetilde f$ is the power series associated with $f$. Of course,
it remains to show that  $\cC$ is well-defined.

\begin{theorem}
\label{Cayley1} The noncommutative Cayley transform  is a bijection
between the unit ball \ $\left[H_0^\infty(B(\cH)^n_1)\right]_{\leq
1}$ and $Hol_0^+(B(\cH)^n_1)$.
\end{theorem}
\begin{proof}
First, we show that the map $\cC$ is well-defined.  Let $f$ be in
$\left[H_0^\infty(B(\cH)^n_1)\right]_{\leq 1}$ and have the
representation $f(X_1,\ldots,
X_n):=\sum_{k=1}^\infty\sum_{|\alpha|=k} A_{(\alpha)}\otimes
X_\alpha$, \ $(X_1,\ldots, X_n)\in [B(\cH)^n]_1. $ We shall prove
that $\cC(f)\in Hol_0^+(B(\cH)^n_1$.
 Due to the
Schwartz type lemma for bounded free holomorphic functions on the
open unit ball of $B(\cH)^n$ (see \cite{Po-holomorphic}), we have $
\|f(X_1,\ldots, X_n)\|\leq \|[X_1,\ldots, X_n]\|$  for $(X_1,\ldots,
X_n)\in [ B(\cH)^n]_1. $
  Since  $f(rS_1,\ldots, rS_n)$ is in $B(\cE)\otimes_{min} \cA_n$
and  $\|f(rS_1,\ldots, rS_n)\|\leq r<1$,  the operator
   $I-f(rS_1,\ldots, rS_n)$ is invertible with its
inverse
 $(I-f(rS_1,\ldots, rS_n))^{-1}$
 in $B(\cE)\otimes_{min} \cA_n\subset B(\cE)\bar\otimes F_n^\infty$, where $F_n^\infty$ is the
 noncommutative analytic Toeplitz algebra \cite{Po-von}, i.e., the weakly closed
 algebra generated by the  left creation operators and the identity.
 Therefore, the operator
$ (I-f(rS_1,\ldots, rS_n))^{-1}f(rS_1,\ldots, rS_n)$ is in
$B(\cE)\otimes\cA_n$ and has a representation
 $\sum_{\alpha\in \FF_n^+} B_{(\alpha)}\otimes
r^{|\alpha|} S_\alpha $ for some operators  $B_{(\alpha)}\in
B(\cE)$. Using the fact that
$$r^{|\alpha|}B_{(\alpha)}=P_{\cE\otimes \CC} (I_\cE\otimes S_\alpha^*)
(I-f(rS_1,\ldots, rS_n))^{-1}f(rS_1,\ldots,
rS_n)|_{\cE\otimes \CC}, $$
 we  deduce that
\begin{equation*}
\begin{split}
(I-f(rS_1,\ldots, rS_n))^{-1}f(rS_1,\ldots,
rS_n)&=  f(rS_1,\ldots, rS_n)+ f(rS_1,\ldots, rS_n)^2+\cdots \\
  &=\sum_{k=1}^\infty\sum_{|\alpha|=k}
B_{(\alpha)}\otimes r^{|\alpha|} S_\alpha,
\end{split}
\end{equation*}
where the coefficients $B_{(\alpha)}$ are given by relation
\eqref{b-alpha}. Hence and  due to the   definition of the
noncommutative Cayley transform, we have
$$
\cC(f)(rS_1,\ldots, rS_n):=g(rS_1,\ldots, rS_n)=(I-f(rS_1,\ldots, rS_n))^{-1}f(rS_1,\ldots,
rS_n), \quad r\in[0,1).
$$
This shows that $g(X_1,\ldots, X_n)
=\sum_{k=1}^\infty\sum_{|\alpha|=k} B_{(\alpha)}\otimes   X_\alpha$
is a free holomorphic function on $[B(\cH)^n]_1$, and $g(0)=0$. Now,
we prove that
\begin{equation}
\label{g*g}
 g(X)^* +I +g(X)\geq 0 \quad \text{   for any  }
X\in [B(\cH)^n]_1.
\end{equation}
 Due to the Poisson mean value property of Theorem \ref{MVP}, it is enough to show that
$$
h(rS_1,\ldots, rS_n):=g(rS_1,\ldots, rS_n)^* + I +g(rS_1,\ldots,
rS_n)\geq  0\quad \text{ for any  } \ r\in [0,1).
$$
Notice that
  \begin{equation*}
  \begin{split}
  f(rS_1,\ldots, rS_n)[I+g(rS_1,\ldots, rS_n)]
  &=f(rS_1,\ldots, rS_n)\left\{I+[I-f(rS_1,\ldots, rS_n)]^{-1}f(rS_1,\ldots,
  rS_n)\right\}\\
  &=f(rS_1,\ldots, rS_n)[I-f(rS_1,\ldots, rS_n)]^{-1}
  =g(rS_1,\ldots, rS_n).
\end{split}
  \end{equation*}
  Using this relation, we deduce that
  \begin{equation*}
  \begin{split}
h(rS_1,\ldots, rS_n)&=[I+g(rS_1,\ldots, rS_n)]^* [I+g(rS_1,\ldots,
rS_n)]-g(rS_1,\ldots, rS_n)^* g(rS_1,\ldots, rS_n)\\
&=[I+g(rS_1,\ldots, rS_n)]^* [I+g(rS_1,\ldots, rS_n)]\\
 &\quad -[I+g(rS_1,\ldots, rS_n)]^* f(rS_1,\ldots, rS_n)^* f(rS_1,\ldots, rS_n)
  [I+g(rS_1,\ldots,
rS_n)]\\
&=[I+g(rS_1,\ldots, rS_n)]^*[I-f(rS_1,\ldots, rS_n)^* f(rS_1,\ldots,
rS_n)] [I+g(rS_1,\ldots, rS_n)].
\end{split}
  \end{equation*}
Since $\|f(rS_1,\ldots, rS_n)\|\leq 1$, we deduce that
$h(rS_1,\ldots, rS_n)\geq 0$ for any $r\in[0,1)$, which proves
relation \eqref{g*g}. Therefore, $\cC(f)\in Hol_0^+(B(\cH)^n_1)$.

To prove injectivity of $\cC$, let $f_1, f_2\in
\left[H_0^\infty(B(\cH)^n_1)\right]_{\leq 1}$ such that $\cC f_1=
\cC f_2$. Then
$$
[I-f_1(rS_1,\ldots, rS_n)]^{-1}f_1(rS_1,\ldots, rS_n)=
[I-f_2(rS_1,\ldots, rS_n)]^{-1}f_2(rS_1,\ldots, rS_n).
$$
Multiplying this equality to the left by $I-f_1(rS_1,\ldots, rS_n)$
and  to the right by $I-f_2(rS_1,\ldots, rS_n)$, we deduce that
$f_1(rS_1,\ldots, rS_n)=f_2(rS_1,\ldots, rS_n)$ for $r\in[0,1)$.
Consequently, $f_1= f_2$.

To prove that the noncommutative Cayley transform  is  surjective,
let $g$ be in $Hol_0^+(B(\cH)^n_1)$ and have the representation
$g(X_1,\ldots, X_n):=\sum_{k=1}^\infty \sum_{|\alpha|=k}
B_{(\alpha)}\otimes  X_\alpha$. First, notice that
\begin{equation*}
\begin{split}
H(S_1^{(m)},\ldots, S_n^{(m)})&:=g(S_1^{(m)},\ldots, S_n^{(m)})^* +I
+ g(S_1^{(m)},\ldots, S_n^{(m)})\\
&=[I+ g(S_1^{(m)},\ldots, S_n^{(m)})]^*[I+ g(S_1^{(m)},\ldots,
S_n^{(m)})] -g(S_1^{(m)},\ldots, S_n^{(m)})^* g(S_1^{(m)},\ldots,
S_n^{(m)}).
\end{split}
\end{equation*}
Since $H(S_1^{(m)},\ldots, S_n^{(m)})\geq 0$, we deduce that $ \|[I+
g(S_1^{(m)},\ldots, S_n^{(m)})]x\|\geq \|g(S_1^{(m)},\ldots,
S_n^{(m)})x\| $ for any $x\in\cE\otimes \cP^{(m)}$. Consequently,
there exists a contraction $ A_m:\cE\otimes\cP^{(m)}\to
\cE\otimes\cP^{(m)}$ such that $ A_m[I+ g(S_1^{(m)},\ldots,
S_n^{(m)})]=g(S_1^{(m)},\ldots, S_n^{(m)}). $ Since
$g(S_1^{(m)},\ldots, S_n^{(m)})$ is lower triangular,
$I+g(S_1^{(m)},\ldots, S_n^{(m)})$ is invertible, and therefore
$$
A_m^*=[I+g(S_1^{(m)},\ldots, S_n^{(m)})^*]^{-1}g(S_1^{(m)},\ldots,
S_n^{(m)})^*.
$$

Now, notice that, for each $m\in \NN$ and $i=1,\ldots, n$, we have
$(S_i^{(m+1)})^*|_{\cP^{(m)}}=(S_i^{(m)})^*=S_i^*|_{\cP^{(m)}}. $
Hence,
 $A_{m+1}^*|_{\cE\otimes\cP^{(m)}}=A_m^*$ for any $m\in \NN$.
Using a standard argument, one can prove that there is a unique
contraction  $A\in  B(\cE\otimes F^2(H_n))$ such that
$A^*|_{\cE\otimes\cP^{(m)}}=A_m^*$ for any $m\in \NN$. Indeed, if $x\in
\cE\otimes F^2(H_n)$ let $q_m:=P_{\cE\otimes \cP^{(m)}} x$ and notice
that $\{A_m^* q_m\}_{m=1}^\infty$ is a Cauchy sequence. Therefore,
we can  define $A^*x:=\lim_{m\to \infty} A_m^* q_m$. Since
$\|A_m\|\leq 1$ for $m\in\NN$, so is the operator $A$.

Taking into account that $R_iS_j=S_jR_i$, $i,j=1,\ldots,n$, and
$\cP^{(m)}$, $m\in\NN$, is an invariant subspace under each operator
$R_1,\ldots, R_n$, $S_1,\ldots, S_n$, we deduce that
$(S_j^{(m+1)})^* R_i^*|_{\cP^{(m+1)}}=R_i^* (S_j^{(m+1)})^*. $ Hence
and due to the form of  the operator $A_m$, we have $ (I_\cE\otimes
R_i^*) A_m^*=A^*_m(I_\cE\otimes R_i^*) $
 for any $m\in \NN$ and $i=1,\ldots, n$.
Now, for each $\alpha\in\FF_n^+$ with $|\alpha|=k$, and
$k=0,1,\ldots$, we have
\begin{equation*}
\begin{split}
(I_\cE\otimes R_i^*)A^*(x\otimes e_{\alpha g_i})&= (I_\cE\otimes
R_i^*)A_{k+1}^* e_{\alpha g_i}=A_{k+1}^* (I_\cE\otimes
R_i^*)(x\otimes e_{\alpha g_i})\\
&=A_{k+1}^* (x\otimes e_\alpha) =A_k^*(x\otimes  e_\alpha)
\end{split}
\end{equation*}
and $ A^* (I_\cE\otimes R_i^*)(x\otimes e_{\alpha g_i})=
A^*(x\otimes e_\alpha) =A^*_k (x\otimes e_\alpha)$. Hence, we deduce
that $$(I_\cE\otimes R_i^*)A^*(x\otimes e_{\alpha
g_i})=A^*(I_\cE\otimes  R_i^*)(x\otimes  e_{\alpha g_i}). $$ On the
other hand, if $\alpha\in \FF_n^+$  has the form $g_{i_1}\cdots
g_{i_p}$ with $g_{i_p}\neq i$, then $A^* (I_\cE\otimes
R_i^*)(x\otimes e_{\alpha})=0$ and
$$(I_\cE\otimes R_i^*)A^*(x\otimes
e_{\alpha})=A_{k+1}^*(I_\cE\otimes  R_i^*)(x\otimes e_{\alpha })=0,
$$
which shows that $A^* (I_\cE\otimes R_i^*) (x\otimes
e_{\alpha})=(I_\cE\otimes R_i^*)A^*(x\otimes e_{\alpha})$.
Therefore, $ A(I_\cE\otimes R_i)=(I_\cE\otimes R_i) A$,
$i=1,\ldots,n.$  According to \cite {Po-analytic}, we deduce that
$A$ is in $B(\cE)\bar\otimes F_n^\infty$, the weakly closed algebra
generated by the spatial tensor product. Due to
\cite{Po-holomorphic}, there is a unique $f\in H^\infty(B(\cH)^n_1)$
having the boundary function $A$, i.e., $A=\text{\rm SOT-}
\lim_{r\to 1} f(rS_1,\ldots, rS_n)$. Hence, and using the fact that
$A^*|_{\cE\otimes \cP^{(m)}}=A_m^*$, we deduce that
\begin{equation*}
\begin{split}
A_m=\text{\rm SOT-} \lim_{r\to 1}P_{\cP^{(m)}} f(rS_1,\ldots,
rS_n)|_{\cE\otimes \cP^{(m)}} =\lim_{r\to 1} f(rS_1^{(m)},\ldots,
rS_n^{(m)}) =f(S_1^{(m)},\ldots, S_n^{(m)}).
\end{split}
\end{equation*}
Therefore, we have $f(S_1^{(m)},\ldots, S_n^{(m)})=
 [I+g(S_1^{(m)},\ldots, S_n^{(m)})]^{-1}g(S_1^{(m)},\ldots,
S_n^{(m)}) $ for any $m\in \NN$, which  is equivalent to
$$
f(S_1^{(m)},\ldots, S_n^{(m)})=g(S_1^{(m)},\ldots,
S_n^{(m)})[I-f(S_1^{(m)},\ldots, S_n^{(m)})].
$$
Consequently, $\cC(f)(S_1^{(m)},\ldots,
S_n^{(m)})=g(S_1^{(m)},\ldots, S_n^{(m)}) $ for any $m\in \NN$. By
Lemma \ref{nilp}, we have $\cC(f)=g$, which proves that the Cayley
transform is surjective.
\end{proof}

Denote by $\CC^{(m)}[Z_1,\ldots, Z_n]$, $m\in \NN$, the set
 of all noncommutative polynomials of degree $\leq m$.
 Let $\cA_{n,0}^{(m)}$ be the set of all operators
 $q(S_1^{(m)},\ldots, S_n^{(m)})\in B(\cE\otimes \cP^{(m)})$, where
  $q\in \CC^{(m)}[Z_1,\ldots, Z_n]$ and $q(0)=0$. We also denote by $\cL_{n,0}^{(m)}$
the set of all operators $p(S_1^{(m)},\ldots, S_n^{(m)})\in \cA_{n,0}^{(m)}$ with
 the property that
  $$p(S_1^{(m)},\ldots, S_n^{(m)})^* + I+ p(S_1^{(m)},\ldots, S_n^{(m)})\geq 0.$$
We introduce now the truncated (or constrained)   Cayley transforms
 $\cC^{(m)}$, $m\in \NN$, defined on the unit ball of  the subalgebras
 $\cA_{n,0}^{(m)}$,
 and point out the connection  with the noncommutative Cayley transform.

\begin{theorem}\label{Cayley2}
The Cayley transform $\cC^{(m)}:[\cA_{n,0}^{(m)}]_{\leq 1}\to
\cL_{n,0}^{(m)}$ defined by \ $ \cC^{(m)}(X):= X(I-X)^{-1}$, $ X\in
[\cA_{n,0}^{(m)}]_{\leq 1}, $ is a bijection and its inverse is
given by \ $ [\cC^{(m)}]^{-1} (Y)=Y(I+Y)^{-1}$, $Y\in
\cL_{n,0}^{(m)}. $ Moreover,
\begin{enumerate}
\item[(i)]\
$ \cC^{(m)}[f(S_1^{(m)},\ldots, S_n^{(m)})]=(\cC
f)(S_1^{(m)},\ldots, S_n^{(m)}) $ \ for any $f\in
[H_0^\infty(B(\cH)^n_1)]_{\leq 1}$ and $m\in\NN$;
\item[(ii)]\
$[\cC^{(m)}]^{-1}[g(S_1^{(m)},\ldots, S_n^{(m)})]=[\cC^{-1)}(g)]
(S_1^{(m)},\ldots, S_n^{(m)})$ \ for any $g\in Hol_0^+(B(\cH)^n_1)$.
\end{enumerate}
\end{theorem}
\begin{proof}
First, note that if $X\in [\cA_{n,0}^{(m)}]_{\leq 1}$,
 then $X^{m+1}=0$ and $I-X$ is
 invertible.
Therefore,
$$Y:=X(I-X)^{-1}=X+X^2+\cdots +X^m
$$
has the form $q(S_1^{(m)},\ldots, S_n^{(m)})$, where
  $q\in \CC^{(m)}[Z_1,\ldots, Z_n]$ and $q(0)=0$. Notice also that
  $X(I+Y)=Y$ and
  $$
  Y+I+Y^*=(I+Y)^*(I+Y)-Y^*Y=(I+Y^*)(I-X^*X)(I+Y)\geq 0.
  $$
Therefore $Y\in \cL_{n,0}^{(m)}$. Conversely, if $Y\in
\cL_{n,0}^{(m)}$, then $Y+I+Y^*\geq 0$ and, as in the proof of
Theorem \ref{Cayley1}, there exists a contraction $A_m:\cE\otimes
\cP^{(m)}\to \cE\otimes \cP^{(m)}$ such that $ A_m=Y(I+Y)^{-1}. $
Since $Y^{m+1}=0$, it is easy to see that $A_m$ has the form
$p(S_1^{(m)},\ldots, S_n^{(m)})$ for some polynomial
  $p\in \CC^{(m)}[Z_1,\ldots, Z_n]$  with $p(0)=0$. Hence,  $A_m\in \cA_{n,0}^{(m)}$.
  As in the proof of Proposition \ref{Cayley-series}, one can prove that
   the Cayley transform  $\cC^{(m)}$  is one-to-one and
   $\cC^{(m)}(A_m)=Y$. Therefore $\cC^{(m)}$ is a bijection.

If $f\in [H_0^\infty(B(\cH)^n_1)]_{\leq 1}$ then
$f(S_1^{(m)},\ldots, S_n^{(m)})$ is in $[\cA_{n,0}^{(m)}]_{\leq 1}$.
Since $f(S_1^{(m)},\ldots, S_n^{(m)})^{m+1}= 0$ and $\cC
f=f+f^2+\cdots$, we have
\begin{equation*}
\begin{split}
(\cC f)(S_1^{(m)},\ldots, S_n^{(m)})=f(S_1^{(m)},\ldots,
S_n^{(m)})[I-f(S_1^{(m)},\ldots, S_n^{(m)})]^{-1}
=\cC^{(m)}[f(S_1^{(m)},\ldots, S_n^{(m)})].
\end{split}
\end{equation*}
Therefore,  item (i) holds.  Setting $g=\cC f$ in (i) and using
Theorem \ref{Cayley1}, one can easily  deduce (ii).
    The proof is  complete.
\end{proof}

The following result is a consequence of Theorem \ref{Cayley2}.

\begin{corollary} Let \ $\cT_{n, I}^{(m)}\subset B(\cE\otimes \cP^{(m)})$
be the set
 of all positive operators of the form
$$p(S_1^{(m)},\ldots, S_n^{(m)})^* +I+ p(S_1^{(m)},\ldots, S_n^{(m)}),
$$
where $p(S_1^{(m)},\ldots, S_n^{(m)})\in \cL_{n,0}^{(m)}$.
  Then, there is a one-to-one correspondence
$$A\mapsto T:=[\cC^{(m)}(A)]^*+ I+ [\cC^{(m)}(A)]
$$
between $\cA_{n,0}^{(m)}$  and $\cT_{n, I}^{(m)}$.
\end{corollary}

\bigskip

  \section{Carath\' edory interpolation for free holomorphic
  functions with positive real parts}

In this section we  solve the  Carath\' eodory  interpolation
problem for free holomorphic functions on $[B(\cH)^n]_1$ with
positive real parts and show that it is equivalent to  the Carath\'
eodory-Fej\' er interpolation  problem  for multi-analytic operators
\cite{Po-analytic} and  to the Carath\' eodory  interpolation
problem for positive-definite multi-Toeplitz kernels on free
semigroups \cite{Po-structure}. Using the results from
\cite{Po-structure}, we can provide
 a parametrization of all solutions in terms of generalized Schur sequences.

Recall that $\cP^{(m)}$ is the set of all polynomials in $F^2(H_n)$
of degree $\leq m$.
 According to \cite{Po-analytic}, an operator
$A_m\in B(\cE\otimes \cP^{(m)})$ is called multi-analytic if there
exists a sequence of  operators  $\{A_{(\alpha)}\}_{|\alpha|\leq m}$
in $B(\cE)$  such that $A_m$ has the matrix representation
$[A_{\alpha,\beta}]_{|\alpha|\leq m, |\beta|\leq m}$, where
\begin{equation}
\label{A-ab} A_{\alpha,\beta}:=\begin{cases} A_{(\alpha\backslash _l
\beta)}& \text{
if } \alpha \geq_l \beta\\
0& \text{ otherwise}.
\end{cases}
\end{equation}
Moreover, the set of all  multi-analytic operators on $\cE\otimes
\cP^{(m)}$ coincide with the commutant of the operators
$I_\cE\otimes S_1^{(m)},\ldots, I_\cE\otimes S_n^{(m)}$, where
$S_i^{(m)}:=P_{\cP^{(m)}} S_i|_{\cP^{(m)}}$,
 $i=1,\ldots,n$.

  The definition of  a
multi-analytic operator $A\in B(F^2(H_n))$ is now clear. Moreover,
we proved in \cite{Po-analytic} that  $A\in B(\cE\otimes F^2(H_n))$
is a multi-analytic operator if and only $A\in
B(\cE)\bar\otimes\cR_n^\infty$, the weakly closed algebra generated
by the spatial tensor product. In this case, there exists a unique
sequence of operators $\{A_{(\alpha)}\}_{\alpha\in \FF_n^+}$ in
$B(\cE)$ such that $A$ has the Fourier representation \
$\sum_{\alpha\in \FF_n^+} A_{(\alpha)}\otimes  R_\alpha$.

The Carath\' eodory-Fej\' er  interpolation problem  for the
noncommutative analytic Toeplitz algebra $\cR_n^\infty$ is the
following: given $\{A_{(\alpha)}\}_{|\alpha|\leq m}\subset B(\cE)$,
find a sequence $\{A_{(\alpha)}\}_{|\alpha|\geq  m+1}\subset B(\cE)$
such that $\sum_{\alpha\in \FF_n^+} A_{(\alpha)}\otimes  R_\alpha$
is the Fourier representation of an element $f\in
B(\cE)\bar\otimes\cR_n^\infty$ with $\|f\|\leq 1$. This problem was
solved in \cite{Po-analytic} where, using the noncommutative
commutant lifting theorem \cite{Po-isometric}, we proved that the
Carath\' eodory-Fej\' er interpolation problem  for the
  $\cR_n^\infty$ has solution if and only if $\|A_m\|\leq 1$, where
  $A_m$ is defined   above.

\begin{lemma}\label{plu-pro}
Let $u$ be a free pluriharmonic function on $[B(\cH)^n]_1$  with operator-valued
coefficients. Then $u$ is positive  on $[B(\cH)^n]_1$ if and only if
$u(S_1^{(m)},\ldots, S_n^{(m)})\geq 0$ for any $m\in \NN$. If the positive free
 pluriharmonic function has the representation
 $$
u(X_1,\ldots, X_n)=\sum_{k=1}^\infty\sum_{|\alpha|=k}
 A_{(\alpha)}^*\otimes X_\alpha^*+
 A_{(0)}\otimes I+ \sum_{k=1}^\infty\sum_{|\alpha|=k}
 A_{(\alpha)}\otimes X_\alpha,
 $$
then $ \left\|\sum_{|\alpha|=k}
A_{(\alpha)}^*A_{(\alpha)}\right\|^{1/2}\leq \|A_{(0)}\| $ for any $
k\geq 0$.
\end{lemma}
\begin{proof}
Assume that  $u$ has the representation
$$
u(X_1,\ldots, X_n)=\sum_{k=1}^\infty\sum_{|\alpha|=k}
 B_{(\alpha)}\otimes X_\alpha^*+
 A_{(0)}\otimes I+ \sum_{k=1}^\infty\sum_{|\alpha|=k}
 A_{(\alpha)}\otimes X_\alpha
 $$
 and $u(S_1^{(m)},\ldots, S_n^{(m)})\geq 0$ for any $m\in \NN$.
 Since $S_\alpha^{(m)}=0$ for $\alpha\in \FF_n^+$ with  $|\alpha|\geq m+1$, the
  latter inequality is equivalent
 to
 $$
 T_m:=\sum_{1\leq |\alpha|\leq m} B_{(\alpha)}\otimes
(S_\alpha^{(m)})^*+A_{(0)}\otimes I+ \sum_{1\leq |\alpha|\leq m}
A_{(\alpha)}\otimes  S_\alpha^{(m)}\geq 0
$$
for any $m\in \NN$. Since $T_m^*=T_m$, we deduce that
$B_{(\alpha)}=A_{(\alpha)}^*$ for any $\alpha\in \FF_n^+$. As in the
proof of Theorem \ref{H=M=S}, one can show that, for any vector of
the form $\sum_{|\beta|\leq m}  h_\beta \otimes e_\beta \in \cE
\otimes\cP^{(m)}$,
$$
\left< T_m\left(\sum_{|\beta|\leq m}  h_\beta \otimes
e_\beta\right), \sum_{|\gamma|\leq m}  h_\gamma \otimes
e_\gamma\right>= \sum_{|\beta|\leq m, |\gamma|\leq m}\left<
K(\gamma, \beta) h_\beta, h_\gamma\right>,
$$
where the operator matrix
 $[K(\alpha, \beta)]_{|\alpha|\leq m, |\beta|\leq m}$
  is defined by
\begin{equation*}
         K(\alpha, \beta):=
         \begin{cases}
          A^*_{ (\beta\backslash_r \alpha)}
      &\text{ if } \beta>_r\alpha\\
          A_{(0)}  &\text{ if } \alpha=\beta\\
          A_{ (\alpha\backslash_r \beta)}
      &\text{ if } \alpha>_r\beta\\
          0\quad &\text{ otherwise}
         \end{cases}
         \end{equation*}
for any $|\alpha|\leq m, |\beta|\leq m$.
Hence $[K(\alpha, \beta)]_{|\alpha|\leq m, |\beta|\leq m}$
is a positive
  multi-Toeplitz matrix for any $m\in \NN$.
  According to \cite{Po-posi} there exists a completely positive map
  $\mu:C^*(R_1,\ldots, R_n)\to B(\cE)$ such that
  $\mu(R_{\widetilde \alpha})=A_{(\alpha)}^*$ for $\alpha\in \FF_n^+$.
  Therefore,
   $
  u(X_1,\ldots, X_n) =(\cP\mu)(X_1,\ldots, X_n)
  $
Since $\mu$ is completely positive, $\cP\mu\geq 0$.
The converse is obvious.

Since $R_1,\ldots, R_n$ are isometries with orthogonal ranges, so are the
operators $R_\alpha$ if $|\alpha|=k$, where  $k=1,2, \ldots$. Therefore, the row
operator $[R_\alpha:\  |\alpha|=k]$ has norm one. Since $\mu$ is a completely positive
linear map, we have
\begin{equation*}
\begin{split}
\left\|\sum_{|\alpha|=k} A_{(\alpha)}^*A_{(\alpha)}\right\|^{1/2}
&=\|[A_{(\alpha)}^*:\ |\alpha|=k]\| =\|[\mu(R_{\widetilde \alpha}):\
|\alpha|=k]\| \leq \|\mu\|_{cb}\|[R_\alpha:\ |\alpha|=k]\| \leq
\|\mu(I)\|.
\end{split}
\end{equation*}
The proof is complete.
\end{proof}

Using our noncommutative Cayley transform (Theorem \ref{Cayley1}),
we  prove  the following Carath\'eodory interpolation result
 for  free holomorphic functions  with positive real parts on $[B(\cH)^n]_1$
 and coefficients in $B(\cE)$, where $\cE$ is a separable Hilbert space.

\begin{theorem} \label{Cara1}
Let $\{B_{(\alpha)}\}_{|\alpha|\leq m}$ be a  sequence of operators
in $B(\cE)$
 with $B_{(0)}\geq 0$.
Then there exists a sequence $\{B_{(\alpha)}\}_{|\alpha|\geq
m+1}\subset B(\cE)$ such that
\begin{equation}
\label{g-frac} g(X_1,\ldots,X_n):=\frac
{B_{(0)}}{2}+\sum_{k=1}^\infty \sum_{|\alpha|=k} B_{(\alpha)}\otimes
X_\alpha,\qquad X:=(X_1,\ldots, X_n)\in [B(\cH)^n]_1,
\end{equation}
 is a free holomorphic  function with  positive real part,
  i.e., $\text{\rm Re\,} g(X_1,\ldots,X_n)\geq 0$ for any
$(X_1,\ldots,X_n)$  in $[B(\cH)^n]_1$,
 if and only if
\begin{equation}\label{posi-top}
\sum_{1\leq |\alpha|\leq m} B^*_{(\alpha)}\otimes
(S_\alpha^{(m)})^*+B_{(0)}\otimes I+ \sum_{1\leq |\alpha|\leq m}
B_{(\alpha)}\otimes  S_\alpha^{(m)}\geq 0.
\end{equation}
\end{theorem}
\begin{proof}
First, assume that $g\in Hol^+(B(\cH)^n_1)$  and has the
representation \eqref{g-frac}. Applying Lemma \ref{plu-pro} to the
free pluriharmonic function $u:=2\text{\rm Re\,} g$, we deduce
condition \eqref{posi-top}.

Conversely, assume that $\{B_{(\alpha)}\}_{|\alpha|\leq m}$ is
a sequence of  operators in $B(\cE)$ such that \eqref{posi-top}
holds, and denote
\begin{equation}\label{Y-sum}
Y:=\sum_{1\leq |\alpha|\leq m} B_{(\alpha)} \otimes S_\alpha^{(m)}.
\end{equation}
First, we consider the case when $B_{(0)}=I_\cE$.
 According to Theorem \ref{Cayley2}, the inverse truncated
Cayley transform $[\cC^{(m)}]^{-1} (Y)$ is a multi-analytic operator
on $\cE\otimes \cP^{(m)}$, of the form $X:=\sum_{1\leq |\alpha|\leq
m} A_{(\alpha)}\otimes  S_\alpha^{(m)}$ for some operators
$\{A_{(\alpha)}\}_{1\leq |\alpha|\leq m}\subset B(\cE)$. Applying
the Carath\' eodory-Fej\' er interpolation result for multi-analytic
operators \cite{Po-analytic}, we find  a sequence
$\{A_{(\alpha)}\}_{|\alpha|\geq m+1}\subset B(\cE)$ such that
$\sum_{\alpha\in \FF_n^+}A_{(\alpha)}\otimes  R_\alpha$ is the
Fourier representation of a an element $\varphi\in B(\cE)\bar
\otimes \cR_n^\infty$ with $\|\varphi\|\leq 1$. Let $f$ be the free
holomorphic function on $[B(\cH)^n]_1$ with boundary function
$\varphi$, i.e., $ f(X_1,\ldots, X_n)=\sum_{k=1}^\infty
\sum_{|\alpha|=k} A_{(\alpha)}\otimes X_\alpha$, \
$(X_1,\ldots,X_n)\in [B(\cH)^n]_1. $ Since  $f\in
[H_0^\infty(B(\cH)^n_1)]_{\leq 1}$,   Lemma  \ref{nilp} implies
\begin{equation}
\label{X-f} f(S_1^{(m)}, \ldots, S_n^{(m)})=X.
\end{equation}
Now, we can use  Theorem \ref{Cayley1} to deduce that the
noncommutative Cayley transform $\psi:=\cC(f)$ is in $
Hol_0^+(B(\cH)^n_1)$ and has a representation    $ \psi(X_1,\ldots,
X_n)=\sum_{k=1}^\infty \sum_{|\alpha|=k} C_{(\alpha)}\otimes
X_\alpha$. Since $\psi=(1-f)^{-1} f$ and $\cC^{(m)}(X)=Y$, we can
use relations \eqref{X-f}  and \eqref{Y-sum} to obtain
\begin{equation*}
\begin{split}
\sum_{1\leq |\alpha|\leq m} C_{(\alpha)}\otimes
S_\alpha^{(m)}&=\psi(S_1^{(m)}, \ldots, S_n^{(m)}) =(I-f(S_1^{(m)},
\ldots, S_n^{(m)}))^{-1} f(S_1^{(m)}, \ldots,
S_n^{(m)})\\
&=(I-X)^{-1} X=Y =\sum_{1\leq |\alpha|\leq m} B_{(\alpha)}\otimes
S_\alpha^{(m)}.
\end{split}
\end{equation*}
Hence, we deduce that $C_{(\alpha)}=B_{(\alpha)}$ for $1\leq
|\alpha|\leq m$.

Now, we consider the general case when $B_{(0)}\geq 0$. Let $\epsilon>0$ and notice that
condition \eqref{posi-top} implies
$$
\sum_{1\leq |\alpha|\leq m} D_{(\alpha)}(\epsilon)^*\otimes
(S_\alpha^{(m)})^*+D_{(0)}\otimes I+ \sum_{1\leq |\alpha|\leq m}
D_{(\alpha)}(\epsilon)\otimes  S_\alpha^{(m)}\geq 0,
$$
where $D_{(0)}:=I$ and $D_{(\alpha)}(\epsilon)^*:=(B_{(0)}+\epsilon
I)^{-1/2} B_{(\alpha)}(B_{(0)} +\epsilon I)^{-1/2},$  \ $
1\leq|\alpha|\leq m. $ Applying the first part of the proof, we find
a sequence of operators $\{D_{(\alpha)}(\epsilon)\}_{ |\alpha|\geq
m+1}\subset B(\cE)$ such that
$$
\varphi_\epsilon (X_1,\ldots,X_n)=\frac{1}{2} I+\sum_{k=1}^\infty
\sum_{|\alpha|=k} D_{(\alpha)}(\epsilon) \otimes X_\alpha,\quad
(X_1,\ldots,X_n)\in [B(\cH)^n]_1,
$$
is a free holomorphic function on $[B(\cH)^n]_1$ with positive real part.
Define
$$\xi_\epsilon(X_1,\ldots,X_n):=\left[(B_{(0)}+\epsilon I)^{1/2}\otimes I\right]
\varphi_\epsilon(X_1,\ldots,X_n)\left[(B_{(0)}+\epsilon
I)^{1/2}\otimes I\right]
$$
for $(X_1,\ldots,X_n)\in [B(\cH)^n]_1$.
 Notice that $\xi_\epsilon$ is a free holomorphic function on
$[B(\cH)^n]_1$ with positive real part. Moreover, we have
$$
\xi_\epsilon(X_1,\ldots,X_n)=\frac{1}{2}(B_{(0)}+\epsilon I)\otimes
I+\sum_{1\leq |\alpha|\leq m}
 B_{(\alpha)}\otimes
X_\alpha +\sum_{k=m+1}^\infty
\sum_{|\alpha|=k}C_{(\alpha)}(\epsilon)\otimes X_\alpha,
$$
where $C_{(\alpha)}(\epsilon):=(B_{(0)}+\epsilon I)^{1/2}
D_{(\alpha)}(\epsilon) (B_{(0)}+\epsilon I)^{1/2}$ for $|\alpha|\geq
m+1$. By Lemma \ref{plu-pro}, we have
\begin{equation}
\label{CCBB}
\left\|\sum_{|\alpha|=k}C_{(\alpha)}(\epsilon)^*C_{(\alpha)}(\epsilon)\right\|
\leq
\frac{1}{2}\|B_{(0)}+\epsilon I\|.
\end{equation}
Therefore, there exists a constant $M>0$ such that
$\|C_{(\alpha)}(\epsilon)\|\leq M$ for any $\epsilon\in (0,1)$ and
$\alpha\in \FF_n^+$ with $|\alpha|\geq m+1$. Due to Banach-Alaoglu
theorem,  the ball $[B(\cE)]_M^-$  is compact
 in the $w^*$-topology. Since $\cE$ is a separable Hilbert space,
  $[B(\cE)]_M^-$ is  a metric space in the $w^*$-topology which coincides with
   the weak operator topology on $[B(\cE)]_M^-$. Consequently,
 the diagonal process guarantees the
existence  of a sequence $\{\epsilon_m\}$ such that $\epsilon_m\to
0$ and $G_{(\alpha)}:=$WOT-$\lim\limits_{\epsilon_m\to
0}C_{(\alpha)}(\epsilon_m) $ exists if  $|\alpha|\geq m+1$. Due to
\eqref{CCBB}, we have $ \left\|\sum_{|\alpha|=k}
G_{(\alpha)}^*G_{(\alpha)}\right\|^{1/2}\leq \frac{1}{2}\|B_{(0)}\|
$ for any $ k\geq 0$. Consequently, $ \limsup_{k\to
\infty}\left\|\sum_{|\alpha|=k} G_{(\alpha)}^*G_{(\alpha)}
\right\|^{1/2k}\leq 1. $ This implies that  the function defined by
$$
\xi(X_1,\ldots,X_n)=\frac{1}{2}B_{(0)}\otimes I+\sum_{1\leq
|\alpha|\leq m}
 B_{(\alpha)}\otimes
X_\alpha +\sum_{k=m+1}^\infty \sum_{|\alpha|=k}G_{(\alpha)}\otimes
X_\alpha
$$
for $(X_1,\ldots,X_n)\in [B(\cH)^n]_1$
 is free holomorphic. Since Re\,$\xi_\epsilon(X_1,\ldots,X_n)\geq 0$ for
  any $\epsilon\in (0,1)$ and any $(X_1,\ldots,X_n)\in [B(\cH)^n]_1$, we have
  Re\,$\xi_\epsilon(S_1^{(m)},\ldots, S_n^{(m)})\geq 0$ for any $m\in \NN$
   and $\epsilon\in (0,1)$. Hence, taking $\epsilon \to 0$,  we deduce that
   Re\,$\xi(S_1^{(m)},\ldots, S_n^{(m)})\geq 0$ for any $m\in \NN$.
    Applying once again Lemma
   \ref{plu-pro}, we conclude that  $\xi$ has positive real part and complete the proof.
\end{proof}

We recall that a multi-Toeplitz matrix $[K(\sigma,\omega)]_{
|\sigma|\le m, |\omega|\le m}$ with  $K(\sigma,\omega)\in B(\cE)$,
admits a positive-definite multi-Toeplitz extension to $\FF_n^+$ if
there exist some operators $K(\sigma,\omega)\in B(\cE)$ for $
|\sigma|\ge m+1$ and $|\omega|\ge m +1$,  such that $K:\Bbb
F_n^+\times \Bbb F_n^+\to B(\cE)$
 is a positive semidefinite multi-Toeplitz kernel.

In this setting, the Carath\' eodory interpolation problem  is to
find all positive semidefinite multi-Toeplitz extensions of a
positive multi-Toeplitz matrix $[K(\sigma,\omega)]_{ |\sigma|\le m,
|\omega|\le m}$.
We proved (\cite{Po-structure}) that a  multi-Toeplitz kernel on
$\{\sigma\in \Bbb F_n^+: |\sigma|\le m\}$, admits a
positive semidefinite multi-Toeplitz extension to $\Bbb F_n^+$ if and
only if the operator matrix
 $M_m:=[K(\sigma,\omega)]_{ |\sigma|\le m, |\omega|\le m}$
is positive.

Using   Theorem \ref{Cara1} and Theorem \ref{H=M=S}, we can obtain
another proof of the above-mentioned  result as well  as the
following.

\begin{theorem}\label{Cara-equiv}
 The following problems are equivalent:
\begin{enumerate}
\item[(i)] Carath\' eodory  interpolation
problem for free holomorphic functions on $[B(\cH)^n]_1$ with
positive real parts;
 \item[(ii)]
Carath\' eodory-Fej\' er interpolation  problem  for multi-analytic
operators;
\item[(iii)]
 Carath\' eodory  interpolation
problem  for positive semidefinite multi-Toeplitz kernels on free
semigroups.
\end{enumerate}
\end{theorem}

\begin{proof}
To prove the implication (i)$\implies$(ii), let
$\{A_{(\alpha)}\}_{|\alpha|\leq m}\subset B(\cE)$ be such that
$\|A_m\|\leq 1$, where $A_m$ has the matrix representation
 $[A_{\alpha, \beta}]_{|\alpha|\leq m, |\beta|\leq m}$, given by  \eqref{A-ab}.
 Notice that $B:=\sum_{|\alpha|\leq m} A_{(\alpha)}\otimes S_{g_1 \alpha}^{(m+1)}$ is
  in $[\cA_{n,0}^{(m+1}]_{\leq 1}$ and
$g(S_1^{(m+1)},\ldots, S_n^{(m+1)}):= \cC^{(m+1)}(B)$ is in
$\cL_{n,0}^{(m+1)}$. Therefore,
 $$g(S_1^{(m+1)},\ldots, S_n^{(m+1)})
=\sum_{1\leq |\sigma|\leq m+1} B_{(\sigma)}\otimes
S_{\sigma}^{(m+1)}
$$
 for
some operators  $\{B_{(\sigma)}\}_{1\leq |\sigma|\leq m+1}$, and $
g(S_1^{(m+1)},\ldots, S_n^{(m+1)})^* +I +g(S_1^{(m+1)},\ldots,
S_n^{(m+1)})\geq 0. $ Since (i) holds, we find a sequence of
operators $\{B_{(\sigma)}\}_{|\sigma|\geq m+2}$ such that the
function
  $g(X_1,\ldots,X_n):= \sum_{|\sigma|\geq 1} B_{(\sigma)}\otimes  X_\sigma$
  is  free holomorphic  and
  $g(X_1,\ldots,X_n)^* + I+g(X_1,\ldots,X_n)\geq 0$
    for   \ $(X_1,\ldots,X_n)\in [B(\cH)^n]_1.
   $
   Due to Theorem \ref{Cayley1}, the function $f:=\cC^{-1} (g)$ is in
   $[H_0^\infty(B(\cH)^n_1)]_{\leq 1}$  and has the form
   $f(X_1,\ldots,X_n)=\sum_{\beta|\geq 1} C_{(\beta)}\otimes  X_\beta$.
   On the other hand,
    using Theorem \ref{Cayley2}, we deduce that
    \begin{equation*}
    \begin{split}
f(S_1^{(m+1)},\ldots, S_n^{(m+1)})&=[\cC^{-1}(g)](S_1^{(m+1)},\ldots,
 S_n^{(m+1)})
=[\cC^{(m+1)}]^{-1}[g(S_1^{(m+1)},\ldots, S_n^{(m+1)})]\\
& =[\cC^{(m+1)}]^{-1}[g_{m+1}(S_1^{(m+1)},\ldots, S_n^{(m+1)})]
=B=\sum_{|\alpha|\leq m} A_{(\alpha)}\otimes  S_{g_1
\alpha}^{(m+1)}.
    \end{split}
    \end{equation*}
Hence, $ C_{(g_1\alpha)}=A_{(\alpha)}$  and $ C_{(g_i\alpha)}=0$
 for any $|\alpha|\leq m$  and  $i=2,\ldots, n$. Consequently,
 the boundary function  of the free holomorphic function $f$ has
 the Fourier representation
 $
 \sum_{|\alpha|\leq m}
  A_{(\alpha)}\otimes  S_{g_1\alpha} +\sum_{|\beta |\geq m+2} C_{(\beta)}\otimes  S_\beta.
  $
Let $\varphi$ be the free holomorphic function on $[B(\cH)^n]_1$
which has the boundary function $(I_\cE\otimes S_1^*)f$. It is clear
now that $\varphi$ is in $[H^\infty(B(\cH)^n_1)]_{\leq 1}$ and has
the form
$$
\varphi(X_1,\ldots,X_n)=\sum_{|\alpha|\leq m} A_{(\alpha)}\otimes
X_\alpha +\sum_{|\gamma|\geq m+1} C_{(\gamma)}\otimes  X_\gamma,
$$
which implies (ii). The proof of the  implication (ii)$\implies$ (i)
is contained in the proof of
 Theorem \ref{Cara1}. We prove now that (i)$\implies$(iii).
 Let $[K(\alpha, \beta)]_{|\alpha|\leq m, |\beta|\leq m}$ be a positive
  multi-Toeplitz matrix  and  set $B_{(\alpha)}:=K(\alpha, g_0)$, $|\alpha|\leq m$.
  Consequently, we have
\begin{equation*}
         K(\alpha, \beta):=
         \begin{cases}
          B^*_{ (\beta\backslash_r \alpha)}
      &\text{ if } \beta>_r\alpha\\
          B_{(0)}  &\text{ if } \alpha=\beta\\
          B_{ (\alpha\backslash_r \beta)}
      &\text{ if } \alpha>_r\beta\\
          0\quad &\text{ otherwise}
         \end{cases}
         \end{equation*}
         for any $|\alpha|\leq m, |\beta|\leq m$.
As in the proof of Theorem \ref{H=M=S},  one can show that the
matrix $[K(\alpha, \beta)]_{|\alpha|\leq m, |\beta|\leq m}$ is
positive if and only if $T_m\geq 0$.

Applying now Theorem \ref{Cara1}, there is  a sequence of  operators
$\{B_{(\alpha)}\}_{|\alpha|\geq m+1}\subset B(\cE)$ such that
\begin{equation}
\label{g-form}
 g(X_1,\ldots,X_n):=\frac{B_{(0)}}{2}+ \sum_{k=1}^\infty
\sum_{|\alpha|=k} B_{(\alpha)}\otimes  X_\alpha, \quad
(X_1,\ldots,X_n)\in [B(\cH)^n]_1,
\end{equation}
 is a free holomorphic function with positive real part.
Thus $g\in Hol^+(B(\cH)^n_1)$ and, due to Theorem \ref{H=M=S}, $g\in
\cS^+(B(\cH)^n_1)$. Therefore (iii) holds. The converse
(iii)$\implies$(i) is based on similar arguments. Assume that
$\{B_{(\alpha)}\}_{|\alpha|\leq m}$ is a sequence of  operators such
that condition \eqref{posi-top} holds. Then  the matrix $[K(\alpha,
\beta)]_{|\alpha|\leq m, |\beta|\leq m}$ is positive and due to
(iii) it admits a positive semidefinite  multi-Toeplitz extension $K:\FF_n^+\times
\FF_n^+\to B(\cE)$. Applying again Theorem \ref{H=M=S}, we find a
free holomorphic function $g$ of the form \eqref{g-form} with
positive   real part. This completes the proof.
\end{proof}

Using the results of this section together with Theorem 3.1 from \cite{Po-posi},
we deduce the following result.

\begin{remark} The Carath\' eodory  interpolation
problem for free holomorphic functions with
positive real parts on $[B(\cH)^n]_1$  has a solution if and only if
 there is a completely  positive linear  map
 $$\nu: {\cA_n^*+ \cA_n}\to B(\cE)\ \text{  such that
 } \
 \nu(S_\alpha)=B^*_{(\alpha)},\quad |\alpha|\leq m,
 $$
 i.e., $\nu$ solves  the
 noncommutative trigonometric moment problem  for
 the operator system
 $ {\cA_n^*+ \cA_n}$, with data
 $\{B^*_{(\alpha)}\}_{|\alpha|\leq m}$.
 \end{remark}

 We say
that a multi-Toeplitz kernel
$ K:\Bbb F_n^+\times \Bbb F_n^+\to B(\cE) $ has a Naimark dilation
if there is a Hilbert space $\cK\supset\cE$
and an $n$-tuple $(V_{1},\ldots,V_{n})$ of isometries on $\cK$ with
orthogonal ranges such that
$ K(g_0,\sigma)=P_\cE V_\sigma|_\cE$ for any $\sigma\in\FF_n^+. $
The Naimark dilation is called minimal if
$~\cK=\bigvee_{\sigma\in\FF_n^+} V_\sigma\cE$.
The $n$-tuple  $(V_{1},\ldots,V_{n})$ is called the minimal
isometric dilation of $K$. Extending on the classical case
\cite{SzF-book}, we proved in \cite{Po-structure} that a
multi-Toeplitz kernel on $\FF_n^+$ is positive semidefinite if and
only if it admits a minimal Naimark dilation. In this case its
minimal Naimark dilation is unique up to an isomorphism.

In \cite{Po-structure}, we showed  that any  positive semidefinite
multi-Toeplitz kernel $ K:\Bbb F_n^+\times \Bbb F_n^+\to B(\cE) $
uniquely determines and is uniquely determined by a sequence of row
contractions $\{\Gamma_j\}_{j=1}^\infty$ called generalized Schur
sequence.
We also  obtained   a concrete matrix representation of the minimal
Naimark dilation   for
 positive semidefinite  multi-Toeplitz kernels  on free
 semigroups, in terms of their generalized Schur sequences,
 extending the noncommutative minimal isometric  dilation theorem
 for row contractions.
This  geometric version of the minimal Naimark dilation was used to
show that there is a one-to-one correspondence between the set of
all positive multi-Toeplitz matrices $M_m:=[K(\sigma,\omega)]_{
|\sigma|\leq m, |\omega|\leq m}$ and the Schur sequences
$\{\Gamma_j\}_{j=1}^m$, and recursively calculate
$\{\Gamma_j\}_{j=1}^m$ from $\{K(g_0,\sigma)\}_{|\sigma|\leq m}$.
We  also obtained a parametrization of all solutions  of the
Carath\' eodory interpolation problem for   positive semidefinite
multi-Toeplitz kernels in terms of generalized Schur sequences.
Consequently,  using  the results of this section,  we  have now  a
parametrization of all solutions of any of the Carath\' eodory type
interpolation problems of Theorem \ref{Cara-equiv}.
%
%

      %\Refs
      %\widestnumber\key{BFPQR}
      %\def\n{\key}
       %


\begin{thebibliography}{99}



\bibitem{Ar} {\sc W.B.~Arveson},
 Subalgebras of $C^*$-algebras,
 {\it Acta Math.}
 {\bf 123} (1969),  141--224.











\bibitem{Be} {\sc F.A.~Berezin},
Covariant and contravariant symbols of operators, (Russian), {\it
Izv. Akad. Nauk. SSSR Ser. Mat.} {\bf 36} (1972), 1134--1167.










\bibitem{B} {\sc H.~Bohr},
      {A theorem concerning power series,}
      {\it Proc. London Math. Soc.} (2) {\bf 13}  (1914), 1--5.







\bibitem{Bu} {\sc J.~W.~Bunce},
 Models for n-tuples of noncommuting operators, {\it J. Funct. Anal.}
  {\bf 57}(1984), 21--30.


\bibitem{Ca} {\sc C.~Carath\' eodory},
 \" Uber den Variabilit\" atsbereich der Koeffzienten von
  Potenzreinen die gegebene Werte nicht annehmen,
{\it Math. Ann.}
 {\bf 64}  (1907), 95--115.


\bibitem{CaFe} {\sc C.~Carath\' eodory and L.~Fej\' er},
\" Uber den  Zusammenhang der Extremen von harmonischen  Funktionen
mit ihren Koeffizienten und \" uber den Picard-Landau'chen Satz,
{\it Rend. Circ. Mat. Palermo} {\bf 32} (1911),  218--239.

\bibitem{Co} {\sc J.B.~Conway},
{\em Functions of one complex variable. I.} Second Edition. Graduate Texts in Mathematics  {\bf 159}.
{ Springer-Verlag, New York}, 1995.




\bibitem{Cu} {\sc J.~Cuntz},
 Simple $C^*$--algebras generated by isometries,
 {\it  Commun. Math. Phys.}
 {\bf 57} (1977), 173--185.




\bibitem{DKP}  {\sc K.~R.~Davidson, E.~Katsoulis, and D.R.~Pitts},
  The structure of free semigroup algebras,
 {\it J. Reine Angew. Math.}
  {\bf 533} (2001), 99--125.







\bibitem{DLP} {\sc K.~R.~Davidson, J.~Li, and D.R.~Pitts},
  Absolutely
continuous representations and a Kaplansky density theorem for free
semigroup algebras, {\it J. Funct. Anal.} {\bf 224} (2005), no. 1,
160--191.




\bibitem{DP2} {\sc K.~R.~Davidson and D.R.~Pitts},
%
 The algebraic structure of  non-commutative
  analytic Toeplitz algebras,
{\it  Math. Ann.}
   {\bf 311} (1998),  275--303.



\bibitem{DP1}  {\sc K.~R.~Davidson and D.R.~Pitts},
 Invariant subspaces and hyper-reflexivity for free semigroup algebras,
{\it  Proc. London Math. Soc.}
 {\bf 78} (1999),  401--430.







\bibitem{ER} {\sc E.G.~Effros and Z.J.~Ruan},
  {\em Operator spaces},
 London Mathematical Society Monographs. New Series, {\bf 23}.
 The Clarendon Press, Oxford University Press, New York, 2000.





\bibitem{ES}  {\sc E.v.~Egerv\' ary and O.~Sz\'asz},
{Einige Extremalprobleme im Bereiche der trigonometrischen
Polynome,}
 {\it Math. Zeitschrift} \,{\bf 27} (1928), 641--652.

\bibitem{Fe} {\sc L.~Fej\' er},
{\" Uber trigonometrische Polynome,} {\it J. Reine Angew. Math.}
{\bf 146} (1916), 53--82.




 \bibitem{F} {\sc  A.~E.~Frazho},
 Models for noncommuting operators, {\it J. Funct. Anal.}
  {\bf 48} (1982), 1--11.



\bibitem{Her} {\sc G.~Herglotz},
\" Uber Potenzreien mit positiven, reelen Teil im Einheitkreis,
Berichte \" uber die Verhaundlungen der k\" oniglich s\" achsischen
Gesellschaft der Wissenschaften zu Leipzig, {\it Math.-Phys. Klasse}
{\bf 63} (1911), 501--511.



\bibitem{H} {\sc K.~Hoffman},  {\em Banach Spaces of Analytic Functions},
Englewood Cliffs: Prentice-Hall, 1962.






\bibitem{N} {\sc M.A.~Naimark}, Positive-definite functions on a
commutative group, {\it Bulletin Acad. Sci. URSS} {\bf 7} (1943),
237--244.




\bibitem{P}  {\sc V.I.~Paulsen}, \emph {Completely Bounded Maps and Dilations},
Pitman Research Notes in Mathematics, Vol.146, New York, 1986.






      \bibitem{PPoS}  {\sc V.I.~Paulsen, G.~Popescu,  and D.~Singh},
         On Bohr's inequality,
        {\it Proc. London Math. Soc.}
        {\bf 85} (2002), 493--512.




\bibitem{Pi} {\sc G.~Pisier}, \emph{ Similarity Problems and Completely Bounded Maps},
Springer Lect. Notes Math., Vol.1618, Springer-Verlag, New York,
1995.


\bibitem{Po-isometric} {\sc G.~Popescu}, Isometric dilations for infinite
sequences of noncommuting operators, {\it Trans. Amer. Math. Soc.}
{\bf 316} (1989), 523--536.


\bibitem{Po-multi} {\sc G.~Popescu},
 Multi-analytic operators and some  factorization theorems,
 {\it Indiana Univ. Math.~J.}
 {\bf 38} (1989),   693--710.



      \bibitem{Po-von} {\sc G.~Popescu},
{Von Neumann inequality for $(B(H)^n)_1$,}
      {\it Math.  Scand.} {\bf 68} (1991), 292--304.
      %%\bigskip






   %\leavevmode\vrule height 2pt depth -1.6pt width 44pt,



      \bibitem{Po-funct} {\sc G.~Popescu},
      {Functional calculus for noncommuting operators,}
       {\it Michigan Math. J.} {\bf 42} (1995), 345--356.



      \bibitem{Po-analytic} {\sc G.~Popescu},
      {Multi-analytic operators on Fock spaces,}
      {\it Math. Ann.} {\bf 303} (1995), 31--46.



      \bibitem{Po-disc}  {\sc G.~Popescu},
      {Noncommutative disc algebras and their representations,}
       {\it Proc. Amer. Math. Soc.} {\bf 124} (1996),  2137--2148.


     \bibitem{Po-posi}   G.~Popescu,
      {Positive-definite functions on free semigroups,}
      {\it Canad. J. Math.} {\bf 48} (1996), no. 4, 887--896.





      \bibitem{Po-poisson} {\sc G.~Popescu},
     {Poisson transforms on some $C^*$-algebras generated by isometries,}
       {\it J. Funct. Anal.} {\bf 161} (1999),  27--61.



\bibitem{Po-structure} {\sc G.~Popescu},
Structure and entropy for positive definite Toeplitz kernels on free
semigroups,
   {\it  J. Math. Anal. Appl.} {\bf 254} (2001), 191-218.

\bibitem{Po-curvature} {\sc  G.~Popescu},
  Curvature invariant for Hilbert modules over free semigroup algebras,
   {\it Adv. Math.}
 {\bf 158} (2001), 264--309.




\bibitem{Po-nehari} {\sc G.~Popescu},
  Multivariable Nehari problem and interpolation,
   {\it  J. Funct. Anal.}
   {\bf 200} (2003), 536--581.


\bibitem{Po-entropy} {\sc G.~Popescu},
Entropy and Multivariable Interpolation, {\it Mem. Amer. Math. Soc.}
{\bf 184} (868) (2006).



      \bibitem{Po-holomorphic} {\sc G.~Popescu},
      {Free holomorphic functions on the unit ball of $B(\cH)^n$},
      {\it J. Funct. Anal.}  {\bf 241} (2006), 268--333.


   \bibitem{Po-Bohr} {\sc G.~Popescu},
     Multivariable Bohr inequalities,
   {\it Trans. Amer.  Math. Soc.}
   {\bf 359} (2007), 5283--5317.


      \bibitem{Po-unitary} {\sc G.~Popescu},
      {Unitary invariants in multivariable operator theory}, {\it Mem. Amer. Math. Soc.}, to appear.

\bibitem{Ri}  {\sc F. Riesz},
Sur certains syst\` emes singuliers d'\' equations int\' egrales,
{\it Ann. Sci. Ecole Norm. Sup. (Paris)} {\bf  28}  (1911), 33-62.



 \bibitem{R} {\sc W.~Rudin},
{\em Real and Complex Analysis}, { McGraw-Hill Book Co.} (1966).

\bibitem{R2} {\sc W.~Rudin},
{\em Function theory in the unit ball of \,$\CC^n$}, {
Springer-verlag, New-York/Berlin}, 1980.

\bibitem{S} {\sc D.~Sarason},
 Generalized interpolation in $H^\infty$,
{\it Trans. Amer. Math. Soc.}
  {\bf 127} (1967),   179--203.

\bibitem{Sc}  {\sc I.~Schur},
%
 \"Uber Potenzreihen die im innern des Einheitshreises beschr\"ankt sind,
{\it  J. Reine Angew. Math.}
 {\bf 148} (1918), 122--145.


\bibitem{St} {\sc W.F.~Stinespring},
      {Positive functions on $C^*$-algebras,}
{\it Proc. Amer. Math. Soc.} {\bf 6}  (1955), 211-216.



\bibitem{SzF-book} {\sc B.~Sz.-Nagy and C.~Foia\c{s}}, {\em Harmonic
Analysis of Operators on Hilbert Space}, North Holland, New York
1970.



        %

      \bibitem{vN}  {\sc J.~von Neumann},
      {Eine Spectraltheorie f\"ur allgemeine Operatoren eines unit\"aren
      Raumes,}
      {\it Math. Nachr.} {\bf 4} (1951), 258--281.
      %%\bigskip


      \bibitem{W1}  {\sc G.~Wittstock},
      Ein operatorwertiger Hahn-Banach Satz,
      {\it J. Funct. Anal.}
      {\bf 40} (1981), 127--150.


      \bibitem{W2}  {\sc G.~Wittstock},
      On matrix order and convexity, in Functional analysis: Surveys
      and Recent Results,
      {\it Math. Studies},
      Vol {\bf 90}, North-Holland, Amsterdam, 1984, 175--188.

       \end{thebibliography}
      \end{document}